\documentclass[1p]{elsarticle}
\usepackage{epsf}  
\usepackage{amsfonts, amssymb, amsthm, amsmath}  
\usepackage{pdfsync}
\usepackage{graphicx}
\usepackage{psfrag}
\usepackage{hyperref}
\newtheorem{thm}{Theorem}  

\newtheorem{cor}[thm]{Corollary}  
\newtheorem{lemma}[thm]{Lemma}  
\newtheorem{remark}[thm]{Remark}  
\newtheorem{defn}[thm]{Definition}

\newtheorem{example}[thm]{Example}  
\numberwithin{thm}{section}  
\def\pf{\noindent\emph{Proof: }}  
\def\stop{\hfill$\square$}  
\providecommand{\C}[2]{\ensuremath {C^{#1,\underline{#2}}}}
\providecommand{\St}[1]{\mathcal S\lrb{#1}}
\def\smod{\mathcal M^{\infty,\underline{1}}}

\providecommand{\totl}[1]{\ensuremath{\lceil #1\rceil }}
\providecommand{\totb}[1]{\ensuremath{\underline {#1}}}
\newcommand{\fun}{\totl}

\newcommand{\ex}{\bold}
\providecommand {\e}[1]{\mathfrak t^{#1}}
\providecommand{\tC}[1]{\ensuremath \mathcal E^\times\left(#1\right)}
\providecommand{\tCp}[1]{\ensuremath{ {}^{+}\mathcal E^{\times}\left(#1 \right)}}

\providecommand{\fp}[2]{{}_{\hspace{3pt}#1\hspace{-2pt}}\times_{#2}}

\DeclareMathOperator{\spec}{Spec}
\DeclareMathOperator{\id}{id}
\DeclareMathOperator{\expl}{Expl}

\providecommand{\et}[2]{\ensuremath{\ex T^{#1}_{#2}}}
\providecommand{\lrb}[1]{\ensuremath{\left(#1\right)}}
\providecommand{\abs}[1]{\left\lvert #1\right\rvert}  
\providecommand{\norm}[1]{\left\lVert #1\right\rVert}

\title{Exploded Manifolds}  
\author{Brett Parker  }  
 \ead{ brettdparker@gmail.com}

\begin{document}

\begin{abstract}
This paper provides an introduction to exploded manifolds.
The category of exploded  manifolds is an extension of the category of smooth manifolds with an excellent holomorphic curve theory. Each exploded manifold has a tropical part which is a union of convex polytopes glued along faces. Exploded manifolds are useful for defining and computing Gromov-Witten invariants relative to normal crossing divisors, and using tropical curve counts to compute Gromov-Witten invariants. 

\end{abstract}

\maketitle

\tableofcontents

\section{Introduction}

The category of exploded manifolds is useful for 
\begin{itemize}
\item defining and computing Gromov-Witten invariants relative to normal crossing divisors,
\item computing Gromov-Witten invariants using normal crossing degenerations and degenerations appearing in tropical geometry such as Mikhalkin's higher dimensional pair of pants decomposition of projective hypersurfaces \cite{mpants}.
\item relating `tropical' curve counts to counts of holomorphic curves in a significantly more general setting than toric manifolds.
\end{itemize}

This introductory section sketches how exploded manifolds are useful for studying holomorphic curves. The rest of
this paper contains an introduction to exploded manifolds, including definitions and differential geometric properties. DeRham cohomology of exploded manifolds is discussed in \cite{dre} and the papers \cite{cem}, \cite{reg} and \cite{egw} go on to construct Gromov-Witten invariants of exploded manifolds.

\

The word `tropical' is in quotes above, because some objects which are called tropical in this paper are much more general than what is traditionally studied in tropical geometry. (For example, our tropical curves will satisfy a balancing condition only when certain topological conditions hold.) Each exploded manifold $\ex B$ has a tropical part $\totb{\ex B}$, which can be thought of as a collection of convex polytopes glued over faces using integral affine maps.  The operation of taking the tropical part of an exploded manifold is functorial. (This contrasts to the common procedure of taking the tropicalization of a subvariety of a torus defined over Puiseux series, which is only functorial under monomial maps). 

A second important functor is called the explosion functor. Given a complex manifold $M$ with normal crossing divisors, the explosion of $M$ is an exploded manifold $\expl (M)$. For example, if $M$ is a toric manifold with its toric divisors, the tropical part of $\expl M$ is the toric fan of $M$. The tropical part of any holomorphic curve in $\expl M$ is a tropical curve in $\expl M$. More generally, if $M$ is a complex manifold with a collection $N_{i}$ of transversely intersecting complex codimension $1$ submanifolds, then the tropical part of $\expl(M)$ has one vertex for each connected component of $M$, a copy of $[0,\infty)$ for each submanifold, a face $[0,\infty)^{2}$ for each intersection and an $n$ dimensional quadrant $[0,\infty)^{n}$ for each $n$-fold intersection. (This is sometimes called the dual intersection complex of $M$.)

\includegraphics{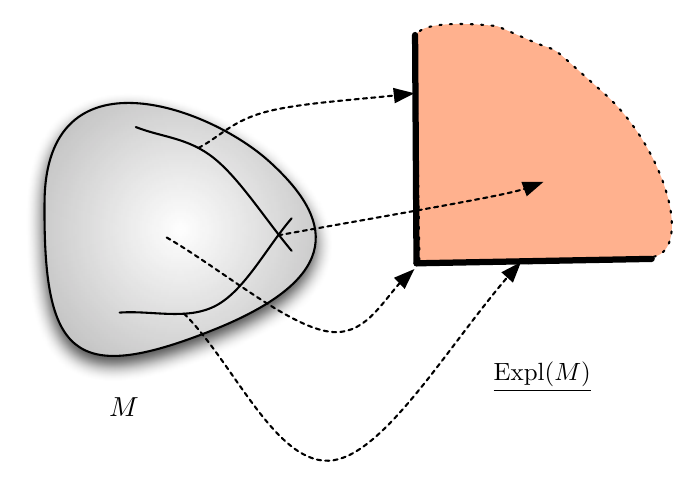}

 This explosion functor can be viewed as a base change in the language of log geometry, and log geometry provides one language for describing some holomorphic exploded manifolds. The paper \cite{elc} explores the link between exploded Gromov-Witten invariants and  log Gromov-Witten invariants defined by Gross and Siebert in \cite{GSlogGW}, and separately by Abramovich and Chen in \cite{acgw}, \cite{Chen}, \cite{Chen2} and \cite{acev}. Defining log Gromov-Witten invariants may be regarded as a step in Gross and Siebert's program on mirror symmetry using tropical geometry. (See for example \cite{gross}.)    Inspired by the Strominger-Yau-Zaslow  conjecture \cite{SYZ}, many other researchers have been exploring the link between tropical geometry and mirror symmetry, such as   Kontsevich and Soibelman in  \cite{kontsevich} and Fukaya in \cite{Fukayagraph}.
 
\

Exploded manifolds have an excellent holomorphic curve theory. As in the smooth setting, holomorphic curves refer to both abstract holomorphic curves (which in our case are one complex dimensional complex exploded manifolds), and holomorphic  maps of abstract holomorphic curves to exploded manifolds with an (almost) complex structure. An abstract exploded curve has a smooth part which corresponds to a nodal Riemann  surface  $\Sigma$ with punctures, and a tropical part which is a metric graph with one vertex for each component, one internal edge for each node, and one semi infinite edge for each puncture.  It is called stable if the nodal Riemann surface $\Sigma$ is stable.

\includegraphics{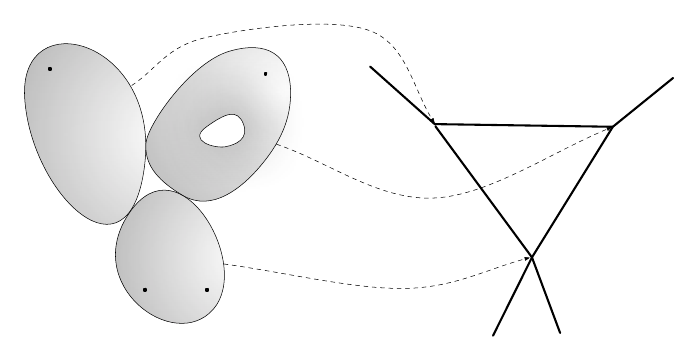}

 It turns out that the moduli stack of stable exploded curves with genus $g$ and $n$ punctures is equal to the explosion of Deligne Mumford space $\expl (\bar{\mathcal M}_{g,n})$, where $\bar{\mathcal M}_{g,n}$ is viewed as a complex orbifold with normal crossing divisors given by the boundary, which consists of the nodal Riemann surfaces. The explosion of the forgetful map $\expl (\bar{\mathcal M}_{g,n+1})\longrightarrow \expl (\bar{\mathcal M}_{g,n})$ is actually a smooth family in the exploded category which has as the fiber over each point in $\expl (\bar{\mathcal M}_{g,n})$ the corresponding exploded curve quotiented by its automorphisms.

The paper \cite{egw} proves that the moduli stack of holomorphic curves in an exploded manifold also has an exploded structure. In nice cases, the moduli space of holomorphic curves is itself an exploded manifold. In the case that holomorphic curves in an exploded manifold are suitably tamed by a symplectic form, the moduli space of holomorphic curves is compact, and has a virtual fundamental class, allowing Gromov-Witten invariants to be defined in the exploded setting. For example, to define Gromov-Witten invariants of a K\"ahler manifold $X$ relative to normal crossing divisors, take the Gromov-Witten invariants of the corresponding exploded manifold $\expl (X)$. As outlined in \cite{elc}, it is expected that the Gromov-Witten invariants of $\expl(X)$ are the same as the log Gromov-Witten invariants of $X$ defined by Gross and Siebert or Abramovich and Chen. Sometimes, these relative Gromov-Witten invariants can be determined by considering only the tropical curves in the tropical part of $\expl X$. For instance, in the case that $X$ is a toric manifold, and we use the toric boundary divisors, work of Mikhalkin, in \cite{Mikhalkin} and  of Siebert and Nishinou in \cite{Sieberttrop} translates to saying that genus zero relative Gromov-Witten invariants can be computed by considering tropical curves in $\totb{\expl X}$. 

A similar (but less functorial) construction can associate an exploded manifold to a symplectic manifold with orthogonally intersecting codimension $2$ symplectic submanifolds, allowing Gromov-Witten invariants relative to these symplectic submanifolds to be defined. The approach of Ionel to relative Gromov-Witten invariants from \cite{IonelGW} should give similar invariants in this case.  

As an example application of these relative Gromov-Witten invariants, suppose that there existed a connected compact symplectic $4$-manifold containing $3$ orthogonally intersecting embedded symplectic spheres, which intersect each other once and have topological self intersection numbers $1$, $1$ and $2$. Then we could construct an exploded manifold $\ex X$ for which the tropical part of holomorphic curves behave like the tropical curve in the following picture:

\includegraphics{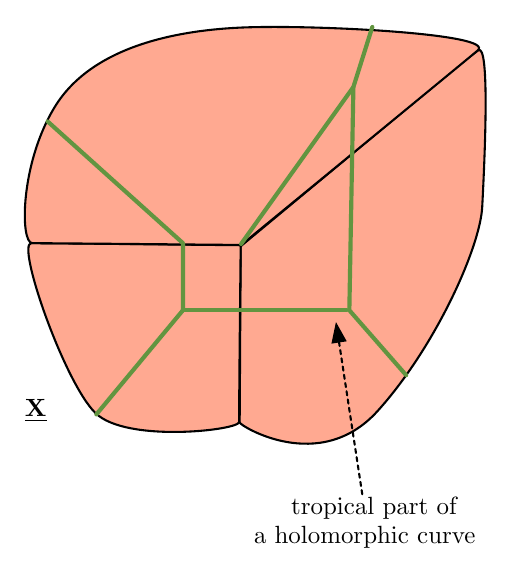} 

In particular, away from the central vertex, the tropical part of holomorphic curves in $\totb{\ex X}$  obey a balancing condition, although when they pass through the left pointing edge in the above picture of $\totb{\ex X}$, they are deflected to the left. 
(In general, the tropical part of a holomorphic curve in an exploded manifold $\ex B$ is always a piecewise linear map of a metric graph into the tropical part of $\ex B$, but it does not always have to satisfy the balancing condition familiar to tropical geometers. For example, let $M$ be a $\mathbb CP^{1}$ bundle over $\mathbb CP^{1}$ with nontrivial Chern class and  divisor given by the zero section. Then the tropical part of $\expl M$ is $[0,\infty)$, and the tropical part of holomorphic curves in $\expl M$ need not satisfy the usual balancing condition.)

  In our hypothetical exploded manifold $\ex X$, consider the following enumerative problem: count the (virtual) number of rigid holomorphic curves that 
\begin{enumerate}
\item 
pass through a given point in $\ex X$
\item have tropical parts which have two seminfinite ends $[0,\infty)$, one of which is mapped in the direction $(0,1)$ in the above picture, and one of which is mapped in the direction $(0,-1)$ in the above picture. 
\end{enumerate}

The paper \cite{egw} proves that this enumerative problem should have as an answer a rational number which does not depend on the point through which our holomorphic curves are required to pass.  If we choose our point to have tropical part in the left hand half of $\totb{\ex X}$, the balancing condition for tropical curves in $\ex X$ implies that there are no holomorphic curve passing through our given point. On the other hand, if we choose our point to have tropical part in the right hand half of $\totb{\ex X}$, the only possible tropical curve satisfying our conditions is drawn below, and a simple calculation implies that there is exactly one holomorphic curve satisfying our conditions passing through this point. The conclusion to be drawn from this contradictory situation is that there does not exist a compact symplectic manifold with embedded symplectic spheres as described above.

\includegraphics{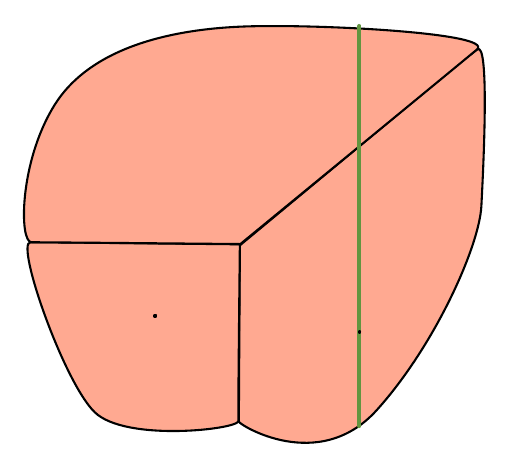}

 The reader may ask how these Gromov-Witten invariants of exploded manifolds correspond to Gromov-Witten invariants of smooth manifolds. Each smooth manifold can also be considered as an exploded manifold, so instead of considering holomorphic curves in a smooth symplectic manifold $X$, we can consider the exploded holomorphic curves. Each exploded manifold $\ex B$ has a `smooth part' $\totl{\ex B}$ which is a topological space. The smooth part of an exploded holomorphic curve in $X$ is a holomorphic curve in $X$ considered as a smooth manifold. The smooth part of the moduli space of exploded curves in $X$ is the usual moduli space of holomorphic curves in $X$.  The moduli space of exploded curves in $X$ therefore contains all the information present in the usual moduli space of  holomorphic curves in $X$, so we may as well work in the exploded category for computing Gromov-Witten invariants. This has the advantage that many degenerations which look nasty from the smooth perspective turn out to be smooth families in the exploded category.
 
 For example, applying the explosion functor to a normal crossings degeneration gives a smooth family of exploded manifolds. As shown in \cite{egw}, holomorphic curve invariants such as Gromov-Witten invariants do not change in connected families of exploded manifolds, so the Gromov-Witten invariants can be computed in any fiber of the resulting family. In particular, suppose that a normal crossing degeneration has a singular fiber which is a union of pieces $X_{i}$. The Gromov-Witten invariants of a smooth fiber can then be computed using the Gromov-Witten invariants of an exploded manifold $\ex X$ corresponding to this  singular fiber, the computation of which involves tropical curve counts in the tropical part of $\ex X$ coupled with relative Gromov-Witten invariants of the pieces $X_{i}$. If these pieces are simple, then this can simplify the computation of Gromov-Witten invariants considerably. 
 
 As a concrete example, consider the $3$ complex dimensional toric manifold with moment map  the polytope defined by the equations
 \[x_{1}\leq 1 ,\ \  x_{2}\leq  1 ,\ \  -x_{1}-x_{2}\leq 2 ,\ \  x_{3}-x_{1}-x_{2}\leq 0,\ \  x_{3}-x_{2}\leq 0,\ \  x_{3}\leq 0\]
 \includegraphics{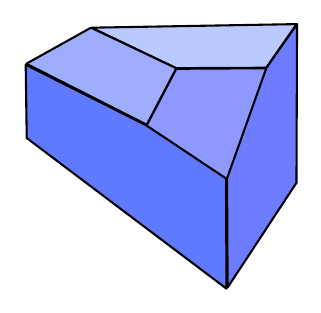}
 The corresponding toric fan is dual to the above polytope, and has edges which are the span of the vectors $(1,0,0)$, $(0,1,0)$, $(-1,-1,0)$, $(-1,-1,1)$, $(0,-1,1)$ and $(0,0,1)$. The corresponding toric manifold can be considered as a partial compactification of $(\mathbb C^{*})^{3}$ where the coordinates $z_{1}, z_{2},$ and $z_{3}$ correspond to projection of the above fan onto the first, second and third components respectively. 
 
  Consider the map $z_{3}$ from our toric manifold to $\mathbb C$. This map $z_{3}$ can be considered as a normal crossing degeneration. The fiber over any point in $\mathbb C-0$ is equal to $\mathbb CP^{2}$, but the fiber over $0$ is singular, and consists of the three toric boundary divisors corresponding to the top faces of the above polytope. (Each is isomorphic to $\mathbb CP^{2}$ blown up at one point.) Consider this toric manifold as a complex manifold with transversely intersecting complex submanifolds given by the toric boundary divisors of the top faces $x_3=0$, $x_{3}=x_{1}+x_{2}$ and $x_{3}=x_{2}$, and call the corresponding exploded manifold $\ex X$. Similarly, we can consider $\mathbb C$ as a complex manifold with the divisor $0$. The explosion of the map $z_{3}$ is a smooth family of exploded manifolds 
  \[\expl z_{3}:\ex X\longrightarrow \expl (\mathbb C,0)\]
 The tropical part of $\ex X$ can be identified with the positive span of $(-1,-1,1)$, $(0,-1,1)$ and $(0,0,1)$, and the tropical part of $\expl z_{3}$ can be identified with projection to the third coordinate.
 
 \includegraphics{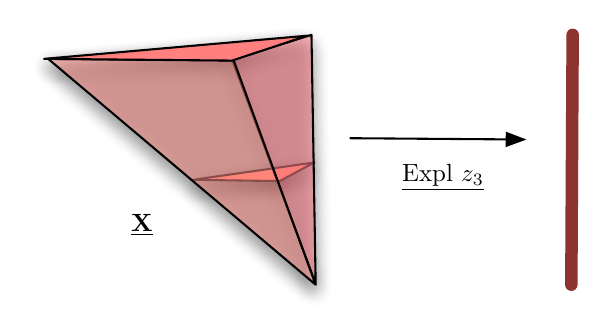}
 
 There are now a large number of exploded manifolds in this family that correspond to the singular fiber of our original degeneration. The following is a picture of the tropical part of one of these exploded manifolds with the tropical part of a holomorphic curve.
 
\includegraphics{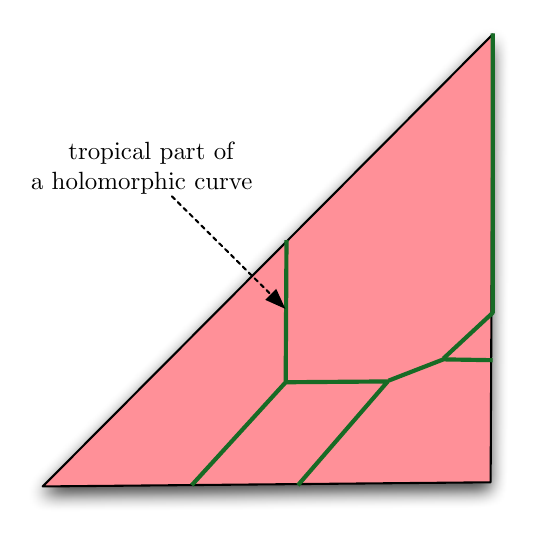}  

In this case, the number of holomorphic curves in the above exploded manifold passing through a given number of points can be calculated using  Mikhalkin's tropical methods from \cite{Mikhalkin2}, by choosing the points to have tropical parts that appear in the interior of the above triangle in `tropical general position'. 

\

 An example of a degeneration `breaking a manifold into simple pieces' is given by Mikhalkin's higher dimensional pair of pants decomposition of projective hypersurfaces from \cite{mpants}. An exploded version of this construction is discussed in example \ref{pair of pants} beginning on page \pageref{pair of pants}. In particular, given any $n$ dimensional projective hypersurface $X$ which intersects the toric boundary of $\mathbb CP^{n+1}$ nicely, one can take the intersection of $X$ with the boundary of $\mathbb CP^{n+1}$ to be a normal crossing divisor of $X$. Then  there is a connected family of exploded manifolds containing both $\expl{X}$ and an exploded manifold $\ex X'$ which has as its tropical part the $n$ dimensional balanced polyhedral complex which appears as the base of Mikhalkin's singular fibration from \cite{mpants}, and which has a smooth part a union of `higer dimensional pairs of pants', which are copies of $\mathbb CP^{n}$ with normal crossing divisors given by $n+2$ hyperplanes. To compute the relative Gromov-Witten invariants of $X$, one can compute the corresponding invariants of $\ex X'$.
 A similar construction gives a connected family of exploded manifolds containing both $X$ and an exploded manifold with tropical part a truncated version of the tropical part of $\ex X'$.

\

For an analytic perspective on degenerations that can be treated nicely in the exploded category, consider the  following example. Think of the following picture as the image of some smooth map from a symplectic manifold to $\mathbb R^{2}$.

\includegraphics{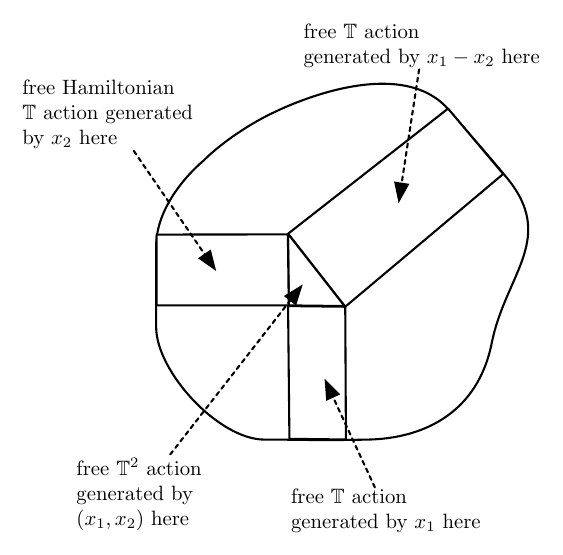}

Assume that in the regions pictured above, the Hamiltonian flow generated by the specified function on $\mathbb R^{2}$ is a free circle action, and that in the central triangle, the two coordinate functions generate commuting free circle actions. From this point on, we shall describe a family of almost complex manifolds, forgetting the original symplectic structure. (The symplectic form is necessary to tame holomorphic curves in the family we shall construct, but we shall not say more about it in this example.) It is possible to choose an almost complex structure $J$ on our symplectic manifold so that the symplectic form is positive on holomorphic tangent planes, and so that $J$ is preserved by the circle actions in the regions where they are defined, and so that $J$ also has  the symmetry  specified in the picture below.

 \includegraphics{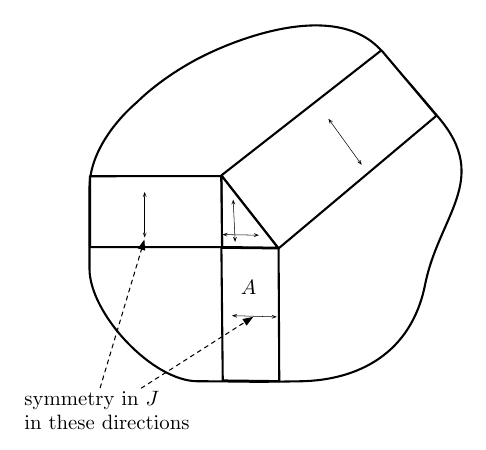}

In particular, it is possible to identify the region labeled $A$ above with $N\times [0,1]$ where $N$ is a manifold with boundary, $x_{1}$ is the coordinate on $[0,1]$, and the circle action is independent of $x_{1}$. In this region, we may choose $J$ to send the vector field generating the circle action to $\frac \partial {\partial x_{1}}$, and require that $J$ is independent of $x_{1}$ and preserved by the circle action. We may specify $J$ similarly in the other regions with circle actions - choosing $J$ independent of both $x_{1}$ and $x_{2}$ in the central triangle region which can be identified with some manifold $M$ times a triangle in $\mathbb R^{2}$ with vertices $(0,0)$, $(1,0)$ and $(0,1)$. 

We are now ready to describe a family of almost complex manifolds for $t\in(0,1]$. We may replace the region $A$ with $N\times [0,\frac 1t]$, extending $J$ symmetrically, replace the other regions with circle actions by similarly stretched regions, and replace the central triangle region with $M$ times a triangle with vertices $(0,0)$, $(\frac 1t,0)$ and $(0,\frac 1t)$, again extending $J$ symmetrically. We may depict this by stretching in all the directions pictured in which $J$ has symmetry.

\includegraphics{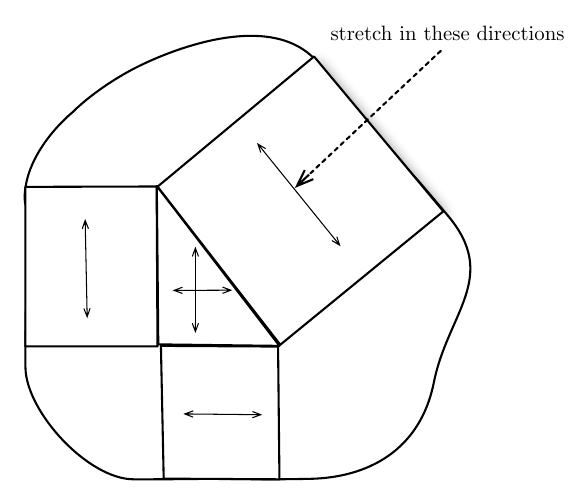} 

The following is a picture of the image of a holomorphic curve in a member of our family with $t\simeq \frac 1 {20}$ . The  picture has been rescaled by a factor of about $5t$. 

\includegraphics{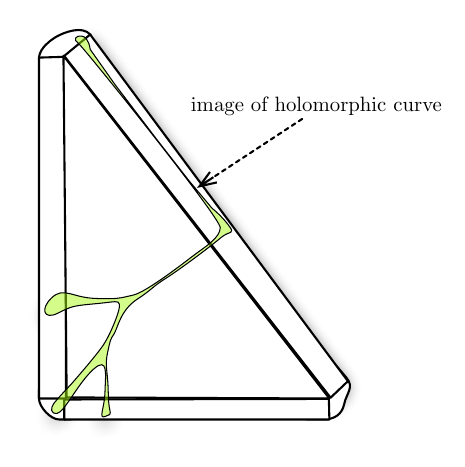}

The important point here is as follows: Consider the image of a holomorphic curve in the picture which has been rescaled by a factor of $t$  for very small $t$. The pieces of a holomorphic curve which look large will approximate a special type of holomorphic cylinder which is preserved by one of our circle actions. The projection of this cylinder to our picture is simply a line with integral slope. In the limit $t\rightarrow 0$, the image of holomorphic curves in the picture rescaled by $t$ will look piecewise linear.  

For extremely small $t$, the image of a holomorphic curve in the picture rescaled by a factor of about $5t$ will look something like this:

\includegraphics{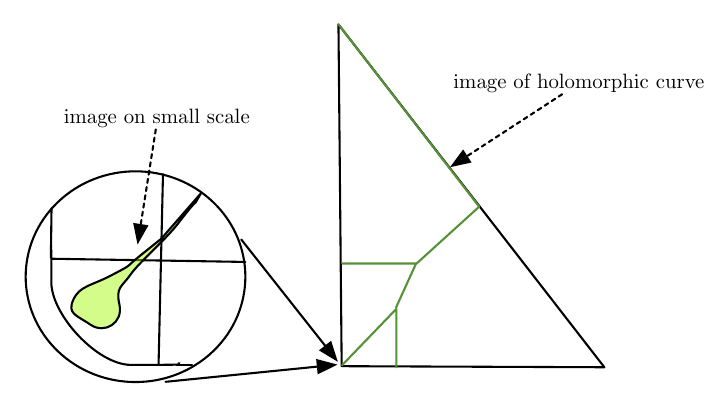}

It is unreasonable to expect that the rescaled image as $t\rightarrow0$ will in general contain all the information necessary to reconstruct holomorphic curve invariants. We must also somehow keep track of information which is on a smaller scale as $t\rightarrow 0$. This is one striking feature that exploded manifolds have. Exploded manifolds have multiple topologies; the tropical part of an exploded manifold could in some sense be thought of as the large scale of the exploded manifold. 

The above family of almost complex manifolds fit into a smooth connected family of exploded manifolds, one of which has a tropical part which looks like the above rescaled picture. (The example containing $\mathbb CP^{2}$ discussed earlier is a concrete example of such a family of exploded manifolds.) This family of exploded manifolds can be thought of as representing a symplectic triple sum. If we think of this family as `stretching over a triangle', then there are analogous families of exploded manifolds `stretching over' any compact convex polytope with integer vertices. 

For example, the usual symplectic sum can be represented by a family of exploded manifolds `stretching over' an interval, and containing an exploded manifold that can be represented pictorially as follows:

\

\includegraphics{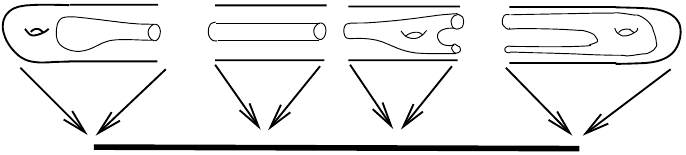}

\

This case of a family representing a symplectic sum has been studied by Eliashberg Givental and Hofer in \cite{sft}, Ionel and Parker in \cite{IP}, Jun Li in \cite{Li}, and Ruan and A. Li in \cite{ruan}. In each of these cases, the researchers did not need knowledge of the large scale to describe the limiting problem of finding holomorphic curves. The relationship between the moduli spaces considered by the above researchers and the exploded moduli space is as follows: Consider the schematic diagram above in which the large scale or tropical part of our exploded manifold is represented as an interval, and the pictures above it represent what is seen of a holomorphic curve in the small scale at different places. All but a finite number of those small scale pieces will have holomorphic curves with translational symmetry in the left right direction, such as depicted in the second picture from the left. The remaining interesting pieces have an order to them. The above researchers record these interesting pieces, along with the order in which they come, and also consider two different pieces of holomorphic curve in the middle to be equivalent if they differ by a  $\mathbb C^{*}$ or  $\mathbb R^{*}$ action (in \cite{IP},\cite{Li}, \cite{ruan}, there is a $\mathbb C^{*}$ action, but in \cite{sft} there is only a $\mathbb R^{*}$ action.) Approaching holomorphic curve from the exploded category can be thought of as a generalization of the above researchers' work to the case of stretching in multiple directions.

 \section{Tropical semiring and exploded semiring}\label{tropical semiring}
 
We shall need the following semirings to describe coordinates on exploded manifolds.

 \ 
  
  The tropical semiring is a semiring $\e {\mathbb R}$ which is equal to $\mathbb R$ with `multiplication' being the operation of usual addition and `addition' being the operation of taking a minimum. Write elements of  $\e {{\mathbb R}}$ as $\e x$ where $x\in\mathbb R$. Then we can write the operations as follows
  \[\e x\e y:=\e{x+y}\]
  \[\e x+\e y:=\e{\min \{x,y\}}\]
  $\e 1$ can be thought of as something which is infinitesimally small, so use the following order on $\e{\mathbb R}$:
  
  \[\e x>\e y\text{ when } x<y\]
 
  Given a ring $R$, the exploded semiring $R\e{ {\mathbb R}}$ consists of elements 
  $c\e x$ with $c\in R$ and $x\in\mathbb R$. Multiplication and addition are as follows:
  \[c_1\e xc_2 \e y=c_1c_2\e{x+y}\]
  \[c_1\e x+c_2\e y=\begin{cases}&c_1\e x\text{ if }x<y
  \\&(c_1+c_2)\e x\text{ if }x=y
  \\&c_2\e y\text{ if }x>y\end{cases}\]
 
  It is easily checked that addition and multiplication are associative and obey the usual distributive rule. The reader familiar with Puiseaux series will  notice that $R\e{\mathbb R}$ may be viewed as arising from  taking the leading term of Puiseaux series  with coefficients in $R$. We will mainly be interested in $\mathbb C\e {\mathbb R}$ and $\mathbb R\e {\mathbb R}$.
 
 There are semiring homomorphisms
 \[ R\xrightarrow{\iota} R\e{\mathbb R}\xrightarrow{\totb{\cdot}}\e{\mathbb R}\]
defined by 
\[\iota(c):=c\e 0\]
\[\totb{c\e x}:=\e x\]
The homomorphism $\totb\cdot:R\e{\mathbb R}\longrightarrow \e{\mathbb R}$ is especially important. We shall call $\totb{c\e x}=\e x$ the tropical part of $c\e x$. There will be an analogous tropical part of  exploded manifolds which can be thought of as the large scale.

Define  $R\e{\mathbb R^+}$ to be the sub semiring of $R\e{\mathbb R}$ consisting of elements of the form $c\e x$ where $x\geq 0$. There is a semiring homomorphism 
\[\totl\cdot:R\e{\mathbb R^+}\longrightarrow R\] 
given by 
\[\totl{c\e x}:=\begin{cases}& c \text{ if }x=0
\\&0\text{ if }x>0\end{cases}\]
Call $\totl{c\e x}$ the smooth part of $c\e x$. Note that this smooth part homomorphism can be thought of as setting $\e 1=0$, which is intuitive when $\e 1$ is thought of as infinitesimally small.

\
 
We shall use the following order on $(0,\infty)\e{\mathbb R}$, which again is intuitive if $\e 1$ is thought of as being infinitesimally small and positive:
\[x_1\e {y_1}<x_2\e{y_2}\text{ whenever }y_1>y_2\text{ or } y_1=y_2\text{ and }x_1<x_2\]

\section{Exploded manifolds}

A smooth manifold can be regarded as a topological space with a sheaf of smooth real valued functions. Similarly, an exploded manifold can be regarded as a (maybe non Hausdorff) topological space with a sheaf of $\mathbb C^{*}\e{\mathbb R}$ valued functions.
The following definition of abstract exploded spaces is far too general, but it allows us to talk about local models for exploded manifolds as abstract  exploded spaces without giving too many definitions beforehand. Think of this as analogous to introducing manifolds by making a definition of an `abstract smooth space' as a topological space with a sheaf of real valued functions, then taking about $\mathbb R^{n}$ with its sheaf of smooth functions as an abstract smooth space, then defining a manifold as an abstract smooth space locally modeled on $\mathbb R^{n}$. 

\begin{defn}[Abstract exploded space]\label{abstract es}

An abstract exploded space  $\ex B$ consists the following:
\begin{enumerate}
\item A (possibly non-Hausdorff) topological space $\ex B$ whose topology is induced from a surjective map to a Hausdorff topological space, $\ex B\longrightarrow \totl{\ex B}$. \footnote{I added this important condition on the topology of $\ex B$ after publishing. A similar condition was present in my original definition of exploded manifolds in \href{http://arxiv.org/abs/0706.3917}{http://arxiv.org/abs/0706.3917}. The published version assumed this was the case for exploded manifolds without it actually being implied by the definition. Sorry. }
\item A sheaf of Abelian groups on this topological space  $\mathcal E^\times(\ex B)$ called the sheaf of exploded functions on $\ex B$ so that:
\begin{enumerate}

\item Each element  $f\in \mathcal E^\times(U)$ is a 
map
\[f: U\longrightarrow  \mathbb C^*\e {\mathbb R}\]
\item  Multiplication is given by pointwise multiplication in $\mathbb C^{*}\e{\mathbb R}$.
\item $\mathcal E^\times(U)$ includes the constant functions if $U\neq\emptyset$.
\item Restriction maps are given by restriction of functions.

\end{enumerate}
\end{enumerate}
\end{defn}

\begin{defn}[Morphism of abstract exploded spaces]
A morphism $f:\ex B\longrightarrow\ex C $ of abstract exploded spaces is a continuous map

\[ f:\ex B\longrightarrow \ex C \]
so that $f$ preserves $\mathcal E^{\times}$ in the sense that if  $g\in\mathcal E^\times(U)$, then  $f\circ g$ is in $\mathcal E^\times\left( f^{-1}(U)\right)$
\[f^*g:=f\circ g\in \mathcal E^\times\left( f^{-1}(U)\right)\] 

\end{defn}

\

 It may occur to the reader that from the perspective of differential geometry, using $(0,\infty)\subset\mathbb R$ instead of $\mathbb C^{*}$ would be a more natural choice for extending the category of smooth manifolds. The choice of $\mathbb C^{*}$ is used so that holomorphic curve theory works out easily in the category of exploded manifolds. Making the same definitions but replacing $\mathbb C^{*}$ with $(0,\infty)$ and $\mathbb C$ with $[0,\infty)$ would give an interesting category worthy of study.
 
 \
 
 The next sequence of examples will give  local models for exploded manifolds.

\

\begin{example}[Smooth Manifold]\label{smooth manifold} \end{example}

Any smooth manifold $M$ can be considered as an abstract exploded space as follows: the toplogical space is just $M$ with the usual topology, and the sheaf $\mathcal E^\times(M)$ consists of all functions of the form $f\e a$ where $f\in C^\infty(M,\mathbb C^*)$ and $a$ is a locally constant $\mathbb R$ valued function.

 Readers should convince themselves that this is just a different way of encoding the usual data of a smooth manifold, and that a morphism between smooth manifolds regarded as abstract exploded spaces is simply a smooth map. In other words,  the category of smooth manifolds a full subcategory of the category abstract exploded spaces.

\

The reader should check that a point considered as an exploded manifold is a final object in the category of abstract exploded spaces: given any abstract exploded space $\ex B$, and a point $p$, there exists a unique map $\ex B\longrightarrow p$.

\

 Example \ref{smooth manifold} should be considered as a `completely smooth' exploded manifold. At the other extreme, we have the following `completely tropical' exploded manifold. 

\begin{example}[$\ex T^{n}$]\end{example}
 The exploded manifold $\ex T^n$ is best described using coordinates $(\tilde z_1,\dotsc, \tilde z_n)$ where each $\tilde z_i\in\mathcal E^\times (\ex T^n)$. The set $\ex T^{n}$ is identified with $\lrb{\mathbb C^*\e {\mathbb R}}^n$ by prescribing the values of $(\tilde z_1,\dotsc,\tilde z_n)$ in $\lrb{\mathbb C^{*}\e{\mathbb R}}^{n}$, and given the trivial topology in which the only open subsets are the empty set and the entire set. We shall regard $\ex T^{n}$ as being $n$ complex dimensional, or $2n$ real dimensional.  
 
  Exploded functions  $f\in \mathcal E^{\times}(\ex T^n)$ can be written in these coordinates as 
  \[f:=c\e y\tilde z^\alpha:=c\e y\prod \tilde z_i^{\alpha_i}\]
   where $c\in\mathbb C^*$, $y\in\mathbb R$ and $\alpha\in \mathbb Z^n$ are all constant. The above function  $f$ takes values as follows: if a point $p\in \ex T^{n}$ has coordinates $(c_{1}\e{x_{1}},\dotsc,c_{n}\e{x_{n}})$, then  $f$ takes the value 
   $f(p)=c\e y\prod c_{i}^{\alpha_{i}}\e{\alpha_{i}x_{i}}$ on $p$.
   
   \
 
 The exploded manifold $\ex T^{n}$ has a `large scale' or  `tropical part' $\totb {\ex T^{n}}$
 which is $\e{\mathbb R^n}$  (given the usual topology  and integral affine structure on $\mathbb R^n$), which can be regarded as giving another, (non Hausdorff) topology on the set of points in $\ex T^{n}$. There is a map $p\mapsto\totb p$ from the set of points in $\ex T^{n}$ to $\totb{\ex T^{n}}$ which is given in coordinates by
 \[\totb{(\tilde z_1,\dotsc,\tilde z_n)}=(\totb{\tilde z_1},\dotsc,\totb{\tilde z_n})\text{ or }\totb{(c_1\e {a_1},\dotsc, c_n\e {a_n})}=(\e{a_1},\dotsc,\e{a_n})\]
 Note that there is a $(\mathbb C^*)^n$ worth of points $p\in \ex T^n$ over every  point $\totb p\in\totb{\ex T^n}$. 
  
  \
  
  Before continuing, the reader should be able to verify the following observations:
  
\begin{enumerate}
   \item A morphism from $\ex T^n$ to a smooth manifold is simply given by a constant map  $\totl{\ex T^n}\longrightarrow M$. 
 
\item A morphism from a connected smooth manifold $M$ to $\ex T^n$  has the information of a map $\totb f$ from $M$ to a point  $\totb p\in\totb{\ex T^n}$ and a smooth map $f$ from $M$ to the $(\mathbb C^*)^n$ worth of points over $\totb p$.  The map $f$ is determined by specifying the $n$ exploded functions \[f^*(\tilde z_i)\in \tC  M\]

\

Note in particular that $\tC  M$ is equal to the sheaf of smooth morphisms of $M$ to $\ex T$. 
This is true in general. A smooth morphism $f:\ex B\longrightarrow \ex T^n$ from any abstract exploded space $\ex B$ is equivalent to the choice of $n$ exploded functions in $\tC  {\ex B}$ corresponding to $f^*(\tilde z_i)$. 

The following special case is worth emphasizing: Given any abstract exploded space $\ex B$,  the sheaf $\mathcal E^{\times}(\ex B)$ can be identified with the sheaf of morphisms to $\ex T$. Think of this as analogous the sheaf of smooth functions on a manifold being identified as the sheaf of smooth maps  to $\mathbb R$.

\item A morphism $f:\ex T^{n}\longrightarrow \ex T$ is given by an exploded function $f=c\e{y}\tilde z^{\alpha}$. This induces an integral affine map called the tropical part of $f$.
\[\totb{\ex f}:\totb{\ex T^{n}}\longrightarrow \totb{\ex T}\] 
\[(\e{x_{1}},\dotsc,\e{x_{n}})\mapsto \e{y+x_{1}\alpha_{1}+\dotsb+ x_{n}\alpha_{n}}\]

\end{enumerate}

\

The next example describes a hybrid object, part `smooth', part `tropical'. 

\begin{example}[$\et 11$]\label{et 11}
\end{example}
 
 The exploded manifold $\et 1{[0,\infty)}:=\et 11$ is more complicated. We shall describe this using the  coordinate $\tilde z\in\mathcal E^\times(\et 11)$. The set of points  $p\in \et 11$ is identified with $\mathbb C^*\e{\mathbb R^+}$ by specifying the value of $\tilde z(p)\in\mathbb C^*\e{\mathbb R^+}$.
 
 Recall the following smooth part homomorphism 
 \[\totl\cdot:\mathbb C^{*}\e{\mathbb R^{+}}\longrightarrow \mathbb C\]
 \[\totl{c\e{x}}:=\begin{cases}0\text{ if }x>0
 \\ c\text{ if }x=0\end{cases}\]
 The above map gives a map from the set of points $p\in \et 11$ to $\mathbb C$. Pulling back the usual topology on $\mathbb C$ defines the topology on ${\et 11}$. We shall refer to $\mathbb C$ as the smooth part of $\et 11$, and use the notation $\totl{\et 11}=\mathbb C$.

 We can write any exploded function $h\in\mathcal E^{\times}\lrb{\et 11}$ as 
 \[h(\tilde z)=f(\totl{\tilde z})\e y\tilde z^\alpha\text{ for }f\in C^\infty(\mathbb C,\mathbb C^*)\text{, and } y\in\mathbb R,\alpha \in\mathbb Z\text{ locally constant.}\]

 The tropical part or large scale of $\et 11$ is $\totb{\et 11}=\e{[0,\infty)}$, and the map from the set of points $\mathbb C^{*}\e{\mathbb R^{+}}$ in $\et 11$ to $\totb{\et 11}$ is given by
 \[c\e{x}\mapsto \totb{c\e{x}}:=\e{x}\] 
 
\includegraphics{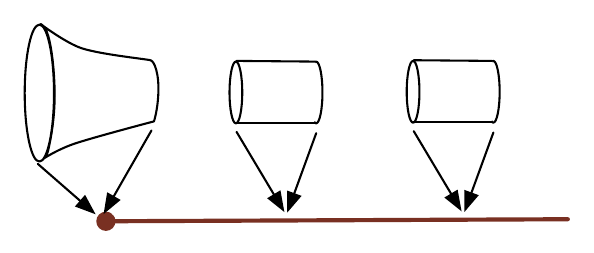}

 This exploded manifold $\et 11$ should be thought of as follows: Over $\e 0\in\totb{\et 11}$, there is  a copy of $\mathbb C^{*}$ which should be considered as a manifold with an asymptotically cylindrical end at $0$, drawn on the left hand side above with this cylindrical end pointing right. Each copy of $\mathbb C^*$ over $\e a\in \totb{\et 11}$ where $a>0$ should be thought of as some `cylinder at infinity'. Note that even though there is a $(0,\infty)$ worth of two dimensional cylinders involved in this exploded manifold, it should still be thought of being two dimensional. This strange feature is essential for the exploded category to have a good holomorphic curve theory. (Actually, we could just as easily had a $\mathbb Q^{+}$ worth of cylinders at infinity, and worked over the semiring $\mathbb C\e{\mathbb Q}$ instead of $\mathbb C\e {\mathbb R}$, but this author prefers the real numbers.)

The exploded manifold $\et 11$ is a kind of hybrid of the last two examples. Restricting to an open set contained inside $\{\totb{\tilde z}=\e 0\}\subset{\et 11}$ we get part of a smooth manifold. Restricting to a subset contained inside $\{\totl{\tilde z}=0\}\subset\et 11$, we get part of $\ex T$.

\

For any abstract exploded space $\ex B$, we can define  $\tCp{\ex B}$ to consist of all functions in $\mathcal E^\times (\ex B)$ which take values only in $\mathbb C^*\e{\mathbb R^+}$. 
Note that a smooth morphism $f:\ex B\longrightarrow\et 11$ from any  exploded manifold $\ex B$  described in examples so far is given by a choice of exploded function 
$f^*(\tilde z)\in \tCp{\ex B}$. This will be true for exploded manifolds in general.

\

\begin{example}[A holomorphic curve in $\ex T^{2}$]
\end{example}
The following is an example of a holomorphic curve in $\ex T^{2}$. Consider the subset of $\ex T^{2}$ where $\tilde z_{1}+\tilde z_{2}+1\in0\e{\mathbb R}$. 

If $\totb{\tilde z_{1}}$ and $\totb{\tilde z_{2}}$ are less than $\totb{1}=\e 0$, then $\tilde z_{1}+\tilde z_{2}+1=1$, so there are no solutions in this region. There are also no solutions where $\totb{\tilde z_{1}}>\totb{\tilde z_{2}}$ and $\totb{\tilde z_{1}}>\e0$ because in this region $\tilde z_{1}+\tilde z_{2}+1=\tilde z_{1}$, which only takes values in $\mathbb C^{*}\e {\mathbb R}$. Similarly, there are no solutions where $\totb{\tilde z_{2}}>\totb{\tilde z_{1}}$ and $\totb{\tilde z_{2}}>\e 0$, because in this region $\tilde z_{1}+\tilde z_{2}+1=\tilde z_{2}$.

Below is a picture of the tropical part of $\ex T^{2}$ with labels indicating simplifications of $\tilde z_{1}+\tilde z_{2}+1$ in various regions.

\includegraphics{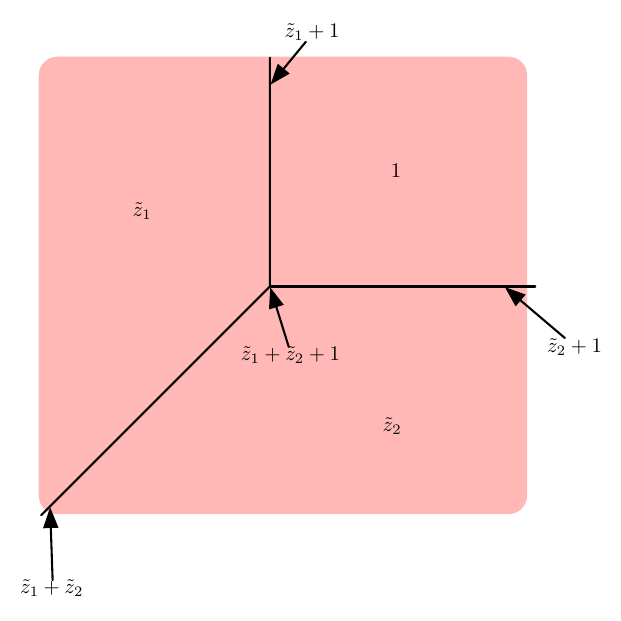}

The tropical part of the subset where $\tilde z_{1}+\tilde z_{2}+1\in 0\e{\mathbb R}$ is the black tropical curve drawn above, which coincides with the subset where the tropical parts of at least two terms of $\tilde z_{1}+\tilde z_{2}+1$ coincide.
This subset where $\tilde z_{1}+\tilde z_{2}+1\in 0\e{\mathbb R}$ is actually a holomorphic curve inside $\ex T^{2}$. We can parametrize this curve using the following 3 coordinate charts modeled on  $\et 11\setminus\{1\}$:
\[\phi_{i}:\et 11\setminus\{1\}\longrightarrow \ex T^{2}\]
\[\phi_{1}(\tilde w)=(-\tilde w,\totl{\tilde w}-1)\]
\[\phi_{2}(\tilde w)=(\totl{\tilde w}-1,-\tilde w)\]
\[\phi_{3}(\tilde w)=(-\tilde w^{-1},\tilde w^{-1}-1)\]

The tropical part of the image of $\phi_{1}$ is the line pointing right in the above picture, the tropical part of $\phi_{2}$ is the line pointing up, and the line pointing left and down is the image of the tropical part of $\phi_{3}$.

 A similar example is considered on page \pageref{tropical curve}.
 
 \

The isomorphism type of a coordinate chart on a smooth manifold only depends on dimension. The isomorphism type of a coordinate chart on an exploded manifold is determined by the dimension and an integral affine polytope. Below we define some integral affine notions.

\begin{defn}[Integral affine]
An integral affine map $\e{\mathbb R^{n}}\longrightarrow \e{\mathbb R^{m}}$ is a map of the form
\[\e{x}\mapsto \e{Mx+y}\]
where $M$ is an integer $n\times m$ matrix and $y\in\mathbb R^{m}$.

An integral affine polytope $P\subset \e{\mathbb R^{m}}$ is a convex polytope with nonempty interior and faces of rational slopes, defined by some finite set of inequalities:
   \[P:=\{\e x \in \e{\mathbb R^{m}}\text{ so that }\e {a_{i}+x\cdot\alpha^{i}}\leq \e{0},\e {a_{i'}+x\cdot\alpha^{i'}}< \e{0} \}\]
where $a_{i}\in\mathbb R$ and $\alpha^{i}\in \mathbb Z^{m}$. 

Integral affine polytopes form a category with morphisms being integral affine maps between polytopes.

Call $P\subset \e{\mathbb R^{m}}$ complete if $P\subset \e{\mathbb R^{m}}$ is complete when given the usual Euclidean metric on $\mathbb R^{m}$ (ie. no strict inequalities are used to define $P$.) Call $P\subset \e{\mathbb R^{m}}$ open if it can be defined using only strict inequalities.

\end{defn}
  
     \begin{example}[$\et mP$]\label{polytope}
     \end{example}
   
     Given an integral affine polytope $P\subset \e{\mathbb R^{m}}$, define the abstract exploded space $\et mP$ as follows:
     
     The set of points in $\et mP$ is equal to the set of points $p\in \ex T^{m}$ so that $\totb p\in P\subset \totb{\ex T^{m}}$. The coordinates $\tilde z_{1},\dotsc,\tilde z_{m}:\ex T^{m}\longrightarrow \mathbb C^{*}\e{\mathbb R}$ then restrict to give coordinate functions
     \[\tilde z_{1},\dotsc,\tilde z_{m}:\et mP\longrightarrow \mathbb C^{*}\e{\mathbb R}\]
     In coordinates, the set of points in $\et mP\subset \ex T^{m}$ is equal to 
     \[\{(c_{1}\e{x_{1}},\dotsc,c_{m}\e{x_{m}})\text{ so that }c_{i}\in\mathbb C^{*} \text{ and }(\e{x_{1}},\dotsc,\e{x_{m}})\in P\}\]
     Consider the collection of exploded functions $\tilde \zeta:=\e {a}\tilde z^{\alpha}$ on $\ex T^{m}$ so that on $\et mP\subset \ex T^{m}$, $\totb{\tilde \zeta}\leq \e 0$. Choose some finite generating set $\{\tilde \zeta_{1},\dotsc,\tilde \zeta_{n}\}$ of these functions so that any other function $\tilde \zeta$ of this type  can be written as $\tilde \zeta=\e{a}\tilde \zeta_{1}^{\beta_{1}}\dotsb \tilde \zeta_{n}^{\beta_{n}}$  where $\beta_{i}\in\mathbb N$ and $a\geq 0$. Recalling that $\totl{ct^{x}}=0$ if $x>0$ and $\totl{ct^{0}}=c$, define the functions $\zeta_{i}:=\totl{\tilde \zeta_{i}}: \et mP\longrightarrow \mathbb C$. Then consider the following map 
     \[(\zeta_{1},\dotsc,\zeta_{n}):=(\totl{\tilde \zeta_{1}},\dotsc,\totl{\tilde \zeta_{n}}):\et mP\longrightarrow \mathbb C^{n}\]
     Give ${\et mP}$ the topology which is the pullback of the usual topology on $\mathbb C^{n}$ under the above map $(\zeta_{1},\dotsc,\zeta_{n})$, and call the image of this map the smooth part of $\et mP$, written as $\totl{\et mP}\subset \mathbb C^{n}$. 
     The sheaf of  exploded functions on $\et mP$ can then be described as functions of the form $f( \zeta_{1},\dotsc, \zeta_{n})\e{a}\tilde z^{\alpha}$, where $f:\mathbb C^{n}\longrightarrow \mathbb C^{*}$ is smooth and $a\in\mathbb R$ and $\alpha\in\mathbb Z^{m}$ are locally constant. (This topology and sheaf of exploded functions is  independent of the choice of generating set $\{\tilde \zeta_{1},\dotsc,\tilde \zeta_{n}\}$.) 
      
      Naturally, the tropical part of $\et mP$ is given by 
     \[\totb{\et mP}=P=\{(\totb{\tilde z_{1}},\dotsc,\totb{\tilde z_{m}})\}\subset \e{\mathbb R^{m}}\]
     
     Simple examples of the above construction are $\et 1{\mathbb R}$ which is equal to  $\ex T$ and $\et 1{[0,\infty)}$ which is equal $\et 11$. (Where no ambiguity is present, we shall write polytopes such as $\e{[0,\infty)}\subset\e{\mathbb R}$ simply as $[0,\infty)\subset\mathbb R$.) 
     
     \begin{remark}\label{monomial description}\end{remark} $\et nn:=\et n{[0,\infty)^{n}}$ should be understood as an $n$-fold product of $\et 11$. If the polytope $P$ contains no entire lines, a good way to describe $\et mP$ is as a subset of $\et nn$. In particular, the monomials $\{\tilde \zeta_{1},\dotsc,\tilde \zeta_{n}\}$ give an injective map
     \[(\tilde \zeta_{1},\dotsc,\tilde \zeta_{n}):\et mP\longrightarrow \et nn\]
     so that any exploded function on $\et mP$ may be described as the pullback of some exploded function on $\et nn$. We may therefore think of $\et mP$  as a subset of $\et nn$. This subset is described by the monomial relations between the $\tilde \zeta_{i}$, which may be written in a finite number of equations of the form 
     \[\tilde \zeta_{1}^{r_{1}}\dotsb\tilde \zeta_{n}^{r_{n}}=\e {a}\]
     so in the case that $P$ contains no entire lines, we may consider $\et mP$ as a subset of $\et nn$ described by setting some monomials equal to $\e {a}$. Similarly, if $P$ is $\mathbb R^{a}$ times a polytope that contains no lines, we may consider $\et mP$ as a subset of $\ex T^{a}\times \et nn$ defined by setting some monomials equal to $\e a$. 
     
     Note also that the smooth part of $\et mP$ may be regarded as the subset of $\mathbb C^{n}$ (which is the smooth part of $\et nn$) defined by setting the same monomials equal to $\totl{\e a}$, with $\zeta_{i}:=\totl{\tilde \zeta_{i}}$ as variables in place of $\tilde\zeta_{i}$.

     \
     
     The functions $\zeta_{i}:=\totl{\tilde \zeta_{i}}$ above which generate the smooth functions will come up often, so make the following definition.
     
     \begin{defn}[Basis of smooth monomials] The smooth monomials on $\et mP$ are the functions $\zeta$ of the form $\totl{c\e a\tilde z^{\alpha}}$. A set $\{\zeta_{1},\dotsc,\zeta_{n}\}$ of smooth monomials is a basis for the smooth monomials on $\et mP$ if every smooth monomial on $\et mP$ can be written as some product of nonnegative powers of the $\zeta_{i}$  times a complex number.
     \end{defn}

     We can construct $\mathbb R^{k}\times \et mP$ similarly to $\et mP$. This has coordinates $(x,\tilde z)\in\mathbb R^{k}\times \et mP$, with the product topology. Similarly, the  smooth part $\totl{\mathbb R^{k}\times \et mP}$ is the product $\mathbb R^{k}\times \totl{\et mP}$.   The exploded manifold $\mathbb R^{n}\times \et mP$ has exploded functions equal to functions of the form $f(x,\totl{\tilde \zeta_{1}},\dotsc,\totl{\tilde \zeta_{n}})\e a\tilde z^{\alpha}$, where $f:\mathbb R^{k}\times\mathbb C^{n}\longrightarrow \mathbb C^{*}$ is smooth and $a\in\mathbb R$ and $\alpha\in\mathbb Z^{m}$ are locally constant.  The tropical part  $\totb{\mathbb R^{n}\times\et mP}$ is equal to $P$, with $\totb{(x,\tilde z)}=\totb{\tilde z}\in P$. The real dimension of $\mathbb R^{k}\times \et mP$ is considered to be $k+2m$.
   
\

\begin{example}[A nontrivial example of $\et mP$]\label{nontrivial}\end{example}
Let $P\subset \e{\mathbb R^{2}}$ be defined by 
\[P:=\{t^{(x,y)}\text{ so that } x\geq 0 \text{ and }x+2y\geq 0\}\] 
A basis for smooth monomials on $\et mP$ is given by \[ \zeta_{1}=\totl{\tilde z_{1}}, \  \zeta_{2}=\totl{\tilde z_{1}\tilde z_{2}^{2}}, \  \zeta_{3}=\totl{\tilde z_{1}\tilde z_{2}}\] The  smooth part $\totl{\et mP}$ of $\et mP$  is the image of the map 
\[({ \zeta_{1}},{ \zeta_{2}},{ \zeta_{3}}):\et 2P\longrightarrow \mathbb C^{3}\]
which is equal to the subset of $\mathbb C^{3}$ where $\zeta_{1}\zeta_{2}=\zeta_{3}^{2}$.

\

\begin{example}[$\et 1{[0,l]}$]\label{t0l}\end{example}

Let $P=\e{[0,l]}\subset \e{\mathbb R}$. We shall use the notation $\et 1{[0,l]}$ for $\et 1P$. In this case a basis for the smooth monomials on $\et 1{[0,l]}$ is given by
$\zeta_{1}=\totl{\tilde z}$ and $\zeta_{2}=\totl{\e{l}\tilde z^{-1}}$. Any smooth morphism $f:\et 1{[0,l]}\longrightarrow \mathbb R$ is equal to $g(\zeta_{1},\zeta_{2})$ for some smooth function $g:\mathbb C^{2}\longrightarrow \mathbb R$. The two functions $\zeta_{1}$ and $\zeta_{2}$ satisfy the relation $\zeta_{1}\zeta_{2}=0$, so $\totl{\et 1{[0,l]}}$ is the union of the two coordinate planes in $\mathbb C^{2}$, which can be considered as two copies of $\mathbb C$ glued over $0$. The function $f$ is also equivalent to the choice of two smooth functions $f_{1}:\mathbb C\longrightarrow \mathbb R$, $f_{2}:\mathbb C\longrightarrow \mathbb R$  so that $f_{1}(0)=f_{2}(0)$. The equivalence is given by 

\[f(\tilde z)=\begin{cases} f_{1}(\zeta_{1})\text{ if }\totb{\tilde z}=\e0
\\ f_{1}(0)=f_{2}(0)\text{ if }\e{0}>\totb{\tilde z}>\e l
\\ f_{2}(\zeta_{2})\text{ if }\totb{\tilde z}=\e l\end{cases}\]

The three different possibilities above correspond to three different `strata' of $\et 1{[0,l]}$ where $\totl{\tilde z}$ is in $\e 0,\ \e{(0,l)}$ or $\e l$. The ability to solve problems such as differential equations on exploded manifolds stratum by stratum is part of the usefulness of the exploded category, as complicated problems can be broken into simple pieces. 

\

Before continuing, the reader should be able to do the following easy exercises:
\begin{enumerate}
\item A morphism  $\mathbb R^{n}\times \et mP \longrightarrow \et kQ$ is equivalent to the choice of  $k$ exploded functions $f_{1},\dotsc, f_{k}\in \mathcal E^{\times}\lrb{\mathbb R^{n}\times \et mP}$ so that 
\[(\totb{f_{1}},\dotsc,\totb{f_{k}})\in Q\]
\item A morphism $ \et mP \longrightarrow \mathbb R^{k}$ is equivalent to a choice of continuous map $f: \et mP\longrightarrow \mathbb R^{k}$ so that there exists some smooth map 
\[\hat f:\mathbb C^{n}\longrightarrow\mathbb R^{k}\]
so that
\[f=\hat f( \zeta_{1},\dotsc,{ \zeta_{n}})\] 

\item Any morphism $f:\mathbb R^{n}\times \et mP\longrightarrow \mathbb R^{k}\times\et lQ$ induces an integral affine map $\totb f:P\longrightarrow Q$ so that the following diagram commutes
\[\begin{array}{ccc}\mathbb R^{n}\times \et mP&\xrightarrow{f} &\mathbb R^{k}\times\et lQ
\\ \downarrow & & \downarrow
\\ P &\xrightarrow{\totb f} &Q\end{array}\]

\item $\et mP$ is isomorphic as an abstract exploded space to $\et mQ$ if and only if $P$ is isomorphic as an integral affine polytope to $Q$. 

\item  Any morphism $f:\mathbb R^{n}\times \et mP\longrightarrow \mathbb R^{k}\times\et lQ$ induces a continuous map $\totl f:\totl{\mathbb R^{n}\times\et mP}\longrightarrow \totl{\mathbb R^{k}\times \et lQ}$
so that the following diagram commutes
\[\begin{array}{ccc}\mathbb R^{n}\times \et mP&\xrightarrow{f} &\mathbb R^{k}\times\et lQ
\\ \downarrow & & \downarrow
\\ \totl{\mathbb R^{n}\times\et mP}&\xrightarrow{\totl f} & \totl{\mathbb R^{k}\times \et lQ}\end{array}\]

\end{enumerate}

\

\begin{defn}[Exploded manifold] A smooth exploded manifold $\ex B$ is an abstract exploded space locally isomorphic to $\mathbb R^{n}\times \et mP$. In other words, for all $p\in\ex B$ there exists some open neighborhood $U$ of $p$ and some $\mathbb R^{n}\times\et mP$  so that $U$ is isomorphic as an abstract exploded space to $\mathbb R^{n}\times\et mP$. 

A smooth map $\ex A\longrightarrow \ex B$ of exploded manifolds is a morphism  $\ex A\longrightarrow\ex B$ of abstract exploded spaces. 
\end{defn}

Lemma \ref{sub coordinates} on page \pageref{sub coordinates} proves that the above definition is equivalent to defining an exploded manifold as an abstract exploded space locally isomorphic to  open subsets of $\mathbb R^{n}\times\et mP$.

\

\begin{defn}[Smooth part] The smooth part $\totl {\ex B}$ of a smooth exploded manifold is the Hausdorff\footnote{After publishing,  I added  an extra condition to Definition \ref{abstract es} to ensure that $\totl{\ex B}$ is Hausdorff. For example, we do not want two copies of $\mathbb R$ glued over an open interval to be a valid exploded manifold. } topological space which is the	quotient of the topological space $\ex B$ by  the equivalence relation $p \simeq q$ if every open subset of $\ex B$ which contains $p$ contains $q$. (The fact that this relation is symmetric follows from the fact that $\ex B$ is locally isomorphic to $\mathbb R^{n}\times \et mP$. ) \end{defn}

It should be clear that $\totl\cdot$ gives a functor from the category of smooth exploded manifolds to the category of Hausdorff topological spaces, and that the smooth part of coordinate charts described earlier agrees with the above definition.

 As every smooth map $\ex B\longrightarrow \mathbb R$ factors through $\totl{\ex B}$, the the smooth part $\totl{\ex B}$ of $\ex B$ also has a kind of smooth structure given by the sheaf of functions on $\totl{\ex B}$ which pull back to smooth maps $\ex B\longrightarrow\mathbb R$.

\

A compact exploded manifold $\ex B$ is an exploded manifold for which $\totl{\ex B}$ is compact. A stronger notion which agrees better with the notion of compactness for smooth manifolds is given below:

\begin{defn}[Complete]\label{complete} An exploded manifold $\ex B$ is complete if $\totl{\ex B}$ is compact, and every smooth map of 
$\et 1{(0,l)}$ into $\ex B$ extends to a smooth map of $\et 1{[0,l]}$ into $\ex B$.

A map $f:\ex B\longrightarrow \ex C$ is proper if $\totl f:\totl{\ex B}\longrightarrow \totl{\ex C}$ is proper.
A map $f:\ex B\longrightarrow \ex C$ is complete if it is proper and every smooth map  $\gamma:\et 1{(0,l)} \longrightarrow \ex B$ extends to a smooth map  $\et 1{[0,l]}\longrightarrow \ex B$ if and only if $f\circ \gamma$ extends to a smooth map $\et 1{[0,l]}\longrightarrow \ex C$.
\end{defn}

For example, $\ex T$ is complete, and any compact manifold is complete when thought of as an exploded manifold. The inclusion $\et 1{[0,1)}$ into $\et 1{[0,\infty)}$ is an example of a map which is proper but not complete. An equivalent condition for an exploded manifold to be complete is that it is compact and locally isomorphic to $\mathbb R^{n}\times \et mP$ where the polytope $P\subset \mathbb R^{m}$ is closed, (and hence complete when $\mathbb R^{m}$ is given the standard metric). The second property of completeness can always be tested locally by looking at the tropical part of a coordinate chart or map.  For a map $f:\mathbb R^{n}\times \et mP\longrightarrow \mathbb R^{n'}\times \et {m'}Q$, this second property holds if and only if 
the inverse image under the tropical part of $f$ of any complete subset of $Q$ is a complete subset of $P$. 


\

\section{Stratified structure}

\begin{defn}[Faces and strata of polytopes]
A face $F$ of an integral affine polytope $P$ is a subset of $P$ which is defined by an equation 
\[F:\{\e{x}\in P \text{ so that }x\cdot\alpha \text{ is minimal}\}\]
 
 A stratum $S$ of an integral affine polytope $P$ is a subset of $P$ which is equal to a face of $P$ minus all proper subfaces. This can also be thought of as the interior of a face of $P$.
 
 \end{defn}
For example under this definition, the faces of a triangle are the entire triangle, the closure of each edge and each vertex.
 The strata of a triangle consist of  the interior of the triangle, the interior of each edge and each vertex.

\begin{lemma}\label{tropical iso}
If $f:\mathbb R^{n'}\times \et {m'}{P'}\longrightarrow \mathbb R^{n}\times\et {m}{P}$ is an isomorphism onto an open subset, then $\totb f: P'\longrightarrow P$ is an isomorphism of $P'$ onto a face of $P$.
\end{lemma}
\pf

Let $U$ be an open subset of $\mathbb R^{n}\times \et mP$, and let $\totb U$ be the image of $U$ in $\totb P$. We shall show that $\totb U$ is a union of faces of $P$. It will follow from this that the image of $\totb f$ must be a face of $P$.

First consider the case of an open subset of $\et 1{[0,l]}$. Any open subset of $\et 1{[0,l]}$ is of the form 
\[\{(\totl{\tilde z},\totl{\e l\tilde z^{-1}})\in U'\}\]
for some open subset $U'\subset \mathbb C^{2}$. Therefore, if $U$ contains some point where $\totb{\tilde z}\in\e{(0,l)}$, then $(0,0)\in U'$, so all points where $\totb{\tilde z}\in\e{(0,1)}$ are in $U$. Also,  for some $\epsilon>0$, $(0,\epsilon)$ and $(\epsilon,0)$ are in $U'$, so the points where $\tilde z=\epsilon\e0$ and $\tilde z=\frac 1\epsilon \e l$ are in $U$, so $\totb U=\e{[0,l]}$.  

Now in the general case, given any point $p\in U$ with coordinates
\\ $(x,c_{1}\e {a_{1}},\dotsc, c_{m}\e {a_{m}})$, suppose that 
\[\totb \gamma:[0,l]\longrightarrow P\] 
is any integral affine map so that there is some point $t_{0}\in(0,l)$ for which $\totb \gamma(t_{0})=\totb p$. Then we can construct a corresponding smooth map $\gamma:\et1{[0,l]}\longrightarrow\mathbb R^{n}\times \et mP$ with tropical part $\totb \gamma$ so that $\gamma(1\e {t_{0}})=p$ as follows: 
\[\text{ if }\totb \gamma(t)=v+tw\]
\[\text{define }\gamma(\tilde z)=(x,c_{1}\e{v_{1}}\tilde z^{w_{1}},\dotsc,c_{m}\e{v_{m}}\tilde z^{w_{m}})\] 

Therefore, as $\gamma^{-1}(U)$ is an open subset of $\et 1{[0,l]}$, $\totb{\gamma^{-1}(U)}$ is the entire interval $[0,l]$. Therefore, if $\totb U$ contains a point $p$ on the interior of some interval in $P$, $\totb U$ contains the entire interval. It follows that if $\totb U$ contains a point $p$ on a face of $P$, then $\totb U$ contains the entire face. Therefore $\totb U$ is a union of faces of $P$.
 If $U$ is isomorphic to $\mathbb R^{n'}\times \et {m'}{P'}$, it follows that $\totb U$ is isomorphic to $P'$. The only way for a union of faces of $P$ to be isomorphic to a polytope $P'$ is to be a single face of $P$, so $\totb U$ is isomorphic to a face of $P$.

\stop

\begin{defn}[Strata]A stratum of $\mathbb R^{n}\times\et mP$ is a subset 
\[\mathbb R^{n}\times\et mS:=\{p\in \mathbb R^{n}\times\et mP \text{ so that }\totb p\in S\subset P \}\]
where $S$ is a stratum of $P$.

A stratum of an exploded manifold $\ex B$ is a subset of $\ex B$ which is an equivalence class of the following equivalence relation: Say that $p$ and $p'$ are in the same stratum of $\ex B$ if there exist a finite sequence of points $p_{0}=p, p_{1},\dotsc,p_{n}=p'\in\ex  B$ and neighborhoods $p_{i},p_{i-1}\in U_{i}$ isomorphic to $\mathbb R^{n}\times \et mP$ so that $p_{i}$ and $p_{i-1}$ are in the same stratum of $U_{i}$ for $i=1,\dotsc n$.
\end{defn}

Observe that Lemma \ref{tropical iso} implies that the strata of $\mathbb R^{n}\times\et mP$ considered as an exploded manifold are the same as the strata described in the initial part of the above definition. 

\

 Lemma \ref{sub coordinates} from page \pageref{sub coordinates}  states that
any open subset $U\subset\ex B$ of an exploded manifold  is an exploded manifold. In particular, given any point $p\in U$, there exists an open neighborhood $U'$ of $p$ contained inside $U$ which is isomorphic to $\mathbb R^{n}\times \et mP$. 

Lemma \ref{tropical iso} together with Lemma \ref{sub coordinates} implies that all points in a given stratum are contained in coordinate charts of the form $\mathbb R^{n}\times \et mP$, where $P$, $n$ and $m$ are fixed. The points in these coordinate charts which are in this stratum correspond to $\mathbb R^{n}\times \et m{P^{\circ}}\subset \mathbb R^{n}\times \et m P$ where $P^{\circ}\subset P$ is the interior stratum of $P$. Therefore, each stratum $\ex B_{i}$ of $\ex B$ is a connected exploded manifold locally isomorphic to $\mathbb R^{n}\times \et m{P^{\circ}}$, where $P^{\circ}$ is a fixed open integral affine polytope. As the smooth part of  $\mathbb R^{n}\times \et m{P^{\circ}}$ is $\mathbb R^{n}$, the smooth part $\totl{\ex B_{i}}$ of the stratum $\ex B_{i}$ is a connected $n$ dimensional smooth manifold.

Open subsets of $\et mP$ correspond to closed subsets of $P$. In particular, consider the stratum $\et m{P^{\circ}}$ of $\et mP$ corresponding to the interior of $P$. This is a closed subset of $\et mP$  which corresponds to a single point in $\totl{\et mP}$. As the tropical part of any open subset of $\et mP$ is a union of faces of $P$, the closure of any stratum of $\totl{\et mP}$ contains this point $\totl{\et m{P^{\circ}}}$, and the closure of any stratum of $\et mP$ contains $\et m{P_{\circ}}$. In general, if $S$ is a stratum of $P$, the closure of $\et m{S}$ is equal to the union of $\et m{S'}$ for all strata $S'$ whose closure contains $S$ (so closure in $\totl{\et mP}$ and $\et mP$ goes the opposite direction to closure in $\totb{\et mP}=P$.)  

Therefore, the closure of each stratum $\ex B_{i}$ in $\ex B$ is a union of strata. This makes the smooth part $\totl{\ex B}$ of $\ex B$ a stratified space, with each stratum a smooth manifold, and the closure of each stratum a union of manifolds with even codimension. 

\

Consider the tropical part of the stratified structure of $\ex B$. Each stratum $\ex B_{i}$ is locally isomorphic to $\mathbb R^{n}\times \et m{P^{\circ}}$ and a neighborhood of $\ex B_{i}$ in $\ex B$ is locally isomorphic to $\mathbb R^{n}\times \et mP$ for some fixed polytope $P$ with interior $P^{\circ}$. The map $\ex B_{i}\longrightarrow \totl{\ex B_{i}}$ is a $\et m{P^{\circ}}$ bundle over the manifold $\totl{\ex B_{i}}$, and  a neighborhood of $\ex B_{i}$ is isomorphic to a $\et m P$ bundle over $\totl{\ex B_{i}}$. (This last fact may be proved using equivariant coordinate charts constructed in \cite{cem}.)

We can therefore associate to $\ex B_{i}$ a flat integral affine $P$-bundle $\totl{\ex B_{i}}\rtimes P$ over the manifold $\totl{\ex B_{i}}$. Monodromy around any loop in $\totl{\ex B_{i}}$ gives an automorphism of $P$. (We shall later restrict to the case that $\ex B$ is basic, in which case this monodromy will always be trivial.) As monodromy around a loop in a single coordinate chart is always trivial, choosing a base point in $\totl{\ex B_{i}}$ gives a homomorphism from the fundamental group of the closure of $\totl{\ex B_{i}}$ to  the group of automorphisms of $P$. If $\ex B_{j}$ is in the closure of $\ex B_{i}$, then any choice of path from a base point in $\totl{\ex B_{i}}$ to a base point in $\ex B_{j}$ gives an identification of $P$ as a face of the polytope associated with $\ex B_{j}$. Again, in the case that $\ex B$ is basic, this identification will not depend on the path chosen.

\begin{defn}[Tropical structure]
The tropical structure of $\ex B$ is a category $\ex  B_{T}$ with a functor $\mathcal P$ to the category of integral affine polytopes  so that
\begin{enumerate}
\item The objects in $\ex B_{T}$ correspond to points in $\totl{\ex B}$.
\item The morphisms from $p$ to $q$ in $\ex B_{T}$ correspond to homotopy classes of continuous paths $\gamma$ from $p$ to $q$ in $\totl{\ex B}$ so that if $s\geq t$, $\gamma(s)$ is contained in the closure of the stratum of $\totl{\ex B}$ containing $\gamma(t)$. Composition of morphisms in $\ex B_{T}$ corresponds to composition of homotopy classes of paths.
\item $\mathcal P(p)$ is the integral affine polygon associated to the stratum containing $p$, and $\mathcal P(\gamma)$ is the inclusion discussed above given by parallel transport along $\gamma$.
\end{enumerate}
\end{defn}

\begin{defn}[tropical part]
The tropical part $\totb{\ex B}$ of an exploded manifold $\ex B$ is a stratified topological space in which each stratum is given the structure of the quotient of an integral affine polytope by some automorphisms.

As a topological space, the tropical part of $\ex B$ is  defined as the quotient of the disjoint union of all $P\in \mathcal P(\ex B_{T})$ by all the inclusions  $P\longrightarrow Q$ in $\mathcal P(\ex B_{T})$.   Each stratum $\ex B_{i}$ of $\ex B$ corresponds to a stratum of $\totb{\ex B}$ equal to the image in $\totb{\ex B}$ of the interior of any polytope in $\mathcal P(\ex B_{i})$.

\end{defn}

Note that the topologies on $\totl{\ex B}$ and $\totb{\ex B}$ are in some sense dual: The closure of $\totl{\ex B_{i}}\subset \totl{\ex B}$ contains $\totl{\ex B_{j}}$ if and only if the closure of $\totb{\ex B_{j}}\subset\totb{\ex B}$ contains $\totb{\ex B_{i}}$. This property is analogous to the fact that the strict inequality $\totb w<\e 0$ is equivalent to the equation $\totl w=0$, and the inequality $\totl w\neq 0$ is equivalent to $\totb w=\e 0$.

\

The construction of $\ex B_{T}$  is functorial: any map $f:\ex B\longrightarrow \ex C$ induces a functor $ f_{T}:\ex B_{T}\longrightarrow\ex C_{T}$ and for every object $x$ of $\ex B_{T}$, a map $\mathcal Pf: \mathcal P(x)\longrightarrow \mathcal P(f_{T}x)$ so that diagrams of the following type commute:
\[\begin{array}{ccc}\mathcal P(x)&\xrightarrow{\mathcal P(\gamma)} &\mathcal P(y)\\
\downarrow\mathcal Pf &&\downarrow\mathcal P f
\\ \mathcal P(f_{T}x)&\xrightarrow {\mathcal P(f_{T}\gamma)}&\mathcal P(f_{T}y)\end{array}\]
 In particular, an object $x$ of $\ex B_{T}$ corresponds to a point $x\in\totl{\ex B}$, and $f_{T}x$ corresponds to the point $\totl{f(x)}\in \totl{\ex C}$.  We may choose a coordinate chart containing $x$ with tropical part $\mathcal P(x)$ and a coordinate containing $\totl{f(x)}$ with tropical part $\mathcal P(f_{T}x)$. Then the tropical part of our map in coordinates is the map $\mathcal Pf:\mathcal P (x)\longrightarrow \mathcal P(f_{T}x)$.
Call $f_{T}$ and $\mathcal Pf$ the tropical structure of $f$.

 The construction of $\totb{\ex B}$ is also functorial. Given a map $f:\ex B\longrightarrow \ex C$, our maps $\mathcal Pf$ give a continuous map $\totb f:\totb{\ex B}\longrightarrow \totb{\ex C}$ called the tropical part of $f$.  

The following is a case when the tropical part $\totb{\ex B}$ is simply  a union of polytopes glued along faces.

\begin{defn}[Basic]
The exploded manifold $\ex B$ is basic if there is at most one morphism between any two polytopes in $\mathcal P(\ex B_{T})$.
\end{defn}

Making the assumption that $\ex B$ is basic very often simplifies combinatorial aspects of an argument. All examples of exploded manifolds discussed so far have been basic.

Observe that if $\ex B$ is basic and every polytope in $\mathcal P(\ex B_{T})$ is a Delzant polytope (in other words it is locally isomorphic as an integral affine space to an open subset of the integral affine space $[0,\infty)^{m}$), then the closure of each stratum of $\totl B$ is a manifold. 

\

One example of an exploded manifold which is not basic is the exploded manifold constructed by taking the quotient of $\mathbb R\times\ex T$ by the action $(x,\tilde z)\mapsto (x+1,\e 1\tilde z)$. This exploded manifold has a single stratum. The smooth part of this exploded manifold is a circle, the tropical structure  has nontrivial monodromy around this circle, so this exploded manifold is not basic.  

A second example is given by gluing $\et 1{[0,1]}$ to itself via the map $\tilde z\mapsto \e 1\tilde z^{-1}$. This exploded manifold  has two strata - one stratum with tropical part a point and smooth part a two punctured sphere, one with tropical part an interval, and smooth part a point. There are two different inclusions of the tropical part of the first stratum into the interval $[0,1]$ associated with the second stratum so this exploded manifold is not basic. 

\includegraphics{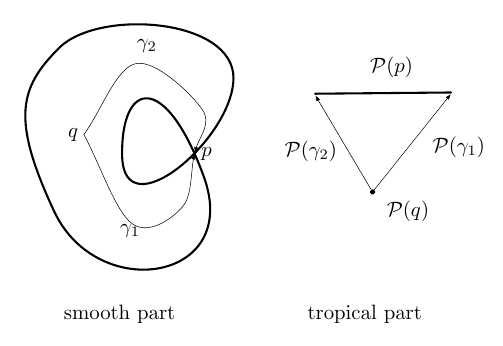}
    
     The following is an example which is a  model for a neighborhood of a stratum with smooth part $M$ and tropical part $P$ in a basic exploded manifold:
     
     \begin{example}[Model for neighborhoods of  strata in basic exploded manifolds]\label{polytope2}\end{example}
     Given a smooth manifold $M$, a polytope $P\subset \mathbb R^{m}$, and  $m$ complex line bundles on $M$, we can construct the exploded manifold $M\rtimes \et mP$ as follows: Denote by $E$ the corresponding total space of our $m$ $\mathbb C^{*}$ bundles over $M$. This has a smooth free $\lrb{\mathbb C^{*}}^{m}$ action. The exploded manifold $\et mP$ also has a $\lrb{\mathbb C^{*}}^{m}$ action given by multiplying the coordinates $\tilde z_{1},\dotsc,\tilde z_{m}$ by the coordinates of $\lrb{\mathbb C^{*}}^{m}$. Construct the exploded manifold $M\rtimes \et mP$  by taking the quotient of $E\times \et mP$ by the action of $\lrb{\mathbb C^{*}}^{m}$ by $(c,c^{-1})$. As the action of $\lrb{\mathbb C^{*}}^{m}$ is trivial on $\totb{\et mP}$, the tropical part of $M\rtimes\et mP$ is still defined, and is equal to $P$.
     
     Alternately, choose coordinate charts on $E$ equal to $U\times\lrb{\mathbb C^{*}}^{m}\subset \mathbb R^{n}\times \mathbb C^{n}$. The transition maps are of the form \[(u,z_{1},\dotsc,z_{m})\mapsto(\phi(u),f_{1}(u)z_{1},\dotsc,f_{m}(u)z_{m})\] Replace these coordinate charts with $U\times \et mP\subset\mathbb R^{n}\times\et mP$, and replace the above transition maps with maps of the form \[(u,\tilde z_{1},\dotsc,\tilde z_{m})\mapsto(\phi(u),f_{1}(u)\tilde z_{1},\dotsc,f_{m}(u)\tilde z_{m})\] The map to the tropical part $\totb{M\rtimes\et mP}:=P$ in these coordinates is given by 
     \[\totb{(u,\tilde z_{1},\dotsc,\tilde z_{n})}=(\totb{\tilde z_{1}},\dotsc,\totb{\tilde z_{n}})\in P\]

\section{The explosion functor and log geometry}\label{expl}

The explosion functor is a functor from a category of complex manifolds with normal crossing divisors to the category of (holomorphic) exploded manifolds. The explosion functor is an important source of examples of exploded manifolds.  At the end of this section we will phrase this functor in the language of log geometry in terms of a kind of base change.

\begin{defn}[Holomorphic exploded manifold]
Given  a connected open subset $U\subset \et mP$, call a smooth map $f:U\longrightarrow \ex T$ holomorphic if it is equal to  $\tilde z^{\alpha}g(\zeta_{1},\dotsc,\zeta_{m})$ where $g$ is holomorphic and $\zeta_{i}$ are smooth monomials.  

A holomorphic exploded manifold is an abstract exploded space locally isomorphic to an open subset of $\et mP$ with the sheaf of holomorphic maps to $\ex T$.
\end{defn}

Suppose that we have a complex manifold $M$ along with a collection of complex codimension $1$ immersed complex submanifolds $N_i$ so that $N_i$ intersect themselves and each other transversely. (We shall call this a complex manifold with normal crossing divisors.) Then there is a complex exploded manifold $\expl (M)$ called the explosion of $M$. 

We define $\expl(M)$ as follows: choose holomorphic coordinate charts on  $M$ which are equal to balls inside $\mathbb C^n$, so that the image of the submanifolds $N_i$ are equal to the submanifolds $\{z_i=0\}$. Then replace a coordinate chart $U\subset \mathbb C^n$ by a coordinate chart $\expl U$ in $(\et 11)^{n}$ with coordinates $\tilde z_i$ so that 
\[\expl U:=\{\tilde z\text{ so that }\totl{\tilde z}\in U\}\]
 Define transition functions as follows: the old transition functions are all of the form \[f_i(z)=z_jg(z)\] where $g$ is holomorphic  and non vanishing. Replace this with \[\expl f_i(\tilde z)=\tilde z_jg(\totl{\tilde z})\] which is then a smooth exploded function. If $f'$ is another transition function with $f'_{j}(z)=z_{k}g'(z)$, then where defined,
 \[\expl(f\circ f')_{i}(\tilde z)=\tilde z_{k}g'(\totl{\tilde z})g(f'(\totl{\tilde z}))=\lrb{\expl f\circ\expl f'}_{i}(\tilde z)\]
  Therefore the explosion of the old transition functions give transition functions which define $\expl (M)$ as a holomorphic  exploded manifold.
 
  \

 In the case that all the submanifolds $N_{i}$ are embedded (as opposed to simply immersed), $\expl (M)$ is a basic exploded manifold. The tropical part of $\expl (M)$ has one vertex for each connected component of $M$, an edge for each submanifold, a face for each intersection, and a $k$ dimensional face for each $k$-fold intersection. 
 
\includegraphics{iec1.pdf}

One natural way to view a complex manifold $M$ with normal crossing divisors $N_{i}$ is as a log space, which in this case means the complex manifold $M$ along with the sheaf of holomorphic functions on $M$ which do not vanish off the divisors $N_{i}$. From this perspective, maps between complex manifolds with normal crossing divisors should be the holomorphic maps which pull back holomorphic functions which do not vanish off divisors to functions of the same type. 

For example a map from a convex open subset $U$ of $\mathbb C^{n}$ with normal crossing divisors given by the coordinate planes to $\mathbb C$ with the divisor $0$ must be in the form \[f(z)=z^{\alpha}g(z)\] where $g$ is holomorphic and nonvanishing and $\alpha\in \mathbb N^{n}$. We can define $\expl f:\expl U\longrightarrow \et 11$ by
\[\expl f(\tilde z)=\tilde z^{\alpha}g(\totl{\tilde z})\]
The explosion of any map between complex manifolds with normal crossing divisors can be defined similarly, and as checked above in the special case of transition functions $\expl(f\circ f')=\expl f\circ\expl f'$.

\

We shall now describe the link between exploded manifolds and log geometry in more detail. Readers not familiar with log geometry may safely skip the remainder of this section, and readers desiring more details on the link between log Gromov-Witten invariants and exploded Gromov-Witten invariants should also consult \cite{elc}.

 Suppose that $\ex B$ is a holomorphic  exploded manifold with tropical structure $\mathcal P(\ex B_{T})$ that contains no polytopes that contain an entire affine line. We can regard $\ex B$ as a log space $L(\ex B)$ as follows: The smooth part, $\totl{\ex B}$
is a kind of singular complex manifold with the sheaf $\mathcal O(\ex B)$ of holomorphic functions given by the sheaf of holomorphic maps of $\ex B$ to $\mathbb C$ regarded as a holomorphic exploded manifold. As well as the  sheaf $\mathcal O(\ex B)$ of holomorphic functions on $\totl{\ex B}$, there is a sheaf   $\tCp{\ex B}$ of holomorphic maps of $\ex B$ to $\et 11$, which gives a sheaf of monoids on $\totl{\ex B}$. The third ingredient needed for a log space is a homomorphism $\tCp{\ex B}\longrightarrow \mathcal O(\ex B)$. The required homomorphism is provided by the smooth part homomorphism $\totl\cdot:\mathbb C^{*}\e{[0,\infty)}\longrightarrow \mathbb C$, so $f\in\tCp{\ex B}$ is sent to $\totl f$. The log space $L(\ex B)$ can be thought of as the data $(\totl{\ex B},\mathcal O(\ex B),\tCp{\ex B},\totl{\cdot})$.

 This construction is functorial: given any holomorphic map $f:\ex B\longrightarrow\ex C$ there is a natural map $L(f):L(\ex B)\longrightarrow L(\ex C)$ of log spaces given by the holomorphic map $\totl f:\totl{\ex B}\longrightarrow\totl{\ex C}$ and the pullback map 

\[\begin{array}{cccc}f^{*}:&\totl{f}^{*}\lrb{\tCp{\ex C}}&\longrightarrow &\tCp{\ex B}
\\ &g &\mapsto & g\circ f\end{array}\]

For example a point $p$ considered as a holomorphic exploded manifold corresponds to the log space $L(p)$ which is a point with the monoid $\mathbb C^{*}\e{[0,\infty)}$. As each exploded manifold $\ex B$ comes with a unique map to $p$, the corresponding log space $L(\ex B)$ comes with a canonical map to $L(p)$. As such, $L(\ex B)$ is correctly regarded as a log space over $L(p)$. 

\

We shall now see that $\ex B$ can be recovered from $L(\ex B)\longrightarrow L(p)$, so $\ex B$ can be regarded as a log space $L(\ex B)$ over $L(p)$. To do this, we must reconstruct the set of points in $\ex B$ and the sheaf $\tC{\ex B}$ of holomorphic exploded functions on $\ex B$. Suppose that $P$ is a $m$-dimensional polytope which contains no affine lines. Then $P$ is isomorphic to a polytope contained entirely inside $[0,\infty)^{m}$. Therefore any integral affine map from $P$ to $\mathbb R$ is a finite sum of integral affine maps to $[0,\infty)$ with integral affine maps to $(-\infty,0]$.  It follows that for any $U\subset\et mP$, the group $\tC U$ of holomorphic exploded functions on $U$ is generated by the set $\tCp U$ of holomorphic maps to $\et 11$. Therefore, if the tropical structure $\mathcal P(\ex B_{T})$ contains no polytopes that contain an affine line, then $\tC{\ex B}$ can be recovered as a sheaf of groups from $\tCp{\ex B}$. In particular, $\tC{\ex B}$ is the sheaf of groups generated by the sheaf $\tCp{\ex B}$ of monoids. We must now recover the set of points in $\ex B$ and be able to interpret $\tC{\ex B}$ as a sheaf of functions on this set of points. 

As we are considering $L(\ex B)$ as a log space over $L(p)$,  the natural candidate for the set of points in $\ex B$ is the set of maps of $L(p)$ into $L(\ex B)$ so that the composition with $L(\ex B)\longrightarrow L(p)$ is the identity map. The information in the map $L(\ex B)\longrightarrow L(p)$ is the inclusion of $\mathbb C^{*}\e{[0,\infty)}$ in $\tCp{\ex B}$ as the set of constant maps to $\et 11$. A map \[f:L(p)\longrightarrow L(\ex B)\] is equivalent to a map $\totl f:p\longrightarrow \totl{\ex B}$ and a homomorphism $f^{*}$ from the stalk of $\tCp{\ex B}$ at $\totl{f}(p)$ to $\mathbb C^{*}\e{[0,\infty)}$ so that $\totl {f^{*}(g)}=\totl g(\totl f(p))$.

 As our map $f$ must be compatible with the map $L(\ex B)\longrightarrow L(p)$, we must restrict to the case of homomorphisms that are the identity on the constant maps, so $f^{*}(c\e a)=c\e a$.
The homomorphism $f^{*}$ extends uniquely to a homomorphism from the stalk of $\tC{\ex B}$ at $\totl f(p)$ to $\mathbb C^{*}\e {\mathbb R}$, therefore $\tC{\ex B}$ can be regarded as a sheaf of functions on this set of points.  

 Choose a coordinate chart on $\ex B$ with image in $\totl{\ex B}$ containing $\totl f(p)$ so that this coordinate chart is isomorphic to an open subset of $\et mP$ where $P$ is contained in $[0,\infty)^{m}$. Then 
\[f^{*}\lrb{c\e a\tilde g(\totl{\tilde z})\tilde z_{1}^{\alpha_{1}}\dotsb \tilde z_{m}^{\alpha_{m}}}=c\e a g(\totl f(p)) f^{*}(\tilde z_{1})^{\alpha_{1}}\dotsb f^{*}(\tilde z_{m})^{\alpha_{m}}\]
So $f^{*}$ is entirely determined by $\lrb{f^{*}(\tilde z_{1}),\dotsc,f^{*}(\tilde z_{m})}\in(\mathbb C^{*}\e{[0,\infty)})^{m}$.
The fact that $f^{*}$ is a homomorphism from the stalk of $\tCp{\ex B}$ at $\totl{f}(p)$ to $\mathbb C^{*}\e{[0,\infty)}$ so that $\totl {f^{*}(g)}=\totl g(\totl f(p))$ implies that $\lrb{f^{*}(\tilde z_{1}),\dotsc,f^{*}(\tilde z_{m})}$ is the coordinates for a point $f(p)\in \et mP$ so that $\totl{f(p)}=\totl f(p)$, so $f^{*}$ is simply evaluation at the point $f(p)$.   Therefore the set of maps  $L(p)\longrightarrow L(\ex B)$ as log spaces over $L(p)$ corresponds to the set of maps  from $p$ into $\ex B$. A bijection is given by the functor $L$. 

Similarly, the functor $L$ gives a bijection from the set of holomorphic maps $\ex B\longrightarrow \ex C$ to the set of maps $L(\ex B)\longrightarrow L(\ex C)$ as log spaces over $L(p)$, so long as no polytopes in $\mathcal P(\ex B_{T})$ or $\mathcal P(\ex C_{T})$ contain an entire affine line.  This is because a map $L(\ex B)\longrightarrow L(\ex B)$ automatically gives a map from the set of points in $L(\ex B)$ to the set of points in $L(\ex C)$ which is  continuous with respect to the topologies on $\totl{\ex B}$ and $\totl{\ex C}$, and pulls back holomorhpic exploded functions to holomorphic exploded functions.

The above discussion implies that a subcategory of  holomorphic exploded manifolds can be regarded as log spaces over $L(p)$. A complex manifold with normal crossing divisors can be regarded as a log space $M^{\dag}$ over a point $\spec \mathbb C$ with `sheaf of monoids' given by $\mathbb C^{*}$. There is a canonical map of log space  $L(p)\longrightarrow \spec \mathbb C$ corresponding to the inclusion of $\mathbb C^{*}$ into $\mathbb C^{*}\e{[0,\infty)}$. The explosion functor can be regarded as a base change from log spaces over $\spec \mathbb C$ to log spaces over $L(p)$ given by this map.
\[\begin{array}{ccc}L(\expl M)&\longrightarrow &M^{\dag}
\\\downarrow & &\downarrow
\\ L(p)&\longrightarrow &\spec \mathbb C\end{array}\]
In other words, $\expl M$ considered as a log space is the fiber product of $M^{\dag}\longrightarrow \spec \mathbb C$ with the map $L(p)\longrightarrow \spec \mathbb C$.

The algebraic geometry of log schemes over $L(p)$ is probably a very interesting direction for further research.

  \section{Tangent space}
   
To define the tangent space of an exploded manifold $\ex B$, we shall need to use the sheaf of smooth real valued functions, and we shall need to be able to add together two exploded functions.

 \begin{defn}[Smooth function]
 The sheaf of smooth functions $C^\infty(\ex B)$ is the sheaf of smooth morphisms of $\ex B$ to $\mathbb R$ considered as a smooth exploded manifold. 
 \end{defn}
%
%
 
 \begin{defn}[Exploded tropical function]
 The sheaf of exploded tropical functions $\mathcal E(\ex B)$ is the sheaf of $\mathbb C\e {\mathbb R}$ valued functions which are locally equal to a finite sum of exploded functions in $\mathcal E^\times(\ex B)$. (`Sum' means sum using pointwise addition in $\mathbb C\e {\mathbb R}$.) 
 
 \end{defn}
 
 
   The operation of addition is needed here to state the usual property of being a derivation. The other reason that addition was mentioned in this paper was to emphasize the links with tropical geometry. 
 
 \
    
   The inclusion $\iota:\mathbb C\longrightarrow \mathbb C\e{\mathbb R}$ defined by $\iota (c)=c\e 0$ induces an inclusion of functions
    \[\iota: C^\infty(\ex B)\longrightarrow \mathcal E(\ex B)\]
    \[\iota(f)(p):=\iota(f(p))\]

\begin{defn}[Vector field]
A smooth  vector field $v$ on an exploded manifold $\ex B$ is determined by maps
\[v:C^\infty(\ex B)\longrightarrow C^\infty(\ex B)\]  
   
  \[v:\mathcal E(\ex B)\longrightarrow \mathcal E(\ex B)\]
so that   
 \begin{enumerate}
   \item \label{d1}\[v(f+g)=v(f)+v(g)\]
   \item \label{derivation condition}\[v(fg)=v(f)g+fv(g)\]
   \item \label{d3}\[v(c\e yf)=c\e yv(f)\]
   \item\label{d4} The action of $v$ is compatible with the smooth part homomorhpism $\totl\cdot$ and the inclusion $\iota:C^{\infty}\longrightarrow \mathcal E$ in the sense that 
   \[v(\iota f)=\iota v(f)\]
and where defined, 
\[v\totl f=\totl{vf}\]

   \end{enumerate}
   
      Smooth exploded vector fields form a sheaf. The action of the restriction of $v$ to $U$ on the restriction of $f$ to $U$ is the restriction to $U$ of the action of $v$ on $f$. 
 
  We can restrict a vector field $v$ to  a point $p\in\ex B$ to obtain a tangent vector $v_{p}$ at that point. This is determined by  maps 
   \[v_{p}:C^\infty(\ex B)\longrightarrow \mathbb R\]
   \[v_{p}:\mathcal E(\ex B)\longrightarrow \mathbb C\e{\mathbb R}\]
    (where in the above maps, $\mathbb R$ and $\mathbb C\e{\mathbb R}$ indicate the corresponding sheaves supported at $p\in \ex B$) satisfying the above conditions  with condition \ref{derivation condition} replaced by
    \[v_{p}(fg)=v_{p}(f)g(p)+f(p)v_{p}(g)\]
Denote by $T_p\ex B$ the vector space of tangent vectors at $p\in\ex B$.
  \end{defn}

  We can add vector fields on $\ex B$ and multiply them by functions in $C^{\infty}(\ex B)$.   We shall now work towards a concrete description of vector fields in coordinate charts, and show that the sheaf of smooth vector fields on an exploded manifold $\ex B$ is equal to the sheaf of smooth sections of $T\ex B$, which is a real vector bundle over $\ex B$. In local coordinates $\{x_i,\tilde z_j\}$, a basis for this vector bundle will be given by $\{\frac \partial {\partial x_i}\}$ and the real and imaginary parts of $\{\tilde z_i\frac\partial{\partial\tilde z_i}\}$. This is part of the reason that the dimension of $\mathbb R^{n}\times \et mP$ is $n+2m$.

  \begin{lemma}\label{differentiation order}
  Differentiation does not change the order of a function in the sense that given any smooth  vector field $v$ and exploded tropical function $f\in\mathcal E$, 
  \[\totb {f}=\totb {vf}\]
  \end{lemma}
  
  \pf 
  
  As $v$ is a derivation, $v(1^{2})=1v(1)+1v(1)$, so $v(1)=0$.
  Using axiom \ref{d4} gives
  \[v\e 0=v(\iota 1)=\iota(v1)=0\e 0\]
  We can apply $v$ to the equation 
   $f=1\e 0f$, so $vf=0\e 0f+1\e 0vf$. Taking the tropical part of this equation gives 
  \[\totb {vf}=\totb{ f}+\totb {vf}\text{ i.e. }\totb {vf}\geq\totb{f}\]
  (Recall that we use the order $\e x<\e y$ if $x>y$ as we are thinking of $\ex t$ as being tiny. So $\e x+\e y=\e x$ means that $\e x\geq\e y$.)

   Now suppose that $f\in \mathcal E^{\times}(\ex B)$. Then
   \[0\e 0=v\lrb{\frac f f}=\frac {vf}{f}+fv\lrb{\frac{1}f}\]
   Therefore, 
   \[\e 0\geq \totb{vf}/{\totb f}\]
   so 
   \[\totb f\geq \totb {vf}\]
   
Therefore for any $f\in\mathcal E^{\times}(\ex B)$, $\totb{vf}=\totb f$. As any exploded tropical function is locally  a finite sum of such invertible functions, we may use axiom \ref{d1} to see that the same equation holds for any exploded tropical function.  
  \stop
  
  \begin{lemma}
  For any smooth exploded manifold $\ex B$, there exists a smooth exploded manifold $T\ex B$, the tangent space of $\ex B$, along with a canonical smooth projection $\pi: T\ex B\longrightarrow \ex B$ that makes  $T\ex B$ into a real vector bundle over $\ex B$. The sheaf of smooth vector fields on $\ex B$ is equal to the sheaf of smooth sections of this vector bundle.
  
  In particular, \[T(\mathbb R^n\times\et mP)=\mathbb R^{2n+2m}\times\et mP\]
  \end{lemma}
  
  \pf
  
  We shall first prove that $T(\mathbb R^n\times\et mP)=\mathbb R^{2n+2m}\times\et mP$. We shall use coordinate functions $x_i$ for $\mathbb R^n$ and $\tilde z_i$ for $\et mP$.  A section of $\mathbb R^{2n+2m}\times\et mP\longrightarrow \mathbb R^n\times\et mP$ is given by $n+2m$ smooth functions on $\mathbb R^n\times \et mP$. To a vector field $v$, associate the $n$ smooth functions $v(x_i)$, and the $m$ smooth complex valued functions given by $\totl{ v(\tilde z_i)\tilde z_i^{-1}}$. (Lemma \ref{differentiation order} tells us 
 that $\totb{\tilde z_i^{-1}v(\tilde z_i)}=\e 0$, so this makes sense.) The functions $(v(x_{i}),\totl{\tilde z_{i}^{-1}v(\tilde z_{i})})$ give us $n+2m$ smooth real valued functions, and therefore give us a section to associate with our vector field.

Now we must show that given an  arbitrary choice of $n+2m$ smooth functions on $\et mP$, there exists a smooth vector field $v$ so that  the functions are  $v(x_i)$ and the real and imaginary parts of  $\totl{ v(\tilde z_i)\tilde z_i^{-1}}$, and we must show that these $n+2m$ functions uniquely determine our vector field.

First, recall that exploded tropical functions are a sum of functions of the form $f\e y \tilde z^\alpha$, where $f$ is some smooth function of $x$ and $\zeta_{j}=\totl{\e{a_{j}}\tilde z^{\alpha^j}}$. Axioms \ref{derivation condition} and \ref{d3} imply that $v(\e y\tilde z^{\alpha^{j}})=\e y\tilde z^{\alpha^{j}}\sum_{i}\totl{v(\tilde z_{i})\tilde z_{i}^{-1}}\alpha^{j}_{i}$. Then axiom \ref{d4} implies that  \[v(\zeta_{j})=\zeta_{j}\sum_{i}\totl{v(\tilde z_{i})\tilde z_{i}^{-1}}\alpha^{j}_{i}\]
 Use the notation
\[\zeta_{j}=\totl{\e{a_{j}}\tilde z^{\alpha^j}}:=e^{t_{\alpha^{j}}+i\theta_{\alpha^{j}}}\]
We have that 
\[v(t_{\alpha^{j}})=\sum_{i}\Re \left(v(\tilde z_i)\tilde z_i^{-1}\right)\alpha_i^j\]
\[v(\theta_{\alpha^{j}})=\sum_{i}\Im\left(v(\tilde z_i)\tilde z_i^{-1}\right)\alpha_i^j\]
As $v$ is a derivation on smooth functions, we may apply the usual rules of differentiation including the chain rule to determine what $v(f)$ should be. 
\[v(f)=\sum v(x_i)\frac {\partial f}{\partial x_i}+\sum_{i,j}\Re \left(v(\tilde z_i)\tilde z_i^{-1}\right)\alpha_i^j\frac {\partial f}{\partial t_{\alpha^j}}+\Im\left(v(\tilde z_i)\tilde z_i^{-1}\right)\alpha_i^j\frac {\partial f}{\partial \theta_{\alpha^j}}\]

Note that this depends only $f$ as a function on $\mathbb R^{n}\times\et mP$, despite the fact that $\frac {\partial f}{\partial t_{\alpha^j}}$ and $\frac {\partial f}{\partial \theta_{\alpha^j}}$ may depend on the extension of $f$ to a smooth function of $\totl{\e{a_{j}}\tilde z^{\alpha^j}}$ and $x$, because $\sum_{j}\alpha_{i}^{j}\frac\partial{\partial t_{\alpha^{j}}}$ and $\sum_{j}\alpha_{i}^{j}\frac\partial{\partial \theta_{\alpha^{j}}}$ are vector fields tangent to the subset where $f$ is defined before extension.  Putting any smooth functions in the above formula in the place of $v(x_{i})$ and the real and imaginary parts of $\totl{v(\tilde z_{i})\tilde z_{i}^{-1}}$ gives a derivation. Note also that $v(f)$ is a smooth function, and is real valued if $f$ is real valued.  Using axioms \ref{d1} and \ref{derivation condition},  the corresponding formula for an exploded function is

\[v\left(\sum_{\alpha}f_\alpha\e {y_\alpha} \tilde z^\alpha\right):= \sum_{\alpha}\left(v(f_\alpha)+f_\alpha\sum_{i}\alpha_i v(\tilde z_i)\tilde z_i^{-1}\right)\e {y_\alpha}\tilde z^\alpha\]

It can be shown that $v$ satisfying such a formula satisfies all the axioms for being a smooth exploded vector field, is well defined, and is zero if and only if $v(x_i)=0$ and $\fun {v(\tilde z_i)\tilde z_i^{-1}}=0$.

 This shows that $T\lrb{\mathbb R^n\times\et mP}=\mathbb R^{2n+2m}\times \et mP$. Then the fact that $T\ex B$ is a vector bundle over $\ex B$ follows from our coordinate free definition of a smooth vector field, and the fact that every smooth exploded manifold is locally modeled on $\mathbb R^{n}\times \et mP$.

\stop

For any smooth exploded manifold $\ex B$, we now have that $T\ex B$ is a real vector bundle. In local coordinates $\{x_i,\tilde z_j\}$, a basis for this vector bundle is given by $\{\frac \partial {\partial x_i}\}$ and the real and imaginary parts of $\{\tilde z_i\frac\partial{\partial\tilde z_i}\}$. The dual of the tangent bundle,  $T^*\ex B$ is the cotangent space. A basis for the cotangent space is locally given by $\{dx_i\}$ and the real and imaginary parts of $\{\tilde z_i^{-1}d\tilde z_i\}$. We can take tensor powers (over smooth real valued functions) of these vector bundles, to define the usual objects found on smooth manifolds. 

\begin{remark}[Metrics]\end{remark} A metric on $\ex B$ is a smooth, symmetric, positive definite section of $T^*\ex B\otimes T^*\ex B$. Note that the inverse image of any point $p$ in $\totb{\ex B}$ has the structure of a $(\mathbb C^{*})^{m}$ bundle $M_{p}$ over some manifold in the sense that a  smooth map $\mathbb R\longrightarrow \ex B$ is equivalent to a choice of $p\in\totb{\ex B}$ and a smooth map $\mathbb R\longrightarrow M_{p}$. Any smooth metric on $\ex B$ gives a $(\mathbb C^{*})^{m}$ invariant metric on $M_{p}$ which is complete if $\ex B$ is compact. We may therefore carry out any local construction familiar from Riemannian geometry. The topology coming from a metric on $\ex B$ is Hausdorff, and with this topology, we may regard $\ex B$ as the disjoint union of $M_{p}$ for all $p$ in $\totb{\ex B}$.

\begin{defn}[Standard metric and basis for tangent space]\label{standard metric} Define the standard basis for $T(\mathbb R^{n}\times \et mP)$ to be the basis given by the vector fields $\frac\partial{\partial x_{i}}$, and the real and imaginary parts of $\{\tilde z_i\frac\partial{\partial\tilde z_i}\}$. Define the standard metric on $\mathbb R^{n}\times \et mP$ to be the metric in which the standard basis is orthonormal, and let the standard connection on $\mathbb R^{n}\times \et mP$ be the connection which preserves the standard basis (the Levi-Civita connection of the standard metric.)
\end{defn}

\begin{defn}[Integral vector]\label{integral}
 An integral vector $v$ at a point $p\longrightarrow\ex B$ is a vector $v\in T_p\ex B$ so that for any exploded function $f\in \mathcal E^\times(\ex B)$, 
 \[v(f)f^{-1}\in\mathbb Z\]
 Use the notation ${}^{\mathbb Z}T_p\ex B\subset T_p\ex B$ to denote the integral vectors at $p\longrightarrow \ex B$.
 
\end{defn}

For example, a basis for ${}^{\mathbb Z}T\ex T^n$ is given by the real parts of $ \tilde z_i\frac \partial{\partial \tilde z_i}$. The only integral vector on a smooth manifold is the zero vector.

\

Given a smooth morphism  $f:\ex B\longrightarrow \ex C$, there is a natural smooth morphism $df:T\ex B\longrightarrow T\ex C$ which is the differential of $f$, defined as usual by \[df(v)g:=v( g\circ f)\] Of course, $df$ takes integral vectors to integral vectors.

\

As usual, the flow of a smooth vector field gives a smooth morphism (with the usual caveats about existence in the noncompact case - the existence theory for flows of vector fields on exploded manifolds is identical to the existence theory on smooth manifolds.)

\begin{thm} \label{smooth flow} If $v$ is a smooth vector field on an exploded manifold $\ex B$, then the flow of $v$ for time $1$ is a smooth map when it exists.
\end{thm}

\pf

First consider the case of a coordinate chart of the form $\et nn:=\lrb{\et 11}^{n}$. We can consider the subset of $\et nn$ over any given point in the tropical part $\totb{\et nn}$ to be a smooth manifold, and the restriction of any smooth vector field on $\et nn$ to this subset is just a smooth vector field. Therefore we can apply the usual existence, uniqueness and regularity results in this context. We shall assume that the time $1$ flow of our vector field exists.

We must show that the time one flow of $v$ composed with any smooth exploded function $g\e a\tilde z^{\alpha}$ is  still a smooth exploded function. Use the notation $g(t,\tilde z) \e a\tilde z^{\alpha}$ to indicate the above exploded function composed with the time $t$ flow of $v$. Then 
\[\frac {\partial g}{\partial t}=v(g)+g\tilde z^{-\alpha}v(\tilde z^{\alpha})\]
Note that if $v$ is a smooth vector field, $\tilde z^{-\alpha}v(\tilde z^{\alpha})$ is a smooth function. We shall now prove the smoothness of $g(t,\tilde z)$ using the fact that the flow of smooth vector fields on smooth manifolds is smooth.
 
In particular, the smooth part of $\et nn$ is just $\mathbb C^{n}$.
 A basis for the tangent space is the real and imaginary parts of $\tilde z_{i}\frac \partial{\partial \tilde z_{i}}$, which correspond to the smooth vector fields $r_{i}\frac \partial{\partial r_{i}}$ and $\frac \partial {\partial{\theta_{i}}}$ in polar coordinates on $\mathbb C^{n}$. As any smooth vector field on $\et nn$ is a sum of smooth functions times the above basis vector fields, any smooth vector field on $\et nn$ corresponds to a smooth vector field on $\mathbb C^{n}$ (tangent to all the coordinate planes.) 
 
 Consider the function $\hat g$ on $\mathbb R\times \mathbb C^{n+1}$ which is equal to $g(t,z_{1},\dotsc,z_{n})z_{n+1}$. $\hat g$ satisfies the differential equation 
 \[\frac {\partial\hat  g}{\partial t}=v(\hat g)+\hat g\tilde z^{-\alpha}v(\tilde z^{\alpha})=\hat v\hat g\]
 where
 \[\hat v=v+x\]
where $x$ is a vector field pointing in the last coordinate direction for which $xz_{n+1}=\tilde z^{-\alpha}v(\tilde z^{\alpha})z_{n+1}$. As $\hat g$ at time $1$ is the composition of $\hat g$ at time $0$ with the time $1$ flow of the smooth vector field $\hat v$, $\hat g$ at time $1$ is smooth, therefore, its restriction to $\{z_{n+1}=1\}$ is also smooth, therefore $g$ at time $1$ is smooth. This proves that smooth exploded functions composed with the time $1$ flow of $v$ are smooth.

The general case now follows quickly. The argument for $\et nn\times \ex T^{m}$ is analogous. As any other coordinate chart is the restriction to a subspace of $\et nn\times \ex T^{m}$, the case of a smooth vector field on a general coordinate chart follows from the observation that we can extend it to a smooth vector field on $\et nn\times \ex T^{m}$. Regularity results for flows which only locally exist in a coordinate chart can as usual be obtained from the above using smooth cutoff functions, and the global result for an exploded manifold follows.

\stop

\

Theorem \ref{smooth flow} allows us to prove that any open subset of an exploded manifold is an exploded manifold.

\begin{lemma}\label{sub coordinates}
Any open subset $U\subset\ex B$ of an exploded manifold  is an exploded manifold. In particular, given any point $p\in U$, there exists an open neighborhood $U'$ of $p$ contained inside $U$ which is isomorphic to $\mathbb R^{n}\times \et mP$. 

\end{lemma}

\pf
As any exploded manifold is locally isomorphic to $\mathbb R^{n}\times\et mP $, we may restrict to the case that $U$ is an open subset of $\mathbb R^{n}\times \et mP$. Recall from example \ref{polytope} on page \pageref{polytope}  that the topology on $\mathbb R^{n}\times \et mP$ can be described as follows: Choose a basis $\{\zeta_{1},\dotsc,\zeta_{m'}\}$ for the set of smooth monomials on $\et mP$. 
Then any open subset $U\subset {\mathbb R^{n}\times \et mP}$ is of the form
\[U:=\{(x,{ \zeta_{1}},\dotsc,{ \zeta_{m'}})\in V \}\] 
where $V\subset\mathbb R^{n}\times \mathbb C^{m'}$ is an open subset, and smooth exploded functions are of the form
\[f(x, \zeta_{1},\dotsc, \zeta_{m'})\e a \tilde z^{\alpha}\text{ for }f\in C^{\infty}(\mathbb R^{n}\times \mathbb C^{m'},\mathbb C^{*})\]
We shall use the fact that $\et mP$ is isomorphic to a  subset of $\et mP$ where all the $ \zeta_{i}$ are small. An isomorphism can be constructed as follows: Use coordinates $\tilde z_{k}=e^{r_{k}+i\theta_{k}}\e {a_{k}}$, and consider the vector field $v$ given by 
\[v=\sum_{j}\sum_{k=1}^{m} \abs{ \zeta_{j}}^{2}\alpha^{j}_{k}\frac \partial{\partial r_{k}} \]
where $\zeta_{j}=\totl{\e a\tilde z^{\alpha^{j}}}$.
This  vector field $v$ is half the gradient of the the smooth function \[W:=\sum_{j}\abs{ \zeta_{j}}^{2}\] using a metric where $\{\frac \partial {\partial r_{k}},\frac\partial{\partial \theta_{k}}\}$ are an orthonormal basis. For any point $p\in \et mP$, if $\totb p$ is in the interior of $P$, then $ \zeta_{j}(p)=0$ for all $j$. If $\totb p$ is not in the interior of $P$, given any vector $r_{p}$ pointing into the interior of $P$, if $ \zeta_{j}(p)\neq 0$, then $r_{p}\cdot \alpha^{j}>0$ because the tropical part of the monomial used to define $\zeta_{j}$ is positive on the interior and $0$ at $\totb p$. It follows that  $\sum_{k} (r_{p})_{k}\frac{\partial W}{\partial r_{k}}>0$ where $W>0$. Therefore, the smooth gradient vector field $v$ defined above is non vanishing where  $W\neq 0$ and $vW>0$ when $W>0$. Therefore, by multiplying $v$ by a smooth function $f$, we may achieve the following:

\[fvW:=\begin{cases}0\text{ where } W<\frac \epsilon 4 
\\ W^{2} \text{ where }W\geq \frac \epsilon 2\end{cases}\]
Theorem \ref{smooth flow} implies that the  flow of $fv$ for time $\frac 1\epsilon$ is an isomorphism from the the set where $W<\epsilon$ to $\et mP$.  This completes the proof of the claim that $\et mP$ is isomorphic to a subset where the $\zeta_{j}$ are small.

 Suppose that at the point $p$, $ \zeta_{i}=0$ for $i\in I$ and $ \zeta_{j}\neq 0$ for $j\notin I$. Using the notation $\zeta_{j}=\totl{\tilde \zeta_{j}}$ where $\tilde \zeta_{j}=\e a\tilde z^{\alpha^{j}}$, consider the face $F\subset P$ defined to be the set where $\totb{\tilde \zeta_{j}}=\e 0$ for all $j\notin I$. This face $F$ is the face which contains $\totb p$ in its interior. By changing coordinates, we may assume that 
\begin{itemize}\item
$F$ is the intersection of $P$ with the coordinate plane on which the first $l$ coordinates of  $\mathbb R^{m}$ vanish, where the dimension of $F$ is $m-l$
\item  $P$ is contained in the quadrant of $\mathbb R^{m}$ where the first $l$ coordinates are nonnegative.
 \end{itemize}

A neighborhood of $p$ contained in $U$ will be isomorphic to $\mathbb R^{n+2l}\times \et {m-l}F$. For now, identify coordinates on $\et {m-l}F$ with the last $m-l$ coordinates on $\et mP$, and choose a generating set $\{ \zeta'_{i}\}$ for the smooth monomials on $\et {m-l} F$. 

Let $\phi_{1}$ be a diffeomorphism of $\mathbb R^{n+2l}$ onto a small ball in $\mathbb R^{n}\times (\mathbb C^{*})^{l}$, and let $\phi_{2}$ be an isomorphism of $\et {m-l}F$ onto an open subset of $\et {m-l}F$ where $\sum_{i}\abs{ \zeta'_{i}}^{2}$ is small. Then by identifying coordinates on $\et {m-l}F$ with the the last $(m-l)$ coordinates on $\mathbb R^{n}\times \et mP$, and identifying the remaining coordinates on the subset of $\mathbb R^{n}\times \et mP$ with tropical part $F$ with $\mathbb R^{n}\times (\mathbb C^{*})^{l}$ ( if $\tilde z_{i}$ for $i\leq l$ is the $i$th coordinate on $\et mP$, use $\totl{\tilde z_{i}}$ as the $i$th coordinate of $(\mathbb C^{*})^{l}$), we can combine these two isomorphisms into a smooth map 
\[\phi:=\phi_{1}\times\phi_{2}:\mathbb R^{n+2l}\times \et {m-l}F\longrightarrow \mathbb R^{n}\times \et mP\]
The image of $\phi$ is an open subset of $\mathbb R^{n}\times\et mP$ which we can choose to contain our point $p$, and be contained in $U$. 

It remains to show that the inverse map, $\phi^{-1}$ is also smooth. To check this, we shall check that $\phi^{-1}$ composed with any smooth exploded function on $\mathbb R^{ n+2l}\times\et {m-l}F$ is a smooth exploded function on our open subset of $\mathbb R^{n}\times\et mP$. The composition of $\phi^{-1}$ with any smooth exploded function is a monomial in $\tilde z$ times a smooth $\mathbb C^{*}$ valued function of the $ \zeta_{i}'$ and of $\phi^{-1}$ composed with coordinates on $\mathbb R^{n+2l}$. It therefore suffices to check that $\zeta_{i}'$ is smooth and $\phi^{-1}$ composed with any coordinate on $\mathbb R^{n+2l}$ is smooth. Let $\tilde z_{i}$ be one of the first $l$ coordinates of $\et mP$. As $\totb{\tilde z_{i}}$ is $\leq \e0$ on $\et mP$, $\totl{\tilde z_{i}}$ is smooth. Therefore, any smooth function of $\mathbb R^{n}\times(\mathbb C^{*})^{l}$ considered as a smooth function of the $\mathbb R^{n}$ times the first $l$ coordinates of $\et mP$ is smooth. It follows that $\phi^{-1}$ composed with any coordinate function of $\mathbb R^{n+2l}$ is smooth.  Using the notation $ \zeta'_{i}=\totl{\tilde \zeta'_{i}}$, we can multiply  $\tilde \zeta_{i}'$ by a product $\tilde z^{\alpha}$ of powers of the first $l$ coordinates on $\et mP$ so that $\totb{\tilde \zeta_{i}'\tilde z^{\alpha}}\leq\e 0$ on $P$. Therefore, $\totl{\tilde \zeta_{i}'\tilde z^{\alpha}}$ is a smooth function on ${\et mP}$. Restricted to our subset (where $\totb{\tilde z^{\alpha}}=0$), $\totl{\tilde z^{\alpha}}$ is a smooth $\mathbb C^{*}$ valued function, so $\totl{\tilde z^{-\alpha}}$ is also smooth, and therefore, $ \zeta_{i}'=\totl{\tilde \zeta_{i}'\tilde z^{\alpha}}\totl{\tilde z^{-\alpha}}$ is smooth.  Therefore $\phi^{-1}$ is smooth and $\phi$ is our required isomorphism.

\stop

\section{$C^{k,\delta}$ regularity}\label{regularity}

This section defines some regularities that are natural to consider on exploded manifolds.  These extra regularities are needed because they are the natural level of regularity of the moduli stack of holomorphic curves in an almost complex exploded manifold (see \cite{reg} and \cite{egw} for details). The reader wishing a simple introduction to exploded manifolds may skip this somewhat technical section on first reading.
   
    For any choice of  smooth metric on the exploded manifold $\et 11$, the subset where  $\totb{\tilde z}=\e 0$ will be the manifold $\mathbb C\setminus\{0\}$ with some metric that has a cylindrical end at $0$. When faced with a manifold with a cylindrical end, one way of defining a class of functions with nice regularity is to `compactify' that cylindrical end, and consider the class of smooth functions on the resulting compactified manifold. This is how we defined `smooth' functions on $\et 11$. This choice of what is `smooth' was chosen simply because it was easy to describe using existing language - it is not the only natural choice.  A $\C\infty\delta$ function is a generalization of a function on a manifold with a cylindrical end which is smooth on the interior, and which decays exponentially along with all its derivatives on the cylindrical end.  

\

Recall that every exploded manifold $\ex B$ is a (possibly non-Hausdorff) topological space, so we may talk of continuous maps from $\ex B$ to any topological space. 

     \begin{defn}[Continuous morphism]
       A $C^0$ exploded function $f\in \mathcal E^{0,\times}(\ex B)$ is a function of the form $fg$ where $f$ is a continuous map from $\ex B$ to  $\mathbb C^{*}$ and $g\in\tC{\ex B}$ is a smooth exploded function. We can define $C^{0}$ morphisms of exploded manifolds to be morphisms of abstract exploded manifolds using the sheaf $\mathcal E^{0,\times}$ instead of $\mathcal E^{\times}$.
      \end{defn}
      
      When referring to a continuous map $\ex B\longrightarrow \ex C$ we shall always mean a $C^{0}$ morphism $\ex B\longrightarrow \ex C$ in the sense of the above definition, which is  stronger than a continuous map from $\ex B$ to $\ex C$ considered as topological spaces. For example, any map to $\ex T$ as a topological space is continuous because $\totl{\ex T}$ is a point, but the only continuous morphisms from $\ex T$ to itself are morphisms of the form $c\e a\tilde z^{\alpha}$. On the other hand, if $M$ is a manifold, a $C^{0}$ morphism $\ex B\longrightarrow M$ is equivalent to a continuous map from $\ex B$ to $M$ as topological spaces.
      
      \

Given a continuous real valued function $f$ on $\et 11$, by saying `$f$ converges exponentially with weight $\delta$' on $\et 11$, we mean that the function
\[(f(\totl{\tilde z})-f(0))\abs {\totl{\tilde z}}^{-\delta}\]
extends to a continuous function on all of $\et 11$ which is zero when $\totl{\tilde z}=0$. The class  of $C^{\infty,\delta}$ functions on $\et 11$ consists of functions which are continuous, and have continuous derivatives to all orders which converge exponentially with weight $\delta$ on $\et 11$. It is this that we must generalize to $\mathbb R^{n}\times \et mP$.

\begin{defn}[The operator $e_{S}$]
Given any real or vector valued $C^0$ function  $f$ on $\mathbb R^n\times \et mP$, and a stratum $S\subset P$ define 
\[e_S(f)(x,\tilde z_{1},\dotsc,\tilde z_{m}):=f(x,\tilde z_{1}\e {\frac {a_{1}}2}\totb{\tilde z_{1}}^{-\frac 12},\dotsc,
\tilde z_{m}\e {\frac {a_{m}}2}\totb{\tilde z_{m}}^{-\frac 12})\]
where $(\e{a_1},\dotsc,\e{a_m})$ is any point  in $S$, and $(\e c)^{-\frac 12}=\e{-\frac c2}$. 
%

\end{defn}

So $e_{S}f(x,\tilde z)$ samples the function $f$ at a point with tropical part half way between $\totb{\tilde z}$ and the point $\e a$ in $S$. Note that $e_{S}f$ does not depend on the choice of the point $\e a$ in $S$. 

\

For example consider $\et 22:=\et 2{[0,\infty)^{2}}$. The polytope $[0,\infty)^{2}$ has two one dimensional strata 
\[S_{1}:=(0,\infty)\times 0 \ \ \ \ \ S_{2}:= 0\times(0,\infty)\] and one two dimensional stratum $S_{3}:=(0,\infty)^{2}$.
 If we have a function $f\in C^0(\et 22)$ (in other words a continuous map from the topological space $\et 22$ to $\mathbb R$), then

 \[e_{S_{1}}f(z_{1},z_{2})=f(0,z_{2})\ \ \ \ \ \ \ e_{S_{2}}f(z_{1},z_{2})=f(z_{1},0) \ \ \ \ \ e_{S_{3}}f(z_{1},z_{2})=f(0,0)\]

As a second example, consider $\et1{[0,1]}$. Smooth or continuous functions on $\et 1{[0,1]}$ are generated by $\zeta_{1}=\totl{\tilde z}$ and $\zeta_{2}=\totl{\e{1}\tilde z^{-1}}$. There are three strata of  $[0,1]$ to consider: $0$, $1$, and $(0,1)$.
\[e_{0}\zeta_{1}=\zeta_{1}\ \ \ \  e_{0}\zeta_{2}=0\ \ \ \ e_{0}f(\zeta_{1},\zeta_{2})=f(\zeta_{1},0)\]
\[e_{1}\zeta_{1}=0\ \ \ \ e_{1}\zeta_{2}=\zeta_{2}\ \ \ \ e_{1}f(\zeta_{1},\zeta_{2})=f(0,\zeta_{2})\]
\[e_{(0,1)}\zeta_{1}=0\ \ \ \ e_{(0,1)}\zeta_{2}=0\ \ \ \ e_{(0,1)}f(\zeta_{1},\zeta_{2})=f(0,0)\]
Note that we can consider $\et 1{[0,1]}$ as the subset of $\et 22$ where $\tilde z_{1}\tilde z_{2}=\e 1$. From this perspective we can relate the above two examples by  $e_{0}=e_{S_{1}}$ and $e_{1}=e_{S_{2}}$.

\

 In general, the smooth or continuous functions on $\et mP$ are generated by functions $\zeta_{i}$ of the form $\totl{\e {a_{i}}\tilde z^{\alpha^{i}}}:=\totl{\tilde \zeta_{i}}$. For any stratum $S\subset P$ one of the following two options hold:
 \begin{enumerate}
 \item $e_{S}\zeta_{i}=0$, $\zeta_{i}$ vanishes on the stratum of $\et mP$ corresponding to $S$ and $\totb{\tilde \zeta_{i}}<\e 0$ on $S$
 
 \item or $e_{S}\zeta_{i}=\zeta_{i}$ and  $\zeta_{i}$ is nowhere $0$ on the stratum of $\et mP$ corresponding to S, and $\totb{\tilde \zeta_{i}}=\e0$ on $S$.
 \end{enumerate}
 
  The operation $e_{S}$ on a continuous function $f$ on $\et mP$ is then given by \[e_{S}f(\zeta_{1},\dotsc, \zeta_{n})=f(e_{S}\zeta_{1},\dotsc, e_{S}\zeta_{n})\] Of course, this implies that if $f$ is smooth or continuous,  $e_{S}f$ is too.

\

 Note that the operations $e_{S_{i}}$ commute and $e_{S_{i}}e_{S_{i}}=e_{S_{i}}$. More generally, $e_{S_{i}}e_{S_{j}}=e_{S'}$ where $S'$ is the smallest stratum of $P$ whose closure contains both $S_{i}$ and $S_{j}$.
 
\begin{defn}[The operators $e_{I}$ and $\Delta_{I}$]If $I$ denotes any collection of strata $\{S_{1},\dotsc, S_{n}\}$ of $P$, we shall use the notation 
\[e_{I}f:=e_{S_{1}}\lrb{e_{S_{2}}\lrb{\dotsb e_{S_{n}}f}}\]

\[\Delta_{I}f:=\lrb{\prod_{S_{i}\in I}(\id-e_{S_{i}})}f\]
\end{defn}
For example on $\et 22$,

\[\begin{split}\Delta_{S_{1},S_{2}}f(z_{1},z_{2})&:=(1-e_{S_{1}})(1-e_{S_{2}})(f)(z_1,z_2)
\\&:=f(z_1,z_2)-f(0,z_2)-f(z_1,0)+f(0,0)\end{split}\]
Note that if $S\in I$, $e_{S}\Delta_{I}=0$. In the above example, this corresponds to $\Delta_{S_{1},S_{2}}f(z_{1},0)=0$ and $\Delta_{S_{1},S_{2}}f(0,z_{2})=0$.

\

The operator $(\id-\Delta_{I})$ gives a nice way to extend the domain of definition of a function $f$ defined only on the closure in $\totl{\et mP}$ of the strata in $I$. The function $g=(\id-\Delta_{I})f$ is defined on all of $\et mP$, is smooth if $f$ is, and $e_{S}g=e_{S}f$ for all $S\in I$. For example, if $f( z_{1}, z_{2})$ is a smooth real function defined on the subset of $\et 22$ for which ${ z_{1} z_{2}}=0$, then $(\id-\Delta_{S_{1},S_{2}})f(z_{1},z_{2}):=f(0,z_{2})+f(z_{1},0)-f(0,0)$ is a smooth function extending the domain of definition of $f$ to all of $\et 22$. 

\ 

We shall need a weight function $w_{I}$ for every collection of nonzero strata $I$. (We shall need this weight function to measure how fast functions `converge' when approaching the strata in $I$.) This will have the property that if $f$ is any smooth function, then $\Delta_{I}f$ will be bounded by a constant times $w_{I}$ on any compact subset of $\mathbb R^{n}\times \et mP$. 

Consider the set $Z_{I}$ of smooth monomials on $\mathbb R^{n}\times\et mP$ of the form $\zeta =\totl{\e a\tilde z^{\alpha}}$ so that $\Delta_{I}\zeta=\zeta$. (This is equivalent to $e_{S}\zeta=0$ for all $S\in I$.) Choose some finite set $\{\zeta_{i}\}$ of generators  for $Z_{I}$ so that any $ \zeta\in Z_{I}$  is a product of one of these $\zeta_{i}$ with another smooth function. Then define 

\[w_{I}:=\sum\abs{\zeta_{i}}\]

Continuing the example of $\et 22$ started above, we can choose  $w_{S_{1}}=\abs{z_{1}}$, $w_{S_{2}}=\abs{z_{2}}$,  
$w_{S_{1},S_{2}}=\abs{z_{1}z_{2}}$, and $w_{S_{3}}=\abs {z_{1}}+\abs{z_{2}}$.

\

Note that for any $\zeta \in Z_{I}$, the size of  $\zeta$ is bounded by a constant times $w_{I}$ on any compact subset. Therefore, given any other choice of generators for $Z_{I}$,  the resulting $w'_{I}$ is bounded by a constant times $w_{I}$ on any  compact subset of $\mathbb R^{n}\times \et mP$. Note also that $w_{I_{1}}w_{I_{2}}$ is bounded by a constant times $w_{I_{1}\cup I_{2}}$ on any  compact subset of $\mathbb R^{n}\times \et mP$, as $w_{I_{1}}w_{I_{2}}$ is a finite sum of absolute values of $\zeta\in Z_{I_{1}\cup I_{2}}$.

\

As mentioned above, these $w_{I}$ have the property that if $f$ is any smooth function, then $\Delta_{I}f$ is bounded by a constant times $w_{I}$ on any compact subset of $\mathbb R^{n}\times \et mP$.

 For example, consider the case of $\et nn:=\et{n}{[0,\infty)^{n}}$. Any smooth function on $\et nn$ is determined by a smooth function $f$ on $\totl{\et nn}=\mathbb C^{n}$. Consider the stratum $S_{i}$ corresponding to the set where $z_{i}:=\totl{\tilde z_{i}}=0$. Then, as $\Delta_{S_{i}}f$ has a continuous derivative and vanishes when $z_{i}=0$, $\Delta_{S_{i}}f$ is bounded by a constant times $\abs {z_{i}}$ on compact subsets of $\mathbb C^{n}$. If $i\neq j$, $\Delta_{S_{i}, S_{j}}f$ vanishes where $z_{i}z_{j}=0$. Therefore, as $\Delta_{S_{i},S_{j}}f$ has a continuous derivative, it is bounded by a constant times $\abs{z_{i}z_{j}}$ on compact subsets of $\mathbb C^{n}$ away from where $z_{i}=z_{j}=0$. On the set where $z_{i}=z_{j}=0$, $\Delta_{S_{i},S_{j}}f$ vanishes and its derivative vanishes, therefore, as $\Delta_{S_{i},S_{j}}f$ has continuous second derivative, it is bounded by a constant times $\abs{z_{i}z_{j}}$. Similarly, if $I=\{S_{i_{1}},\dotsc, S_{i_{k}}\}$, $\Delta_{I}f$ has continuous derivatives up to order $k$ and vanishes on the set where $z_{i_{1}}\dotsb z_{i_{k}}=0$, so $\Delta_{I}f$ is bounded by a constant times $\abs{z_{i_{1}}\dotsb z_{i_{k}}}$ on compact subsets of $\et nn$. In the lemma below, we prove the general case.

\begin{lemma}\label{smooth ckdelta}Given any $C^{k}$ function $f$ on $\mathbb R^{n}\times \et mP$, and collection $I$ of at most $k$ strata of $P$,  the function  $\Delta_{I}f$ is bounded by a constant times $w_{I}$ on compact subsets of $\mathbb R^{n}\times \et mP$.
\end{lemma}

\pf

We shall first consider the case of $\et nn$. Use coordinates $z$ for $\totl{\et nn}$.  Introduce a real variable $t_{S}$ for each $S\in I$, and let $t$ denote the vector of all these variables. Define
\[\phi_{I}(t,z):=\lrb{\prod_{S\in I}(e_{S}+t_{S}\Delta_{S})}z\]
This function has the property that if $t_{S}=0$ for all $S\in I_{1}$ and $t_{S'}=1$ for all $S'\in I\setminus I_{1}$, then $\phi_{I}(t,z)=e_{I_{1}}z$.

Use the notation $D_{I}:=\prod_{S\in I}\frac \partial{\partial t_{S}}$. We can rewrite $\Delta_{I}f$ as follows.
\[\Delta_{I}f(z)=\int_{0}^{1}\dotsb \int_{0}^{1}D_{I}f(\phi_{I}(t,z))dt\]
To bound $\Delta_{I}f$, we shall bound the above integrand. 
\[\begin{split}D_{I}f(\phi(t,z))&=\sum_{\coprod_{i=1}^{l}I_{i}=I}D^{l}f(\phi_{I}(t,z))(D_{I_{1}}\phi_{I}(t,z))\dotsb(D_{I_{l}}\phi_{I}(t,z))
\\&=\sum_{\coprod_{i=1}^{l}I_{i}=I}D^{l}f(\phi_{I}(t,z))(\Delta_{I_{1}}\phi_{I-I_{1}}(t,z))\dotsb(\Delta_{I_{l}}\phi_{I-I_{l}}(t,z))
\end{split}\]
The above sum is over all partitions of $I$. The notation $D^{l}f$ indicates the $l$th derivative of $f$ considered as a function on $\mathbb R^{2n}$. On compact subsets, the first $k$ derivatives of $f$ are bounded by a constant. The term $\Delta_{I_{i}}\phi_{I-I_{i}}(t,z)$ is a finite sum of  monomials in $z$ which vanish on all strata in $I_{i}$ multiplied by terms dependent on $t$ which are bounded by $1$. Therefore $\Delta_{I}f$ is bounded on compact subsets by some constant times a finite sum of monomials in $z$ which vanish on all strata in $I$, which in turn are bounded by a constant times $w_{I}$. 

Our lemma therefore holds for $\et nn$. The same argument works for $\mathbb R^{n}\times \ex T^{a}\times\et bb$. Recall from Remark \ref{monomial description} that $\mathbb R^{n}\times \et mP$ is a subset of some $\mathbb R^{n}\times \ex T^{a}\times\et bb$ defined by some monomial equations in $\tilde z$. Therefore any $C^{k}$ function $f$ on $\mathbb R^{n}\times\et mP$ can be extended to a $C^{k}$ function $f'$ on $\mathbb R^{n}\times \ex T^{a}\times\et bb$, which must satisfy our lemma. Note that $e_{S}f$ is the restriction of $e_{S'}f'$ where $S'$ is the strata of $\mathbb R^{n}\times \ex T^{a}\times\et bb$ containing $S$. Therefore, $\Delta_{I}f$ is the restriction of $\Delta_{I'}f'$ where $I'$ is the corresponding  collection of strata containing the strata in $I$. It follows that $\Delta _{I}f$ is bounded on compact subsets by a constant times a sum of absolute values of monomials which vanish on all the strata in $I$, which are in turn bounded on compact subsets by $w_{I}$ times a constant.

\stop

\

We shall now define $C^{k,\delta}$ for any $0<\delta<1$:

\begin{defn}[$C^{k,\delta}$ and $\C\infty\delta$ regularity]
Define $C^{0,\delta}$ to be the same as $C^0$. A sequence of smooth functions $f_i\in C^\infty(\mathbb R^n\times\et mP)$ converge to a continuous function $f$ in  $C^{k,\delta}(\mathbb R^n\times \et mP)$ if the following conditions hold:
\begin{enumerate}
\item Given any collection $I$ of at most $k$ nonzero strata, the sequence of  functions
\[\abs{w_{I}^{-\delta}\Delta_{I}(f_i-f)}\]
converges to $0$ uniformly on compact subsets of $\totl{\mathbb R^{n}\times\et mP}$ as $i\rightarrow\infty$. (This includes the case where our collection of strata is empty and $f_{i}\rightarrow f$ uniformly on compact subsets.)
\item For any smooth vector field $v$, $v(f_i)$ converges to some function $vf$ in $C^{k-1,\delta}$.
\end{enumerate} 
Define $C^{k,\delta}(\mathbb R^n\times \et mP)$ to be the closure of $C^\infty$ in $C^0$ with this topology. Define $C^{\infty,\delta}$ to be the intersection of $C^{k,\delta}$ for all $k$. Define $\C\infty\delta$ to be the intersection of $C^{\infty,\delta'}$ for all $\delta'<\delta$.
\end{defn}

In particular, $C^{\infty}\subset\C\infty1\subset C^{k,\delta}$ for $0<\delta< 1$. Functions in $C^{k,\delta}$ can be thought of as functions which converge a little slower than $C^k$ functions when they approach different strata. Thinking of a single stratum as being analogous to a cylindrical end, this is similar to requiring exponential convergence (with exponent $\delta$) on the cylindrical end.

The $C^{k,\delta}$ topology is given by the following norm restricted to compact subsets on which the operations $e_{S}$ are still defined.
\begin{defn}[The norm $\abs\cdot_{k,\delta}$] Define $\abs{f}_{0,\delta}$ to be the the supremum of $\abs f$. Then define $\abs {f}_{k,\delta}$ to be 
\[\abs{f}_{k,\delta}:=\abs{\nabla f}_{k-1,\delta}+\sup \sum_{\abs I\leq k}\abs{w_{I}^{-\delta}\Delta_{I}f}\]
where $\nabla $ indicates the covariant derivative using the standard connection and the absolute value of tensors is measured using the standard metric (both defined on page \pageref{standard metric}). 
\end{defn}

\

\begin{lemma} $C^{k,\delta}$ is an algebra over $C^{\infty}$ for any $0<\delta<1$.

\end{lemma}

\pf

We already know that $C^{\infty}\subset C^{k,\delta}$. The sum of any two $C^{k,\delta}$ functions is clearly $C^{k,\delta}$, so it remains to prove that the product of two $C^{k,\delta}$ functions is $C^{k,\delta}$.

The following is a formula for  $\Delta_{S}$ of a product:
\[\Delta_{S}fg=\lrb{\Delta_{S}f}g+\lrb{e_{S}f}\lrb{\Delta_{S}g}\]
This formula generalizes to the case of a collection $I$ of strata as follows:
\[\Delta_{I}fg=\sum_{I'\subseteq I}\lrb{e_{I'}\Delta_{I-I'}f}\lrb{\Delta_{I'}g}\]
Therefore we can bound $\abs{w_{I}^{-\delta}\Delta_{I}fg}$ on compact subsets by expressions in $f$ and $g$ using the fact that $w_{I'}w_{I-I'}$ is bounded by a constant times $w_{I}$ as follows:
\[\abs{w_{I}^{-\delta}\Delta_{I}fg}\leq c\sum_{I'\subseteq I}\abs{w^{-\delta}_{I-I'}\Delta_{I-I'}f}\abs{w^{-\delta}_{I'}\Delta_{I'}g}\]
The above inequality is valid on compact subsets, and the constant $c$ depends on the compact subset, but is independent of $f$ and $g$. A similar inequality may be derived for derivatives of $fg$ using the product rule. It follows that $\abs{fg}_{k,\delta}$ can be bounded by a constant times $\abs f_{k,\delta}\abs g_{k,\delta}$ on compact subsets, and that the required constant is independent of $f$ and $g$.

 It follows that if  $f_{i}\rightarrow f$ and $g_{i}\rightarrow g$ in $C^{k,\delta}$, then $f_{i}g_{i}\rightarrow fg$ in $C^{k,\delta}$, because on compact subsets we can estimate 
\[\begin{split}\abs{f_{i}g_{i}-fg}_{k,\delta}\leq\abs{f(g_{i}-g)}_{k,\delta}&+\abs{g(f_{i}-f)}_{k,\delta}\\&+ \abs{(f_{i}-f)(g_{i}-g)}_{k,\delta}\end{split}\]
and  restricted to compact subsets, each of the terms in the right hand side of the above converge  to $0$ as $i\rightarrow \infty$. 

\stop

\

The following lemma allows us to define $C^{k,\delta}$ regularity without reference to smooth functions

\begin{lemma} \label{ckdelta}A continuous function $f$ on $\mathbb R^{n}\times\et mP$ is $C^{k,\delta}$  for some $k\geq 1$ if and only if the following two conditions hold.
\begin{enumerate}
\item $\nabla f$ exists and is $C^{k-1,\delta}$.
\item For all collections $I$ of at most $k$ strata, $w_{I}^{-\delta}\Delta_{I}f$ extends to a continuous function which vanishes on all strata in $I$.
\end{enumerate}
\end{lemma}

\pf

The fact that the above two conditions hold for $C^{k,\delta}$ functions follows immediately from Lemma \ref{smooth ckdelta} and the definition of $C^{k,\delta}$.
We need to show that any function $f$ satisfying the above condition can be approximated in $C^{k,\delta}$ by a smooth function.

Choose a smooth cutoff function $\rho:[0,\infty)\longrightarrow [0,1]$ so that $\rho(x)$ is $1$ for all $x\in [0,\frac 12]$ and $\rho(x)=0$ for all $x\geq 1$. Given a smooth monomial $\zeta$ on $\mathbb R^{n}\times\et mP$ in the form $\zeta=\totl{\e a\tilde z^{\alpha}}$, let $I_{\zeta}$ be the set of strata on which $\zeta$ vanishes.
Consider \[h=f-\rho(t\abs\zeta)\Delta_{I_{\zeta}}f\] for $t$ large. As  $e_{S}\Delta_{I_{\zeta}}f=0$ for any stratum $S\in I_{\zeta}$, the function $h$ agrees with $f$ when $\zeta=0$ and when $\zeta\geq \frac 1t$, but   $\Delta_{I_{\zeta}}h=0$  where $\abs\zeta\leq \frac 1{2t}$. We shall now show that $\abs{f-h}_{k,\delta}$ is small restricted to compact subsets when $t$ is large.

We may expand $\Delta_{I}(f-h)$ as follows: 
\[\Delta_{I}(f-h)=\Delta_{I}\lrb{\rho(t\abs\zeta)\Delta_{I_{\zeta}}f}=\sum_{I_{1}\coprod I_{2}=I}\lrb{\Delta_{I_{1}}\rho(t\abs\zeta)}e_{I_{1}}\Delta_{I_{\zeta}\cup I_{2}}f\]

For any stratum $S$, either $S$ is disjoint from the face on which $\zeta\neq 0$ so  $e_{S}\zeta=0$ or $S$ is contained in the face where $\zeta\neq 0$, so  $e_{S}\zeta=\zeta$. (This dichotomy follows from the fact that any face is a union of strata.) In other words, either $e_{S}\Delta_{I_{\zeta}}=0$ or $\Delta_{S}\rho(t\abs\zeta)=0$. Therefore, the terms in the above sum for which $I_{1}\neq \emptyset$ are $0$, and we get that
\[\Delta_{I}(f-h)=\rho(t\abs\zeta)\Delta_{I\cup I_{\zeta}}f\]

As $w_{I}^{-\delta}\Delta_{I}f$ is continuous and vanishes on the strata in $I$, $e_{S}w_{I}^{-\delta}\Delta_{I}f$ is continuous and vanishes on all strata in $I$. Therefore, as $ w_{I}^{-\delta}e_{S}w_{I}^{\delta}$ is bounded by $1$ and continuous everywhere apart from strata in $I$, $w_{I}^{-\delta}e_{S}\Delta_{I}f$ is continuous and vanishes on all strata in $I$. Therefore, $w_{I}^{-\delta}\Delta_{I\cup I_{\zeta}}f$ must also be continuous and vanish on all the strata in $I\cup I_{\zeta}$. It follows that given any compact subset and $\epsilon>0$,  for $t$ large enough, $\abs{w^{-\delta}_{I}\Delta_{I}(f-h)}<\epsilon$ on the given compact subset. As $\nabla\zeta$ is proportional to $\zeta$, the derivatives of $\rho(t\abs\zeta)$ are bounded independent of $t$ using the standard metric and connection defined on page \pageref{standard metric}. It follows that $h\rightarrow f$ in $C^{k,\delta}$ as $t\rightarrow \infty$.

Note that the new function $h$ still satisfies the two conditions of our lemma. We may therefore repeat the above process for different $\zeta$ and approximate $f$ in $C^{k,\delta}$ by a function $h$ satisfying the conditions of our lemma with the following extra property: for all $\zeta$ in the appropriate form $\zeta=\totl{\e a\tilde z^{\alpha}}$, where $\zeta$ is small enough,  $\Delta_{I_{\zeta}}h=0$.  This $h$ therefore has the property that it is independent of $\zeta $ for $\zeta$ small enough. As $\nabla^{k} h$ exists and is continuous, it follows that $h$ is $C^{k}$. Then the estimates from Lemma \ref{smooth ckdelta} imply that we may approximate $h$ in $C^{k,\delta}$ by a smooth function, so if $f$ satisfies  the conditions of our lemma, there exists a sequence of smooth functions converging in $C^{k,\delta}$ to $f$.

\stop

Lemma \ref{ckdelta} has the following immediate corollary:

\begin{cor}
A continuous function $f$ is $\C k\delta$ if and only if $\nabla^{k}f$ exists and is continuous, and restricted to compact subsets $\abs f_{k,\delta'}$ is finite for all $\delta'<\delta$.
\end{cor}

Note that this Corollary holds for $\C k\delta$ which is the intersection of all $C^{k,\delta'}$ for $\delta'<\delta$. It does not hold for $C^{k,\delta}$. For example if $f$ is continuous and $\abs f_{0,\delta'}$ is finite, this implies that for all $\delta''<\delta'$, $w_{I}^{-\delta''}\Delta_{I}f$ extends to a continuous function which vanishes on all strata in $I$, but does not imply that $w_{I}^{-\delta'}\Delta_{I}f$ extends to a continuous function.

\

\

  In the following pages, we show that we can replace smooth functions with $C^{\infty,\delta}$ functions in the definition of exploded manifolds to create a category of $C^{\infty,\delta}$ exploded manifolds.

\begin{defn}[$C^{k,\delta}$ exploded function]
A $ C^{k,\delta}$ exploded function $f\in \mathcal E^{k,\delta,\times}(\mathbb R^{n}\times\et mP)$ is a function of the form 
\[f(x,\tilde z):= g(x,\tilde z)\tilde z^{\alpha}\e a\text{ where }g\in C^{k,\delta}\lrb{\mathbb R^{n}\times\et mP,\mathbb C^{*}},\ \alpha\in\mathbb Z^{m},\ \e a\in\e{\mathbb R}\]

Similarly, $\mathcal E^{\infty,\delta,\times}=\bigcap_{k}\mathcal E^{k,\delta,\times}$ and $\mathcal E^{\infty,\underline{1},\times}=\bigcap_{\delta<1}\mathcal E^{\infty,\delta,\times}$.

Say that a sequence of exploded functions $g^{i}\tilde z^{\alpha_{i}}\e {a_{i}}$ converge with a given regularity if the sequence of functions $g^{i}$ does and $\alpha_{i}$ is eventually constant. Note that this is a non Hausdorff topology because there is no condition on the sequence $a_{i}$. Say that the sequence converges strongly, if $g^{i}$ converges and the sequences $\alpha_{i}$ and $a_{i}$ are eventually constant.

\

A $C^{k,\delta}$, $C^{\infty,\delta}$, or $\C\infty1$ exploded manifold is an  abstract exploded space locally isomorphic to $\mathbb R^{n}\times \et mP$ with the sheaf $\mathcal E^{k,\delta,\times}$, $\mathcal E^{\infty,\delta,\times}$, or $\mathcal E^{\infty,\underline 1,\times}$ respectively.
\end{defn}

\

\begin{lemma}\label{linear change}

If 
\[\alpha:\mathbb R^{n}\times\et mP\longrightarrow \mathbb R^{n'}\times\et {m'}Q\]
is a `linear' map, so
\[\alpha (x,\tilde z):=\lrb{Mx,\tilde z^{\alpha^{1}},\dotsc,\tilde z^{\alpha^{m'}}}\] 
where $M$ is a $n$ by $n'$ matrix and $\alpha^{i}_{j}$ is a $m$ by $m'$ matrix with integer entries,\\ then the map $\alpha$ preserves $C^{k,\delta}$ in the sense that given any function $f\in C^{k,\delta}(\mathbb R^{n'}\times\et {m'}Q)$,  
\[ f\circ \alpha\in C^{k,\delta}\lrb{\mathbb R^{n}\times\et mP}\]
\end{lemma}

\pf

The important observation in this proof is that if $f$ is a continuous function on $\mathbb R^{n'}\times \et {m'}Q$, $S$ is a stratum of $P$ and $S'$ is the stratum of $Q$ which contains $\totb\alpha(S)$, then 
 \[e_{S}(f\circ \alpha)=(e_{S'}f)\circ\alpha\]
It follows that if $I$ is any collection of strata of $P$ and $I'$ the corresponding collection of strata of $Q$, then 
\[\Delta_{I}(f\circ\alpha)=(\Delta_{I'}f)\circ\alpha\]
Also, $w_{I'}\circ\alpha$ is a finite sum of absolute values of monomials which vanish on the strata in $I$, therefore, $w_{I'}\circ\alpha$ is bounded on compact subsets by a constant times $w_{I}$. Therefore, $\abs{w_{I}^{-\delta}\Delta_{I}(f\circ \alpha)}$ is bounded on compact subsets by a constant times $\abs{(w_{I'}^{-\delta}\Delta_{I'}f)\circ \alpha}$. We may bound derivatives similarly. It follows that if $f_{i}$ is a sequence of smooth functions converging to $f$ in $C^{k,\delta}$, then $f_{i}\circ \alpha$ is a sequence of smooth functions converging to $f\circ \alpha$ in $C^{k,\delta}$. So $f\circ \alpha$ is $C^{k,\delta}$ if $f$ is.

\stop

\begin{lemma}\label{exponentiate vector field}
 Any $C^{k,\delta}$ section of $T\lrb{\mathbb R^{n}\times \et mP}$ may be considered in standard coordinates as a $C^{k,\delta}$ map  
 \[(f_{\mathbb R^{n}},f_{1},\dotsc, f_{m}):\mathbb R^{n}\times \et mP\longrightarrow \mathbb R^{n}\times \mathbb C^{m}\]
  We can define an exponential map of the form 
\[\exp(f)(x,\tilde z_{1},\dotsc,\tilde z_{m}):=\lrb{x+f_{\mathbb R^{n}}(x,\tilde z), e^{f_{1}(x,\tilde z)}\tilde z_{1},\dotsc,e^{f_{m}(x,\tilde z)}\tilde z_{m}}\]

 If $h$ is in $C^{k,\delta}$, then $h\circ\exp f$ is.

\end{lemma}

\pf

We shall show that  $\abs{h\circ \exp f}_{k,\delta}$ on a given compact subset $U$ can be bounded by $\abs h_{k,\delta}c(f)$ on a $\sup_{U} \abs f$ neighborhood of $U$ in the standard metric defined on page \pageref{standard metric}. 

Note that 
\[e_{S}(h\circ\exp f)=(e_{S}h)\circ\exp (e_{S}f) \]
Therefore, 
\[\begin{split}\Delta_{S}(h\circ\exp f)&=h\circ \exp f-(e_{S}h)\circ \exp f+(e_{S}h)\circ \exp f-(e_{S}h)\circ\exp(e_{S}f)
\\&:=(\Delta_{S}h)\circ \exp f +(e_{S}h\circ \exp \circ \Delta_{S})f
\end{split}\]
Induction on the number of strata in $I$ implies that  we can rewrite $\Delta_{I}(h\circ\exp f)$ in the following form
\begin{equation}\label{comp1}\Delta_{I}(h\circ\exp f)=\sum_{I'\coprod I''=I}\lrb{\lrb{ e_{I''}\Delta_{I'}h}\circ\exp\circ \Delta_{I''}}f  \end{equation}

As an example for interpreting the notation above, we write
\[\lrb{h\circ\exp\circ \Delta_{S}}f:=h\circ\exp f-h\circ\exp (e_{S}f)\]
as opposed to 
\[h\circ\exp\lrb{\Delta_{S}f}:=h\circ\exp(f-e_{S}f)\]

The weight function $w_{I}\circ \exp f$ differs from $w_{I}$ by a factor which is bounded by a constant to the power of the size of $f$ - this is because the weight function is a sum of absolute values of monomials in $\tilde z$, and the size of $f$ is measured in the standard metric in which the real and imaginary parts of $\tilde z_{i}\frac{\partial}{\partial \tilde z_{i}} $ are orthonormal. Therefore, there exists some constant so that 
\begin{equation}\label{comp2}\sup_{U}\abs{ w_{I}^{-\delta}(\Delta_{I}h)\circ \exp f}\leq \sup_{\exp f(U)}\abs{w_{I}^{-\delta}\Delta_{I}h}c^{\abs f}\end{equation}
Now, to bound expressions of the form:
\[\lrb{\lrb{\Delta_{I'}h}\circ \exp \circ \Delta_{I''}}f\]
introduce a real variable $t_{S}$ for every stratum  $S\in I''$, and define
\[\Phi_{I''}(t):=\lrb{\prod_{S\in I''} (t_{S}\Delta_{S}+e_{S})}f\]
Using the notation $D_{I''}:=\prod_{S\in I''}\frac \partial{\partial t_{S}}$, we can rewrite
\begin{equation}\label{comp6}\lrb{\lrb{\Delta_{I'}h}\circ \exp \circ \Delta_{I''}}f:=\int_{0}^{1}\dotsb\int_{0}^{1}D_{I''}\lrb{\Delta_{I'}h(\exp \Phi_{I''}(t)) }dt \end{equation}
In order to bound the integrand of the above, we expand:
\begin{equation}\label{comp4}\begin{split}D_{I''}\lrb{\Delta_{I'}h(\exp \Phi_{I''}(t)) }&=\sum \nabla^{l}(\Delta_{I'}h)(D_{I_{1}}\Phi_{I''})\dotsb(D_{I_{l}}\Phi_{I''})
\\& =\sum \Delta_{I'}\lrb{\nabla^{l}h}(\Delta_{I_{1}}\Phi_{I''-I_{1}})\dotsb(\Delta_{I_{l}}\Phi_{I''-I_{l}}) \end{split}\end{equation}
where the sum above is over all partitions $I_{1},\dotsc,I_{l}$ of $I''$, and $\nabla^{l}$ indicates the $l$th covariant derivative using the standard connection defined on page \pageref{standard metric}. (For notational simplicity, the fact that $ \nabla^{l}(\Delta_{I'}h)$ is to be evaluated at $\exp \Phi_{I''}(t)$ has been suppressed in the right hand side of the above expression.) As on compact subsets $w_{I}^{-\delta}$ is bounded by a constant times $w_{I'}^{-\delta}w_{I_{1}}^{-\delta}\dotsc w_{I_{l}}^{-\delta}$,  we get the following estimate on compact subsets:

\begin{equation}\label{comp3}\abs{w_{I}^{-\delta}D_{I''}\lrb{\Delta_{I'}h(\exp \Phi_{I''}(t)) }}\leq c\sum\abs{w_{I'}^{-\delta}\Delta_{I'}\lrb{\nabla^{l}h}}\prod_{i=1}^{l}\abs{w_{I_{i}}^{-\delta}\Delta_{I_{i}}\Phi_{I''-I_{i}}}\end{equation}
A similar inequality for $\abs{w_{I}^{-\delta}D_{I''}\lrb{\Delta_{I'}\nabla^{m}\lrb{h(\exp \Phi_{I''}(t))} }}$ can be obtained by differentiating equation \ref{comp4}.
\begin{equation}\label{comp5}\abs{w_{I}^{-\delta}D_{I''}\lrb{\Delta_{I'}\nabla^{m}\lrb{h(\exp \Phi_{I''}(t)) }}}\leq c\sum\abs{w_{I'}^{-\delta}\Delta_{I'}\lrb{\nabla^{l}h}}\prod_{i=1}^{l}\abs{w_{I_{i}}^{-\delta}\Delta_{I_{i}}\nabla^{m_{i}}\Phi_{I''-I_{i}}}\end{equation}
The above sum is  now over all partitions of $I''\coprod \{1,\dotsc,m\}$, where the number of sets in the partition is $l$, the number of integers in each set is $m_{i}$ and the intersection of the set with $I''$ is $I_{i}$.
 We can bound $w^{-\delta}_{I_{i}}\Delta_{I_{i}}\nabla^{m_{i}}\Phi_{I''-I_{i}}$ by  a sum of terms of the form $\abs{e_{S}w^{-\delta}_{I_{i}}\Delta_{I_{i}}\nabla^{m_{i}}f}$, therefore equations \ref{comp1}, \ref{comp2}, \ref{comp6} and \ref{comp5} give that for all $I$ and $m$ so that $\abs I+m\leq k$, $w_{I}^{-\delta}\Delta_{I}\nabla^{m}(h\circ \exp f)$ is continuous and vanishes on all the strata in $I$. Lemma \ref{ckdelta} then implies that $h\circ \exp f$ is $C^{k,\delta}$.
  
\stop

Consider a map between coordinate charts for which the pull back of exploded coordinate functions is in $\mathcal E^{k,\delta,\times}$ and the pullback of real coordinate functions is $C^{k,\delta}$. Any map of this form factors as a composition of maps in the form of Lemma \ref{linear change} and Lemma \ref{exponentiate vector field}; (first a map in the form of Lemma \ref{linear change} to the product of the domain and target, then a map of the form of Lemma \ref{exponentiate vector field}, then a projection to the target, which is of the form of Lemma \ref{linear change}.) We therefore have the following: 

\begin{cor}\label{coordinate ckdelta}
A morphism is $ C^{k,\delta}$ if and only if the  pull back of exploded coordinate functions are $\mathcal E^{k,\delta,\times}$ functions, and the pull back of real coordinate functions are $ C^{k,\delta}$ functions.

\end{cor}

We can of course define a $C^{k,\delta}$ vector field on an exploded manifold $\ex B$ as a section of $T\ex B$ which is a $C^{k,\delta}$ morphism. Corollary \ref{coordinate ckdelta} implies that in coordinates, $C^{k,\delta}$ vector fields are the vector fields which are $C^{k,\delta}$ functions times the standard basis vector fields.

\

\

The definition of convergence for functions now generalizes in a straightforward way to convergence for $C^{k,\delta}$ maps.

\begin{defn}[$C^{k,\delta}$ convergence]
A sequence of $C^{k,\delta}$ exploded maps 
$f^i:\ex A\longrightarrow \ex B$ converges (strongly) to $f:\ex A\longrightarrow \ex B$ in $C^{k,\delta}$ if  the pullback  under $f^{i}$ of any local coordinate function on $\ex B$ converges (strongly) in $C^{k,\delta}$ to the pullback under $f$.

\end{defn}

\section{Almost complex structures}

\begin{defn}[(Almost) complex structure]
An almost complex structure $J$ on a smooth exploded  manifold $\ex B$ is an endomorphism of $T\ex B$ given by a smooth section of $ T\ex B\otimes T^*\ex B$ which squares to become multiplication by $-1$, so that given any exploded function $\tilde z\in \mathcal E^{\times} (\ex B)$ and  integral vector $v\in {}^{\mathbb Z}T_p\ex B$  
\[(Jv)(\tilde z)=i(v(\tilde z))\]

An almost complex structure $J$ is a complex structure if there exist local coordinates $\tilde z\in\et mP$ so that for all vector fields $v$, $iv(\tilde z_{j})=(Jv)(\tilde z_{j})$. 
\end{defn}

This differs from the usual definition of an almost complex structure only in the extra requirement that integral vectors satisfy $(Jv)\tilde z=iv\tilde z$. Integral vectors are defined on page \pageref{integral}. They are the vectors which satisfy $v(\tilde z)\tilde z^{-1}$ is an integer for all exploded functions $\tilde z$.

 For example, on $\et 11$ the integral vectors are the integer multiples of the real part of $\tilde z\frac\partial{\partial \tilde z}$ in the region where $\totl{\tilde z}=0$. Our definition requires that $J$ of the real part of $\tilde z\frac\partial{\partial \tilde z}$ is the imaginary part in this region.  This extra requirement makes holomorphic curves $\C\infty 1$ exploded maps, and makes the tropical part of holomorphic curves piecewise linear one complexes. If it did not hold, then we would need to use a different version of the exploded category using $\mathbb R^{*}$ instead of $\mathbb C^{*}$ to explore holomorphic curve theory. The analysis involved would be significantly more difficult. 

The following assumption allows us to use standard (pseudo)holomorphic curve results.

\begin{defn}[Civilized almost complex structure]\label{civilized}
An almost complex structure $J$ on $\ex B$ is civilized if it induces a smooth almost complex structure on the smooth part of $\ex B$. In other words, every point in $\ex B$ has a neighborhood $U$ with some map 
\[\psi:U\longrightarrow \mathbb R^{N}\]
so that
\begin{itemize}
\item There exists a smooth almost complex structure $J'$ on $\mathbb R^{N}$ so that $\psi$ is holomorphic in the sense that $d\psi\circ J=J'\circ d\psi$.
\item $\psi$ is an embedding of the smooth part of $U$ in the sense that any smooth map $U\longrightarrow \mathbb R$ is the pullback of some smooth function from $\mathbb R^{N}$.  
\end{itemize}

\end{defn}

The word civilized should suggest that our almost complex structure is well behaved in a slightly unnatural way. Any complex structure is automatically civilized, and there are no obstructions to modifying an almost complex structure to civilize it.  If we assume our almost complex structure is civilized, then the removable singularity theorem may be used on the smooth part of holomorphic curves to prove that holomorphic curves are smooth maps. With an uncivilized almost complex structure, the best that can be expected of holomorphic curves is that they will be $\C\infty1$.

\begin{defn}[Exploded curve]
An exploded  curve is a $2$ real dimensional, complete exploded manifold with a complex structure $j$. 

A holomorphic curve is a holomorphic map of an exploded curve to an almost complex exploded manifold.

A smooth or $\C\infty1$ exploded curve is a smooth or $\C\infty1$ map of an exploded curve to a exploded manifold. 
\end{defn}

\begin{itemize}
\item By a smooth component of a holomorphic curve $\ex C$, we shall mean a stratum of $\ex C$ which is a connected, punctured Riemann surface. 
\item By an internal edge or node of $\ex C$ we shall mean a stratum of $\ex C$ isomorphic to $\et 1{(0,l)}$.
\item  By a puncture or end of $\ex C$, we shall mean a stratum of $\ex C$ isomorphic to $\et 1{(0,\infty)}$.  
\end{itemize}

The information in a holomorphic curve $\ex C$ is equal to the information of a nodal Riemann surface plus gluing information for each node parametrized by $\mathbb C^{*}\e{(0,\infty)}$ (for more details of this gluing information see example \ref{local family model} on page \pageref{local family model}). 

With one exception, all strata of exploded curves are either smooth components, edges or punctures. The exception is $\ex T$. (Actually, there would be further exceptions if we had defined exploded manifolds differently, and allowed the quotient of $\ex T$ by $\tilde z\mapsto c\e a\tilde z$ to be an exploded manifold.) 

\includegraphics{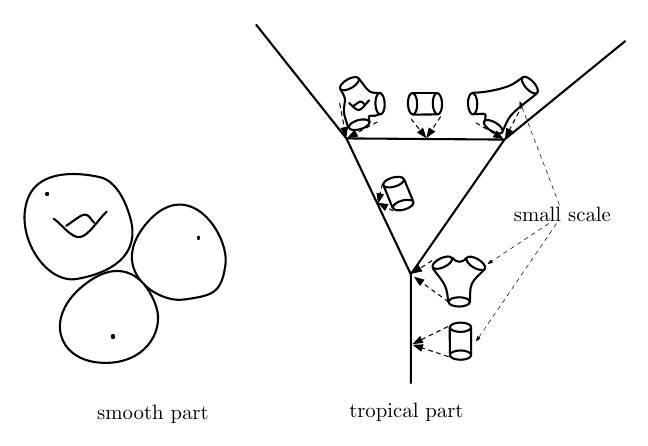}

The above is a picture of an exploded  curve with $3$ smooth components, $3$ punctures and $3$ internal edges. On the left hand side is the smooth part, where the smooth components are the $3$ pictured Riemann surfaces, the 3 punctures are the $3$ points marked, and the  $3$ internal edges correspond to the three nodes where the Riemann surfaces are joined. On the right hand side is a composite picture of the tropical part of our curve, and some of  the small scale of our curve, which is the topological space obtained from putting a smooth metric on our curve. The punctures correspond to the free edges of the tropical curve, the internal edges correspond to the internal edges of the tropical curve, and the smooth components correspond to the vertices of the curve. This exploded curve has genus $2$ - part of the genus can be seen in the tropical part, and part of the genus can be seen in one of the smooth components.

\begin{example}[Balancing condition for exploded curves in $\ex T^{n}$]\label{curve in tropical plane}
\end{example}
Consider a smooth curve $f:\ex C\longrightarrow \ex T^{n}$. This is given by $n$ exploded functions $f^{*}(\tilde z_{1}),\dotsc,f^{*}(\tilde z_{n})\in\mathcal E^{\times}(\ex C)$.

 Each smooth component of $\ex C$ is sent to the $\lrb{\mathbb C^{*}}^{n}$ worth of points over a particular point in the tropical part $\totb{\ex T^{n}}$. In particular, $f$ restricted to a smooth component gives a smooth map of the corresponding punctured Riemann surface to $\lrb{\mathbb C^{*}}^{n}$. Around each puncture of a smooth component, there is some   homology class  $\alpha\in H_{1}\lrb{\lrb{\mathbb C^{*}}^{n},\mathbb Z}$ of  a loop around the puncture.   Of course, the sum of all such homology classes from punctures of a smooth component is zero.
 
Now consider $f$ in a $\et 11$ coordinate chart around a puncture. In these coordinates,

 \[f(\tilde w)=(g_{1}(\totl{\tilde w})\e {a_{1}}\tilde w^{\alpha_{1}},\dotsc,g_{n}(\totl{\tilde w})\e{a_{n}}\tilde w^{\alpha_{n}})\ \ g_{i}\in C^{\infty}(\mathbb C,\mathbb C^{*}),\ \alpha\in\mathbb Z^{n}\]
 
 Similarly, $f$ in a $\et 1{[0,l]}$ coordinate chart around an internal edge can be written as

 \[f(\tilde w)=(g_{1}(\zeta_{1},\zeta_{2})\e {a_{1}}\tilde w^{\alpha_{1}},\dotsc,g_{n}(\zeta_{1},\zeta_{2})\e{a_{n}}\tilde w^{\alpha_{n}})\]
\[\text{ where } \zeta_{1}=\totl w,\ \zeta_{2}=\totl{\e lw^{-1}},\  g_{i}\in C^{\infty}(\mathbb C^{2},\mathbb C^{*}),\ \alpha\in\mathbb Z^{n}\]
 
Again $\alpha\in \mathbb Z^{n}$ can be regarded as the homology class in $H_{1}\lrb{\lrb{\mathbb C^{*}}^{n},\mathbb Z}$ of a loop around the puncture. The tropical part of this map is $x\mapsto a+\alpha x$, so $\alpha$ determines the derivative of the tropical part of $f$.
  Therefore the sum of  all the derivatives of $\totb f$ exiting a vertex sum to $0$. This can be viewed as some kind of conservation of momentum condition for the tropical part of our curve, $\totb{f}$. In tropical geometry, this is called the balancing condition.

\includegraphics{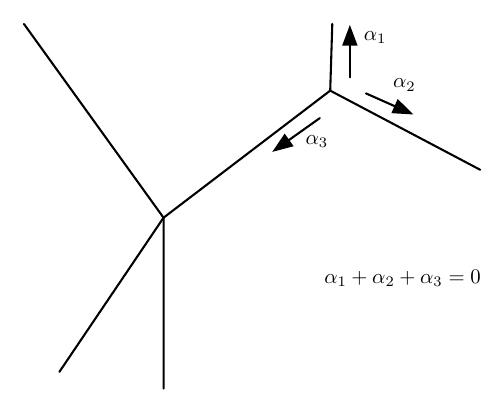}

\begin{example}[Curves as locus of non-invertibility  of polynomials on $\ex T^{n}$]\label{tropical curve}
\end{example}
One way to consider the image of some holomorphic curves  in $\ex T^{n}$ is as the `locus of noninvertiblity' of some set of polynomials 
\[P_{i}(\tilde z):=\sum \tilde c_{i,\alpha}\tilde z^{\alpha}\ i=1,\dotsc,n-1\]
We can consider the set
\[ Z_{\{P_{i}\}}:=\left\{\tilde z\text{ so that }P_{i}(\tilde z)\in 0\e{\mathbb R}\ \forall i\right\}\]
Suppose that for all points $p\longrightarrow Z_{\{P_{i}\}}$, the differentials $\{dP_{i}\}$ at $p$ are linearly independent. Then theorem \ref{zeroset} proved on page \pageref{zeroset} implies that $Z_{\{P_{i}\}}$ is the image of some holomorphic curve.

Let us examine the set $Z_{\{P_{i}\}}$ more closely. For any point $\tilde z_{0}$, denote by $S_{i,\tilde z_{0}}$ the set of exponents $\alpha$ so that $\totb{P_{i}(\tilde z_{0})}=\totb{\tilde c_{i,\alpha}\tilde z^{\alpha}}$. Then there exists some neighborhood of $\totb{\tilde z_{0}}$ in $\totb{\ex T^{n}}$ so that 
\[P_{i}=\sum_{\alpha\in S_{i,\tilde z_{0}}}\tilde c_{i,\alpha}\tilde z^{\alpha}\]
The points inside $Z_{\{P_{i}\}}$ over $\totb{\tilde z_{0}}$ are then given by solutions of the equations 
\[\sum_{\alpha\in S_{i,\tilde z_{0}}} c_{i,\alpha}z^{\alpha}=0\text{ where } \tilde c_{i,\alpha}=c_{i,\alpha}\totb{\tilde c_{i,\alpha}}\text{ and }\tilde z=z\totb{\tilde z}\]
Note that the above equation has solutions for $z\in (\mathbb C^{*})^{n}$ if and only if $S_{i,\tilde z_{0}}$ has more than $1$ element. This corresponds to the tropical function $\totb{P_{i}}$ (which is continuous, piecewise integral affine, and convex) not being smooth at $\totb{\tilde z_{0}}$. We therefore have that $\totb{Z_{\{P_{i}\}}}$ is contained in the intersection of the non-smooth locus of the tropical polynomials $\totb{P_{i}}$.

\section{Fiber products}

\begin{defn}[Transverse]
Two smooth (or $C^{k,\delta}$) exploded morphisms 
\[\ex A\xrightarrow{f}\ex C\xleftarrow{g}\ex B\]
are transverse if for every pair of points $p_1\longrightarrow\ex A$ and $p_2\longrightarrow \ex B$ so that $f(p_1)=g(p_2)$,  $df(T_{p_1}\ex A)$ and $dg(T_{p_2}\ex B)$ span $T_{f(p_1)}\ex C$.
\end{defn}

\begin{defn}[Fiber product]
If $f$ and $g$ are transverse smooth (or $C^{k,\delta}$) exploded morphisms,
\[\ex A\xrightarrow{f}\ex C\xleftarrow{g}\ex B\]
the fiber product $\ex A\fp fg\ex B$ is the unique smooth (or $C^{k,\delta}$) exploded  manifold   with maps to $\ex A$ and $\ex B$ so that the following diagram commutes
\begin{displaymath}
\begin{array}{cll}
\ex A\fp fg\ex B & \longrightarrow & \ex A \\ 
 \downarrow & × & \downarrow \\ 
\ex B & \longrightarrow & \ex C
\end{array}
\end{displaymath}
and with the usual universal property that given any commutative diagram
\begin{displaymath}
\begin{array}{lll}
\ex D & \longrightarrow  & \ex A \\ 
\downarrow & × & \downarrow \\ 
\ex B & \longrightarrow & \ex C
\end{array}
\end{displaymath}
there exists a unique morphism $\ex D\longrightarrow \ex A\fp fg\ex B$ so that the following diagram commutes
\begin{displaymath}
\begin{array}{llc}
\ex D & \rightarrow & \ex A \\ 
\downarrow & \searrow & \uparrow \\ 
\ex B & \leftarrow & \ex A\fp fg\ex B
\end{array}
\end{displaymath}
\end{defn}

The universal property of fiber products implies that they are unique if they exist. We shall prove their existence in the case of transversality in the next few lemmas.

\begin{lemma}\label{smooth intersect} Let $U$ be a standard exploded coordinate chart, and let $f:U\times \mathbb R^{n}\longrightarrow \mathbb R^{n}$ be a smooth function so that there exists a constant $c<1$ so that if $v$ is any vector in the $\mathbb R^{n}$ direction, 
\[\norm{df(v)-v}\leq c\norm v\]
 Then there exists a unique smooth map $g:U\longrightarrow \mathbb R^{n}$ so that $f(u,g(u))=0$. 
\end{lemma}

\pf

This is a version of the implicit function theorem which follows from the smooth case. We can consider the smooth part $\totl{U}$ of $U$ as a subset of $\mathbb C^{m}$, and extend $f$ to be a smooth function on $\mathbb C^{m}\times \mathbb R^{n}$ still obeying the condition that the derivative of $f$ in the $\mathbb R^{n}$ direction is close to the identity, so there exists some constant $c<1$ so that for any vector $v$ in the $\mathbb R^{n}$ direction, 

\[\norm{df(v)-v}\leq c \norm v\]

Then set $g_{0}:\mathbb C^{m}\longrightarrow \mathbb R^{n}$ to be $0$, and let 
\[g_{i}(z):=g_{i-1}(z)-f(z,g_{i-1}(z))\]
 Note that $\abs{f(z,g_{i}(z))}\leq c\abs{f(z,g_{i-1}(z))}$, so $g_{i}$ is a Cauchy sequence which converges to a continuous map $g:\mathbb C^{m}\longrightarrow \mathbb R^{n}$ so that $f(z,g(z))=0$. The implicit function theorem gives that this map $g$ must be smooth. The restriction of $g$ to $\totl U\subset \mathbb C^{m}$ gives the required solution. Uniqueness follows from the fact that for all $z$ there is a unique $x$ so that $f(z,x)=0$ because the map $x\to x-f(z,x)$ is a contraction.

\stop

\begin{lemma}\label{ckdelta intersect} Suppose that  $f:\ex B\times \mathbb R^{n}\longrightarrow \mathbb R^{n}$ is a $C^{k,\delta}$ map   so that for some point $p\in \ex B\times\mathbb R^{n}$, $f(p)=0$ and  the derivative of $f$ at $p$ restricted to the $\mathbb R^{n}$ direction is bijective. Then there exists an open neighborhood of $p$, $U\times U'\subset \ex B\times \mathbb R^{n}$ so that for each $u\in U$, there exists  a  unique point $g(u)\in U'$ so that $f(u,g(u))=0$. The resulting map $g:U\longrightarrow U'$ is $C^{k,\delta}$. 

\end{lemma} 

\pf

The existence and uniqueness of some map $g$  so that $f(u,g(u))=0$  on an appropriately small neighborhood $U\times U'$ follows from the usual inverse function theorem applied to $f$ restricted to each $\mathbb R^{n}$ slice.  We must verify that this $g$ is $C^{k,\delta}$. First, note that  $f$ restricted to each $\mathbb R^{n}$ slice depends continuously on $u$, so $g$ is continuous. The usual implicit function theorem also implies that if we choose our neighborhood small enough,  $\nabla ^{k}g$ exists and is continuous.

From now on, we shall assume that $U\times U'$ is chosen within a single standard coordinate chart so that the operations $e_{S}$ make sense on functions defined on $U$ or $U\times U'$. We shall now prove that for any collection of strata $I$ containing at most $k$ strata,  $w_{I}^{-\delta}\Delta_{I}g$ is continuous and vanishes on strata in $I$. Suppose for induction that this holds for all collections of strata containing at most $\abs I-1$ strata.
Introduce a real variable $t_{S}$ for every strata  $S\in I$, and define
\[\phi_{I}(t):=\lrb{\prod_{S} (t_{S}\Delta_{S}+e_{S})}g\]
Then using the argument following equation (\ref{comp1}) on page \pageref{comp1},
\begin{equation}\label{kjh}\begin{split}0&=\Delta_{I}\lrb{f(u,g)}=\sum_{I_{1}\coprod I_{2}=I}\int_{0}^{1}\dotsb\int_{0}^{1}D_{I_{2}}\lrb{e_{I_{2}}\Delta_{I_{1}}f(u,\phi_{I_{2}})}dt_{I_{2}}
\\& =\sum_{\coprod_{i=1}^{m} I_{i}=I}\int_{0}^{1}\dotsb\int_{0}^{1}e_{I-I_{1}}\Delta_{I_{1}}\nabla^{m-1}f(u,\phi_{I-I_{1}})(D_{I_{2}}\phi_{I-I_{1}})\dotsb (D_{I_{m}}\phi_{I-I_{1}})dt_{I-I_{1}}
\\ &= \int_{0}^{1}\dotsc \int_{0}^{1}e_{I}\nabla f(u,\phi_{I})(\Delta_{I}g)dt_{I}
\\ & + \text{ terms not involving $\Delta_{I}g$}\end{split}\end{equation}

Suppose that $\int_{0}^{1}\dotsc \int_{0}^{1}e_{I}\nabla f(u,\phi_{I})dt_{I}$ has a uniformly bounded inverse when restricted to the $\mathbb R^{n}$ direction and considered for each $u$ as a linear transformation $\mathbb R^{n}\longrightarrow \mathbb R^{n}$. As $\nabla f$ is continuous and invertible restricted to the $\mathbb R^{n}$ direction at $p$, this is always true if we have chosen our neighborhood $U\times U'$ small enough. Then $w_{I}^{-\delta}\Delta_{I}g$ will be continuous and vanish on all strata in $I$ if and only if the expression \[w_{I}^{-\delta}\int_{0}^{1}\dotsc \int_{0}^{1}e_{I}\nabla f(u,\phi_{I})(\Delta_{I}g)dt_{I}\] is continuous and vanishes on all strata in $I$. We shall show that this is the case using the other terms in equation \ref{kjh}. 

First, note that $w_{I}^{-\delta}x$ is continuous and vanishes on all strata of $I$ if $w_{I_{1}}^{-\delta}\dotsb w_{I_{m}}^{-\delta}x$ does. We have that $w^{-\delta}_{I_{1}}\Delta_{I_{1}}\nabla^{m-1}f$ is continuous and vanishes on all strata in $I_{1}$, therefore the same is true of $w_{I_{1}}^{-\delta}e_{I-I_{1}}\Delta_{I_{1}}\nabla^{m-1}f(u,\phi_{I-I_{1}})$. Also, so long as $I_{j}\subsetneq I$, our inductive hypothesis implies that $w^{-\delta}_{I_{j}}\Delta_{I_{j}}g$ is continuous and vanishes on all strata of $I_{j}$. The term $D_{I_{j}}\phi_{I-I_{1}}$ is a signed sum of terms in the form of $e_{I'}\Delta_{I_{j}}g$, therefore, $w_{I_{j}}^{-\delta}D_{I_{j}}\phi_{I-I_{1}}$ is continuous and vanishes on all strata in $I_{j}$.  Therefore, our inductive hypothesis combined with equation \ref{kjh} implies that 
$w_{I}^{-\delta}\int_{0}^{1}\dotsc \int_{0}^{1}e_{I}\nabla f(u,\phi_{I})(\Delta_{I}g)dt_{I}$ is continuous and vanishes on all strata in $I$. Therefore, so long as we have chosen $U\times U'$ small enough, $w_{I}^{-\delta}\Delta_{I}g$ will be continuous and vanish on all strata in $I$. 

We may now complete the proof by induction on the number of derivatives, $k$. Suppose $k=1$. As $w_{S}^{-\delta}\Delta_{S}g$ is continuous and vanishes on $S$, and $\nabla g$ is continuous by the usual implicit function theorem, $g$ is $C^{1,\delta}$, and the lemma holds for $k=1$. Now suppose that the lemma holds for $k-1$, and that $f$ is $C^{k,\delta}$. Our inductive hypothesis implies that $g$ is $C^{k-1,\delta}$, and we have proven that for any collection $I$ of at most $k$ strata $w_{I}^{-\delta}\Delta_{I}g$ will be continuous and vanish on all strata in $I$.  Lemma \ref{ckdelta} implies that it remains to prove that $\nabla g$ is $C^{k-1,\delta}$. The implicit function theorem gives us the following formula for $\nabla g$, where $\nabla_{\mathbb R^{n}}f$ indicates the derivative of $f$ restricted to the $\mathbb R^{n}$ direction and $\nabla_{U}f$ indicates the derivative of $f$ restricted to the $U$ direction in $U\times U'$.
\[\nabla g= \lrb{\nabla_{\mathbb R^{n}}f}^{-1}(u,g)(-\nabla_{U}f)\]
As $\nabla_{\mathbb R^{n}} f$, $g$ and $\nabla_{U}f$ are all $C^{k-1,\delta}$, the above equation gives that $\nabla g$ is $C^{k-1,\delta}$, therefore $g$ is $C^{k,\delta}$ as required.

\stop

\

The following is an example of interesting behavior that can happen in a fiber product.

\begin{example}[A fiber product]
 
\end{example}

 Consider the map $f:\et mP\longrightarrow\ex T^n$
given by 
\[f(\tilde z)=(\tilde z^{\alpha^1},\dotsc,\tilde z^{\alpha^n})\]
Denote by $\alpha$ the $m\times n$ matrix with entries $\alpha^i_j$. The derivative of $f$ is surjective if $\alpha:\mathbb R^m\longrightarrow\mathbb  R^n$ is surjective. Denote by $\abs\alpha\in\mathbb N$ the index of $\alpha^1\wedge\dotsb\wedge\alpha^n\in \bigwedge^n(\mathbb Z^m)$. (In other words, $\abs{\alpha}\in\mathbb N$ is the largest nonnegative integer so that the above wedge is  $\abs\alpha$ times a nonzero element of $\bigwedge^{n}(\mathbb Z^{m})$.) The fiber product of $f$ with the point $(1,\dotsc,1)$ corresponds to the points in $\et mP$ so that $\tilde z^{\alpha^i}=1$ for all $i$.  This is then equal to $\abs \alpha$ copies of $\et {m-n}{\e{\ker \alpha}\cap P}$ where we identify $\e{\mathbb R^{m-n}}=\e{\ker\alpha}$. 

As a simple example, consider  $f(\tilde z):=\tilde z^{2}:\ex T\longrightarrow \ex T$. Then $f^{-1}(1)$ consists of the two points $\tilde z=\pm 1$. The subspace topology on $f^{-1}(1)$ is the trivial topology, but the correct topology on this fiber product is the discrete topology.

 This example shows that although $\ex A\fp fg \ex B$  as a set is equal to the fiber product of $\ex A$ with $\ex B$ as sets, $\ex A\fp fg \ex B$ as a topological space is not always the fiber product of $\ex A$ with $\ex B$ as topological spaces. In contrast, in the special case of $\mathbb Z$-transversality defined below, $\ex A\fp fg \ex B$ as a topological space is the fiber product of $\ex A$ with $\ex B$ as topological spaces

\begin{defn}[$\mathbb Z$-transverse]
Two smooth (or $C^{k,\delta}$) exploded maps 
\[\ex A\xrightarrow{f}\ex C\xleftarrow{g}\ex B\]
are $\mathbb Z$-transverse if they are transverse and if for every pair of points $p_1\longrightarrow\ex A$ and $p_2\longrightarrow \ex B$ so that $f(p_1)=g(p_2)$,  the $\mathbb Z$-linear span of the  image of $d(f)({}^{\mathbb Z}T_{p_1}(\ex A))$ with   $dg({}^{\mathbb Z}T_{p_{2}}( \ex B))$ is a lattice $L$ in ${}^{\mathbb Z}T_{f(p_1)}\ex C$ which is saturated in the sense that it is the intersection of a $\mathbb R$-linear subspace  with  ${}^{\mathbb Z}T_{f(p_1)}\ex C$. 
\end{defn}

Note that although transversality is generic, the property of $\mathbb Z$-transversality is not in general generic

\

\begin{lemma}\label{existence of fp}
If $f$ and $g$  transverse smooth or $C^{k,\delta}$ maps, 
\[\ex A\xrightarrow{f}\ex C\xleftarrow{g}\ex B\]
then the fiber product $\ex A\fp fg\ex B$ exists, and is smooth or $C^{k,\delta}$ respectively.
This fiber product shares the following properties with fiber products of smooth manifolds:

The dimension of $\ex A\fp fg\ex B$ is the sum of the dimensions of $\ex A$ and $\ex B$ minus the dimension of $\ex C$, and the map $\ex A\fp fg\ex B\longrightarrow \ex A\times \ex B$ is an injective map with an injective derivative.

Moreover, if  $f$ and $g$ are $\mathbb Z$-transverse, then $\ex A\fp fg\ex B$ has the topology of a subspace of $\ex A\times\ex B$, so  as a topological space $\ex A\fp fg\ex B$ is the fiber product of $\ex A$ and $\ex B$ as topological spaces. 

\end{lemma}

\pf

If we can construct fiber products locally in $\ex A$ and $\ex B$, then the universal property of fiber products will provide transition maps for the global construction of $\ex A\fp fg\ex B$. We may therefore restrict to the case that $\ex A$, $\ex B$ and $\ex C$ are standard coordinate charts. As the fiber product will be unaffected by considering the composition of $f$ and $g$ with an equidimensional submersion of $\ex C$, we may further assume that $\ex C=(0,\infty)^{n}\times \ex T^{m}$. We may then specialize further by noting the fiber product of the map $f/g:\ex A\times \ex B\longrightarrow (0,\infty)^{n}\times \ex T^{m}$ with the map sending a point to $(1,\dotsc,1,1\e 0,\dotsc, 1\e 0)$ is the same as the fiber product of $f$ with $g$. (Note that exchanging these two models does not affect the property of ($\mathbb Z$-)transversality). 

We have reduced to the case of considering the fiber product of a map $\pi:\ex X\longrightarrow \mathbb R^{n}\times \ex T^{m}$ with a single point. We may then reduce to the case where $\totb \pi$ is trivial as follows: Consider the subset of $\ex X$ with tropical part sent to the tropical part of our point. There is a bijective  equidimensional submersion of some $\mathbb R^{a}\times\et bP\times (\mathbb C^{*})^{c}$ onto this subset of $\ex X$. (Note that we can't say that our subset of $X$ is isomorphic to $\mathbb R^{a}\times \et bP\times (\mathbb C^{*})^{c}$ as the $(\mathbb C^{*})^{c}$ factor within our subset of $X$ has the trivial topology.)  As maps into this subset of $\ex X$ are equivalent to maps into $\mathbb R^{a}\times\et bP\times (\mathbb C^{*})^{c}$ followed by this map, we may replace $\ex X$ with this space for the purposes of calculating the fiber product. Then by choosing smaller open subsets of this space which are standard coordinate charts and relabeling, we have reduced to the case of a map $\pi:\ex X\longrightarrow \mathbb R^{n}\times\ex T^{m}$ so that the image of $\totb{\pi}$ is $0$. This is the point where $\mathbb Z$ transversality comes in. The topology on our fiber product will be the topology given by being a subset of  $\mathbb R^{a}\times\et bP\times (\mathbb C^{*})^{c}$. Any open subset of the original space $\ex X$ will correspond to an open subset of $\mathbb R^{a}\times\et bP\times (\mathbb C^{*})^{c}$ which is invariant in the $(\mathbb C^{*})^{c}$ direction. For the topology as a subset of our original space to match the topology as a subset of $\mathbb R^{a}\times\et bP\times (\mathbb C^{*})^{c}$, we require that each $(\mathbb C^{*})^{c}$ slice contains at most $1$ point of our fiber product. The map $\pi$ on each $(\mathbb C^{*})^{c}$ slice is given by monomials with exponents determined by the tropical part of the original map, and $\pi$ is injective on each slice if and only if our map is $\mathbb Z$ transverse.

We can now simplify  even further by restricting to the case where the target is $\mathbb R^{n}$. To see that this is no loss of generality, choose an injective equidimensional submersion of $\mathbb R^{n+2m}$ into the old target $\mathbb R^{n}\times \ex T^{m}$ with our point the image of $0$. Any map with image contained in the image of $\mathbb R^{n+2m}$ must then factor through this map $\mathbb R^{n+2m}\longrightarrow \mathbb R^{n}\times \ex T^{m}$. The fact that the tropical part of $\pi$ is constant implies that the $\pi^{-1}$ applied to the image of $\mathbb R^{n+2m}$ is an open neighborhood of $\pi^{-1}(0)$, so restricting to this open neighborhood, we may factor the map $\pi$ through this embedding.

 Finally, we have simplified the calculation of the fiber product to the local case of a the fiber product of a map $\pi:\ex X\longrightarrow \mathbb R^{n}$ with a map of a point to $0\in\mathbb R^{n}$. Transversality in this case corresponds to $d\pi$ being surjective. We may assume by restricting to an open subset of $\ex X$ if necessary that $\pi^{-1}(0)$ contains the largest stratum in $\totb{\ex X}$. Then, the fact that $d\pi$ is surjective and $\totb{\pi}$ is trivial imply that we can split the standard coordinate chart $\ex X$ as $ U\times \mathbb R^{n}$ where $\pi(u,0)=0$ for $u$ in the largest stratum in $\totb{U}$ and the derivative of $\pi$ in the $\mathbb R^{n}$ direction at these points is an isomorphism. Then if $\pi$ is smooth,  Lemma \ref{smooth intersect} implies that in an open neighborhood, $\pi^{-1}(0)$ is smooth, and if $\pi$ is $C^{k,\delta}$, Lemma \ref{ckdelta intersect} implies that in an open neighborhood,  $\pi^{-1}(0)$ is $C^{k,\delta}$. This smooth or $C^{k,\delta}$ exploded manifold $\pi^{-1}(0)$ has the required universal property that maps to $\pi^{-1}(0)$ are equivalent to maps to $U$ which when composed with $\pi$ give $0$. We have therefore constructed a local model for the fiber product. Note that the map of this model into $\ex A\times \ex B$ has an injective derivative, and this model has the expected dimension.
 
 The fiber product  $\ex A\fp fg \ex B$ can now be described as an exploded manifold as follows: As a set, $\ex A\fp fg \ex B$ is just the fiber product of $\ex A$ and $\ex B$ as sets. By a smooth or $C^{k,\delta}$ map to $\ex A\fp fg\ex B$, we shall mean a map $\ex F\longrightarrow \ex A\fp fg \ex B$ as sets which comes from smooth or $C^{k,\delta}$ maps $h_{\ex A}:\ex F\longrightarrow \ex A$ and $h_{\ex B}:\ex F\longrightarrow \ex B$ so that $f\circ h_{\ex A}=g\circ h_{\ex B}$. Given any point $p\in \ex A\fp fg\ex B$, the above local construction gives  some subset  $U\subset \ex A\fp fg \ex B$ containing $p$ the structure of an exploded manifold with the following two properties:
 
 \begin{enumerate}
 \item The inclusion $U\subset \ex A\fp fg \ex B$ is smooth or $C^{k,\delta}$ in the above sense,
 \item  Given any smooth or $C^{k,\delta}$ map $h:\ex F\longrightarrow \ex A\fp fg\ex B$, the subset $h^{-1}(U)\subset \ex F$ is open and the corresponding map $h^{-1}(U)\longrightarrow U$ is smooth or $C^{k,\delta}$ as a map of exploded manifolds.
\end{enumerate} 
  Repeating the construction around every point, we can describe $\ex A\fp fg\ex B$ as a topological space using the topology generated by the open subsets of these $U$. In the case of $\mathbb Z$-transversality, this agrees with the topology of $\ex A\fp fg\ex B$ as a subspace of $\ex A\times \ex B$, and therefore agrees with the fiber product of $\ex A$ with $\ex B$ as topological spaces. Given two subsets $U$ and $U'$ satisfying the above two properties,  the exploded structure on  $U\cap U'$ considered as an open subset of $U$ is the same as its exploded structure considered as an open subset of $U'$. Therefore, these $U$ give coordinate charts for $\ex A\fp fg \ex B$ as a smooth or $C^{k,\delta}$ exploded manifold.
  
  \stop    

The following lemma tells us that the tangent space to a fiber product acts in the same way as in the category of smooth manifolds. This implies that we may orient fiber products in the usual fashion. (See \cite{dre} for a discussion on orienting fiber products.) 

\begin{lemma}\label{fp sequence}
If $f:\ex A\longrightarrow \ex C$ is transverse to $g:\ex B\longrightarrow\ex C$, then for any point $(p_{1},p_{2})\in\ex A\fp fg \ex C$, the derivatives of the maps in the following commutative diagram
\[\begin{array}{ccc}\ex A\fp fg\ex B&\xrightarrow{\pi_{2}} &\ex B
\\ \downarrow \pi_{1} & & \downarrow g
\\ \ex A&\xrightarrow{f} &\ex C \end{array}\]
give a short exact sequence
\[0\longrightarrow T_{(p_{1},p_{2})}\lrb{\ex A\fp fg\ex B}\xrightarrow{(d\pi_{1},d\pi_{2})} T_{p_{1}}\ex A\times T_{p_{2}}\ex B\xrightarrow{df-dg} T_{f(p_{1})}\ex C\longrightarrow 0\]
\end{lemma}

\pf

The fact that $(d\pi_{1},d\pi_{2}):T_{(p_{1},p_{2})}\lrb{\ex A\fp fg\ex B}\longrightarrow T_{p_{1}}\ex A\times T_{p_{2}}\ex B$ is injective is part of Lemma \ref{existence of fp}. The fact that $(df-dg)\circ(d\pi_{1},d\pi_{2})=0$ follows from the fact that $f\circ \pi_{1}=g\circ \pi_{2}$. The surjectivity of $(df-dg)$ follows from transversality, and then the exactness of our sequence of maps follows from the fact proved in Lemma \ref{existence of fp} that the dimension of the middle term in the sequence is equal to the sum of the dimensions of the first and last terms.

\stop   

\section{Families and refinements}

\

\begin{defn}[Family]\label{family defn}
  A family of exploded manifolds over $\ex F$ is a map $f:\ex C\longrightarrow \ex F$ so that:
 \begin{enumerate}
  \item $f$ is complete
  \item for every point $p\in \ex C$, 
  \[df:T_p\ex C\longrightarrow T_{f(p)}\ex F\text{ is surjective }\]
 \[\text{ and }df:{}^{\mathbb Z}T_p\ex C\longrightarrow {}^{\mathbb Z}T_{f(p)}\ex F\text{ is surjective }\]
 \end{enumerate}

\end{defn}

(The definition of complete is found on page \pageref{complete}. Recall also that integer vectors in ${}^{\mathbb Z}T_p\ex C$ are the vectors $v$ so that for any exploded function, $vf$ is an integer times $f$. For example, on $\ex T$ the integer vectors are given by integer multiples of the real part of $\tilde z\frac \partial{\partial \tilde z}$, so the map $\ex T\longrightarrow\ex T$ given by $\tilde z^{2}$ is not a family as it does not obey the last condition above.)

Unlike in the smooth category, smooth families in the exploded category need not be locally in the form of a product. For example, there exists a smooth family of exploded manifolds which in some coordinate chart is given by a map $\et mP\longrightarrow\et {k}Q$ given by $(\tilde z,\tilde w)\mapsto \tilde z$, so long as the
the polytope $Q$ is given by the projection of $P$ to the first $k$ coordinates, the projection $P\longrightarrow Q$ is a complete map, the projection of each stratum of $P$ is a stratum of $Q$, and restricted to each stratum, the derivative of the projection map applied to integer vectors is surjective onto the lattice of integer vectors in the image stratum. This may differ from a product because  $P$ may not be a product of $Q$ with something. 

\includegraphics{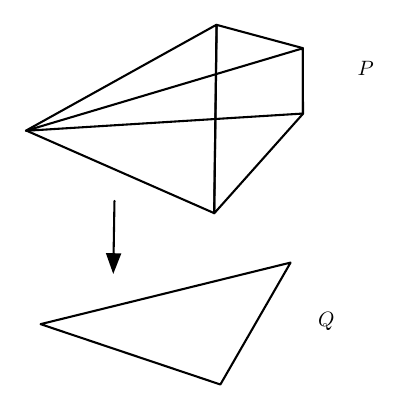}

\begin{example}[Moduli space of stable exploded curves.]
\end{example}
We can represent the usual compactified moduli space  of stable curves $\bar{\mathcal M}_{g,n}$ as a complex orbifold. There exist local holomorphic coordinates so that the boundary of $\bar{\mathcal M}_{g,n}$ in these coordinates looks like $\{z_{i}=0\}$. As in example \ref{expl} on page \pageref{expl}, we can replace these coordinates $z_{i}$ with $\tilde z_{i}$ to obtain a complex  exploded orbifold $\expl \left(\bar{\mathcal M}_{g,n}\right)$. The forgetful map $\pi:\bar{\mathcal M}_{g,n+1}\longrightarrow\bar{\mathcal M}_{g,n}$ induces a map 
\[\pi:\expl\left(\bar{\mathcal M}_{g,n+1}\right)\longrightarrow\expl\left(\bar{\mathcal M}_{g,n}\right)\]
This is a family, and each stable exploded curve with genus $g$ and $n$ marked points corresponds to the fiber over some point $p\longrightarrow\expl\left(\bar{\mathcal M}_{g,n}\right)$. Actually, it is proved in \cite{egw} that $\expl M_{g,n}$ represents the moduli stack of $\C\infty1$ stable curves.

\begin{example}[Model for node formation]\label{local family model}
\end{example}
The following example contains all interesting local behavior of the above example. It is not a family only because it fails to be complete (it is not  proper). Consider the map 
\[\pi:\et 22\longrightarrow\et 11\text{ given by }\pi^{*}\tilde z=\tilde w_{1}\tilde w_{2}\]
The derivative is surjective, as can be seen by the equation 
\[\pi^{*}(\tilde z^{-1}d\tilde z)=\tilde w_{1}^{-1}d\tilde w_{1}+\tilde w_{2}^{-1}d\tilde w_{2}\]

The fibers of this map over smooth points $\tilde z=c\e 0$ are smooth manifolds equal to $\mathbb C^{*}$ considered just as a smooth manifold with coordinates $w_{1}$ and $w_{2}\in\mathbb C^{*}$ related by $w_{1}w_{2}=c$. (Note that there is no point with $\tilde z=0\e0$.)

In contrast, the fibers of this map over points $\tilde z=c\e x$ with $x>0$ are isomorphic to $\et 1{[0,x]}$.



\begin{lemma} \label{base change} If $\pi:\ex C\longrightarrow \ex F$ is a family, and $f:\ex G\longrightarrow \ex F$ is any map, then $\pi':\ex C\fp \pi f\ex G\longrightarrow \ex G$ is a family.  

\end{lemma}

\pf The fact that $\pi$ is a family implies that $\pi$ and $f$ are $\mathbb Z$-transverse, so the fiber product exists, and has the same topology as the fiber product of $\ex C$ and $\ex G$ as topological spaces. Therefore $\pi'$ is proper.  

 We shall  now verify that $\pi'$ is complete:
 \[\begin{array}{ccccc}
 \et 1{(0,l)}&\xrightarrow {\gamma}&\ex C\fp\pi f \ex G&\longrightarrow &\ex C
\\ \downarrow & & \downarrow\pi' &&\downarrow \pi
\\ \et 1{[0,l]}&\xrightarrow {\gamma'}&\ex G&\xrightarrow {f}&\ex F
 \end{array}\]
  Suppose that $\gamma:\et1{(0,l)}\longrightarrow \ex C\fp \pi f\ex G$ is a map which when projected to $\ex G$ extends to a map of $\gamma':\et 1{[0,l]}\longrightarrow \ex G$. The composition of $\gamma'$ with $f$ is equal on $\et 1{(0,l)}$ to the projection  of $\gamma$ to $\ex C$ followed by $\pi$. As $\pi$ is a family, $f\circ \gamma'$ may be lifted to $\ex C$ on some neighborhood of $\et 1{(0,l)}\subset\et 1{[0,l]}$ to agree with  the projection of $\gamma$ to $\ex C$ on $\et 1{(0,l)}$. This then gives an extension of $\gamma$ to this neighborhood of  $\et 1{(0,l)}\subset\et 1{[0,l]}$. Then reparametrizing gives an extension of $\gamma$ to $\et 1{[0,1]}$. Therefore $\pi'$ is a complete map.  

The short exact sequence from Lemma \ref{fp sequence}
\[0\longrightarrow T_{(p_{1},p_{2})}\lrb{\ex C\fp \pi f\ex G}\xrightarrow{(d\pi_{\ex C},d\pi')} T_{p_{1}}\ex C\times T_{p_{2}}\ex G\xrightarrow{d\pi-df} T_{\pi(p_{1})}\ex F\longrightarrow 0\]
implies that if $d\pi$ is surjective, $d\pi'$ must also be surjective.

The last property that we need to check is that $d\pi'\lrb{{}^{\mathbb Z}T_{(p_{1},p_{2})}\ex C\fp \pi f\ex G}={}^{\mathbb Z}T_{p_{2}}\ex G$. This is equivalent to the requirement that given any map \[\gamma :\et 1{(-\epsilon,\epsilon)}\longrightarrow \ex G\] so that $\gamma(1\e0)=p_{2}$, for small enough $\epsilon'$, there exists a map \[\gamma':\et 1{(-\epsilon',\epsilon')}\longrightarrow\ex C\fp\pi f\ex G\] so that  $\pi'\circ\gamma'=\gamma$ and $\gamma'(1\e 0)=(p_{1},p_{2})$. Given such a $\gamma$, we may use the analogous property for the family $\pi$ on the map $f\circ \gamma:\et 1{(-\epsilon,\epsilon)}\longrightarrow \ex F$ to construct a map $\hat\gamma:\et 1{(-\epsilon',\epsilon')}\longrightarrow \ex C$ with the property that $\pi\circ \hat \gamma=f\circ\gamma$ and $\hat\gamma(1\e 0)=p_{2}$. Therefore $(\gamma,\hat \gamma)$ defines our map $\gamma':\et 1{(-\epsilon',\epsilon')}\longrightarrow\ex C\fp\pi f\ex G$ so that  $\pi'\circ\gamma'=\gamma$ and $\gamma'(1\e 0)=(p_{1},p_{2})$ as required. Therefore $\pi':\ex C\fp\pi f\ex G\longrightarrow \ex G$ is a family.

\stop

     \begin{defn}[Refinement]
      A refinement of $\ex B$ is an exploded manifold $\ex B'$ with a map $f:\ex B'\longrightarrow \ex B$  so that
      \begin{enumerate}
 \item $f$ is proper
 \item $f$ gives a bijection between points in $\ex B'$ and $\ex B$
 \item $df$ is surjective
\end{enumerate}

     \end{defn}

\begin{example}[Local model of refinement]\label{refinement example}\end{example}
All refinements are locally of the following form: Suppose that the polytope $P\subset\mathbb R^m$ is subdivided into a union of polytopes with nonempty interior $\{ P_i\}$ so that $P_{i}\subset P$ are closed and the intersection of $P_{i}$ and $P_{j}$ is a (possibly empty) face of $P_{i}$. Then using the standard coordinates $(x,\tilde z)$ on $\mathbb R^{n}\times \et m{P_{i}}$ for each $ P_i$, we can piece the coordinate charts $\mathbb R^n\times \et m{ P_i}$ together using the identity map as transition coordinates. In these coordinates, the map down to $\mathbb R^n\times\et mP$ is the identity map. 

To see that a given refinement $f$ is locally in the above form, restrict to a coordinate chart $\mathbb R^{n}\times\et mP$ on $\ex B$. The fact that $f$ is bijective and proper implies that the image of the tropical part of $f$ subdivides $P$ into polytopes $P_{i}$ which are closed sub polytopes of $P$. The inverse image of the set of points with tropical part $P_{i}$ is an open subset $U$ of $\ex B'$ with tropical part $P_{i}$, and the map $f$ restricted to this open subset can be regarded as a bijective submersion to $\mathbb R^{n}\times \et m{P_{i}}$ so that $\totb{f}:\totb U\longrightarrow P_{i}$ is an isomorphism. It follows that $f:U\longrightarrow \mathbb R^{n}\times\et m{P_{i}}$ is an isomorphism, and our refinement is locally of the type described above.

\

 Note that the above local model implies that if $\ex B$ is basic, a refinement is simply determined by a subdivision of the tropical part $\totb{\ex B}$. 

\

The effect of refinement on the smooth part $\totl{\ex B}$ should remind the reader of the correpsondence between toric blowups and subdivision of toric fans. For example, if $\ex B$ is the refinement of $\ex T^{m}$ given by subdividing $\totb{\ex T^{m}}$ into a toric fan, then $\totl{\ex B}$ is the corresponding toric manifold.

\begin{lemma}

Given any refinement $f:\ex B'\longrightarrow \ex B$, and a map $g:\ex C\longrightarrow \ex B$, the fiber product $\ex C':=\ex B'\fp fg\ex C$ is a refinement of $\ex C$.
\[\begin{array}{ccc}\ex C'&\xrightarrow{g'} &\ex B'
\\\downarrow & & \downarrow
\\ \ex C&\xrightarrow{g} &\ex B\end{array}
\]

\end{lemma}

\pf
As $f$ is a refinement, $f$ and $g$ are $\mathbb Z$-transverse, so the fiber product $\ex C'$ does exist, and has the topology of the fiber product of $\ex B'$ with $\ex C$ as topological spaces. It follows that as $f$ is proper, the corresponding map $\ex  C'\longrightarrow \ex C$ is proper. Similarly, as $f$ is a bijection, $\ex C'\longrightarrow\ex C$ is a bijection. Lastly, as argued in the proof of Lemma \ref{base change} above, the short exact sequence from Lemma \ref{fp sequence} gives that the  derivative of $\ex C'\longrightarrow \ex C$ is surjective because the $df$ is surjective.
\stop

\

The following lemma follows from the above standard local form for refinements. 

\begin{lemma}
 Given a refinement of $ \ex B$, 
 \[f:\ex B'\longrightarrow \ex B\]
 any smooth  vector field  $v$ on $\ex B$ lifts uniquely to a smooth vector field $\tilde v$ on $\ex B$ so that $df(\tilde v)=v$.
\end{lemma}

The above lemma tells us that any smooth tensor field (such as an almost complex structure or metric) lifts uniquely to a smooth tensor field on any refinement. 

\begin{defn}[Stable curve] Call a holomorphic curve stable if it has a finite number of automorphisms, and is not a nontrivial refinement of another holomorphic curve. Call a map of a holomorphic curve $f:\ex C\longrightarrow \ex B$ stable if it has a finite number of automorphisms, and it does not factor  as $\ex C\longrightarrow \ex C_{0}\longrightarrow \ex B$ where $\ex C\longrightarrow\ex C_{0}$ is a refinement map. 
\end{defn}

An equivalent definition of a stable curve is a curve $f:\ex C\longrightarrow \ex B$ which has smooth part $\totl{f}:\totl{\ex C}\longrightarrow \totl{\ex B}$ with a finite number of automorphisms. Note that that if a curve $f$ factors through a nontrivial refinement map, its smooth part will have an unstable spherical component with an infinite number of automorphisms which will not lift to automorphisms of $f$. 

If $\ex B$ has an almost complex structure,  there is a bijection between stable holomorphic curves in $\ex B$ and stable holomorphic curves in any refinement $\ex B'$, given by the fiber product of curves mapping to $\ex B$ with the refinement map $\ex B'\longrightarrow \ex B$. In fact, when the moduli space of stable holomorphic curves in $\ex B$ is smooth, the moduli space of stable holomorphic curves in $\ex B'$ is a refinement of the moduli space of curves in $\ex B$.

\

Recall that exploded tropical functions $\mathcal E(\ex B)$ are $\mathbb C\e{\mathbb R}$ valued functions which are locally a finite sum of exploded functions $\ex B\longrightarrow \ex T$.
\begin{thm}\label{zeroset}
If an exploded tropical function $f\in \mathcal E(\ex B)$ is transverse to $0\e{\mathbb R}$, then the subset 
\[f^{-1}\lrb{0\e{\mathbb R}}\subset\ex B\]
is a codimension $2$ exploded submanifold of $\ex B$.
\end{thm}

\pf
 Given an exploded tropical function $f\in \mathcal E(\ex B)$, we can construct a section $s_{f}$ of a $\mathbb C$ bundle over a refinement of $\ex B$ as follows:

\begin{enumerate}
\item On each coordinate chart $\mathbb R^{n}\times \et mP$,  the tropical part of $f$  is a convex piecewise linear function $\totb f$ on $P$. The regions of linearity $P_{i}\subset P$ of $\totb f$ determine a subdivision of $P$, which determines a refinement $\ex B'$ of $\ex B$. 
\item On each coordinate chart  $U$ of $\ex B'$, $\totb f$ is linear, so we may choose some monomial in coordinate functions $\tilde w=\e{a}\tilde z^{\alpha}$ which has tropical part equal to $\totb f$. To $(U,\tilde w)$, we may assign a coordinate chart $U\times \mathbb C$ on a $\mathbb C$ bundle over $\ex B'$ so that over $U$, $f$ corresponds to the section $s_{f}=f\tilde w^{-1}$ of $U\times \mathbb C$. Transition maps from $(U,\tilde w)$ and $(U',\tilde w')$ are given by the transition map from  $U$ to $U'$ combined with multiplication by $\tilde w\tilde w'^{-1}$ on the $\mathbb C$ factor. It follows that $s_{f}$ defines a global section of the line bundle with these coordinate charts.
\end{enumerate}

The exploded function $f$ being transverse to $0\e {\mathbb R}$ is equivalent to the section $s_{f}$ being transverse to the zero section, therefore if $f$ is transverse to $0\e{\mathbb R}$, the subset of $\ex B$ where $f\in 0\e{\mathbb R}$ is an exploded submanifold of $\ex B$ with codimension $2$.
  
  \stop

\begin{example}[Mikhalkin's pair of pants decomposition of a toric hypersurface as an exploded family]\label{pair of pants}\end{example}
Let \[p=\sum_{\alpha\in S} c_{\alpha}\tilde z^{\alpha}\] be a tropical exploded polynomial on $\ex T^{n}$, where $S$ indicates some finite set of exponents in $\mathbb Z^{n}$, and $c_{\alpha}\in\mathbb C^{*}$. For generic choice of $c_{\alpha}$, $p$ will be transverse to $0\e{\mathbb R}$, so 
\[Z_{p}:=p^{-1}\lrb{0\e{\mathbb R}}\subset \ex T^{n}\]
will be an exploded manifold. The tropical part of $Z_{p}$ is equal to the set where $\totb p$ is not smooth, which corresponds to the subset of points  in $ \totb{\ex T^{n}}$ where $\totb{\tilde z^{\alpha}}=\totb{\tilde z^{\alpha'}}$ for some $\alpha\neq \alpha'$ in $S$. The tropical part of $Z_{p}$ can also be determined by taking all polytopes of dimension $<n$ in the dual fan to the Newton polytope of $p$, which is the convex hull of $S$.

 Let $M$ be the toric space with toric fan dual to the Newton polytope of $S$, and suppose that $M$ is a manifold. Another way of viewing $Z_{p}$ is as follows: The zero set $Z'$ of $\sum_{\alpha\in S}c_{\alpha}z^{\alpha}$ defines a complex submanifold of $M$ which intersects the toric boundary divisors of $M$ nicely because the Newton polytope of $p$ is dual to the toric fan of $M$. Therefore we can regard $Z'$ as a complex manifold with normal crossing divisors given by the toric boundary divisors of $M$. The explosion of $Z'$ is equal to $Z_{p}$. The map $\expl{Z'}\longrightarrow \expl M\longrightarrow \ex T^{n}$ is an isomorphism onto $Z_{p}$.

Suppose now that $S$ contains every lattice point in its convex hull. We shall construct a family of exploded manifolds which corresponds to Mikhalkin's  higher dimensional pair of pants decomposition of the toric hypersurface $Z'$ from \cite{mpants}. This requires choosing a convex function $v:S\longrightarrow\mathbb Z$ so that the convex hull of the set of points over the graph of $v$ in $\mathbb R^{n}\times \mathbb R$ has faces which project to standard simplices in $\mathbb R^{n}$ with volume $\frac 1{n!}$. Given such a function, we can consider the polynomial 

\[\hat p=\sum_{\alpha\in S}c_{\alpha}\tilde w^{v(\alpha)}\tilde z^{\alpha}\ \ \ \in \mathcal E(\et 11\times \ex T^{n})\]

We shall see that the set $Z_{\hat p}$ of noninvertibility of $\hat p$ restricted to a subset where $\tilde w$ is small enough will be an exploded manifold. Our family shall be given by the projection of $Z_{\hat p}$ to (an open subset of) $\et 11$ given by the coordinate $\tilde w$. The set $Z_{\hat p}$ has as its tropical part the points in $\totb{\et 11\times \ex T^{n}}$ where $\totb{\tilde w^{v(\alpha)}\tilde z^{\alpha}}=\totb{\tilde w^{v(\alpha')}\tilde z^{\alpha'}}$ for some $\alpha\neq \alpha'$ in $S$. As before, this corresponds to the polytopes of dimension $\leq n$ in the dual fan to the Newton polytope of $\hat p$ although in this case we have to intersect this dual fan with the half space $\totb{\et 11\times \ex T^{n}}$ - the upshot of this is that we get the polytopes of dimension $\leq n$ in the dual fan to the convex hull of the set of points over the graph of $v$. The edges of $\totb{Z_{\hat p}}$ which are in the interior of $\totb{\et 11\times \ex T^{n}}$ are in correspondence with the downward pointing faces of the Newton polytope of $\hat p$. These downward pointing faces project to standard lattice simplices in $\mathbb Z^{n}$ with corners $\alpha_{0}\dotsc \alpha_{n}$ which an invertible affine transformation of $\mathbb Z^{n}$ can transform to $0$ and the standard basis vectors. Over such an edge, $\hat p$ reduces to $\sum_{i=0}^{n}c_{\alpha_{i}}\tilde w^{v(\alpha_{i})}\tilde z^{\alpha_{i}}$. Similarly, restricted to the interior of $\totb{\et 11\times \ex T^{n}}$,  over the interior of  any $k$ dimensional face of $\totb{Z_{\hat p}}$ the polynomial  $\hat p$ reduces to $\sum_{i=0}^{n+1-k}c_{\alpha_{i}}\tilde w^{v(\alpha_{i})}\tilde z^{\alpha_{i}}$ where the $\alpha_{i}$ can be transformed using an invertible  affine transformation of $\mathbb Z^{n}$ to $0$ and the first $(n+1-k)$ standard basis vectors. 

It follows that when $\tilde w$ is sufficiently small, $\hat p$ will be transverse to $0$, and the map $\tilde w:Z_{\hat p}\longrightarrow \et 11$ will be a family of exploded manifolds because $\hat p$ restricted to the slices where $\tilde w$ is constant is transverse to $0\e{\mathbb R}$. As this transversality fails in real codimension $2$ and $\tilde w$ is also a family around $\tilde w=1\e 0$ we may choose a connected open subset $\ex F\subset \et 11$ containing both $1\e 0$ and $1\e 1$ so $\tilde w:Z_{\hat p}\longrightarrow \et 11$ restricted to $\ex F$ is a connected holomorphic  family of exploded manifolds. The inverse image of $1\e 0$ is our original exploded manifold $Z_{p}$ obtained by exploding a toric hypersurface. The inverse image of $1\e 1$ corresponds to Mikhalkin's pair of pants decomposition of this hypersurface. Because of the standard form of $\hat p$ mentioned in the previous paragraph, each maximal dimensional stratum of the smooth part of $\tilde w^{-1}(1\e{1})$ can be identified with  the subset $\{z_{1}+z_{2}\dotsc +z_{n}=1\}\subset (\mathbb C^{*})^{n}$, which is equal to $\mathbb CP^{n-1}$ minus $n+1$ generic hyperplanes - in other words a higher dimensional pair of pants. Moreover, the closure of each maximal dimensional stratum is equal to $\mathbb CP^{n-1}$ where the lower dimensional strata are given by these hyperplanes and their intersections.

\section{Moduli stack of  exploded  curves}\label{perturbation theory}

We shall use the concept of a stack without giving the general definition. (See the article \cite{stacks} for a readable introduction to stacks). The reader unfamiliar with stacks may just think of our use of stacks as a natural way of encoding information about families of holomorphic curves.

 When we say that we shall consider an exploded manifold $\ex B$ as a stack, we mean that we replace $\ex B$ with a category $\St{\ex B}$  over the category of exploded manifolds (in other words a category $\St{\ex B}$ with a functor to the category of exploded manifolds) as follows: 
 objects in $\St{\ex B}$ are maps into $\ex B$:
\[\ex A\longrightarrow\ex B\]
and morphisms are commutative diagrams
\[\begin{array}{lll}\ex A & \longrightarrow &\ex B
\\ \downarrow & &\downarrow\id
\\ \ex C &\longrightarrow &\ex B
\end{array}\]
The functor from $\St{\ex B}$ to the category of exploded manifolds is given by sending $\ex A\longrightarrow\ex B$ to $\ex A$, and the above morphism to $\ex A\longrightarrow\ex C$.

Note that  maps $\ex B\longrightarrow\ex D$ are equivalent to  functors  $\St{\ex B}\longrightarrow\St{\ex D}$ which commute with the functor down to the category of exploded manifolds. Such functors are  morphism of categories over the category of exploded manifolds, and this is the correct notion of maps of stacks. 

For example, a point thought of as a stack is equal to the category of exploded manifolds itself, with the functor down to the category of exploded manifolds the identity. We shall refer to points thought of this way simply as points.

We shall define the moduli stack of $\C\infty1$ curves below. The regularity $\C\infty1$ is used because that is the natural regularity that the moduli space of holomorphic curves has in the case of transversality.  Similar definitions can be made using smooth instead of $\C\infty1$. In what follows, all morphisms will be assumed $\C\infty1$. In particular, $\St{\ex B}$ will refer to $\ex B$ considered as a stack using the category of $\C\infty1$ exploded manifolds. 

\begin{defn}[Moduli stack of $\C\infty1$ curves]
 The moduli stack $\smod(\ex B)$ of  $\C\infty1$ exploded curves in $\ex B$ is a category over the category of $\C\infty 1$ exploded manifolds with objects being families of $\C\infty1$ exploded curves consisting of the following
  
  \begin{enumerate}
 \item 
 A $\C\infty1$ exploded  manifold $\ex C$
 \item A pair of $\C\infty1$ exploded  morphisms
 
 \[\begin{split}
    &\ex C\longrightarrow \ex B
    \\ \pi &\downarrow
    \\&\ex F
   \end{split}
\]
\item A  $\C\infty1$ section $j$ of $\ker(d\pi)\otimes\left(T^*\ex C/\pi^*(T^*\ex F)\right)$ 
\end{enumerate}
so that
\begin{enumerate}
\item $\pi:\ex C\longrightarrow \ex F$ is a family (definition \ref{family defn} on page \pageref{family defn}).

\item The inverse image of any point $p\in \ex F$ is an exploded  curve with complex structure $j$.

\end{enumerate}

\

A morphism between families of curves is given by $\C\infty1$  morphisms $f$ and $c$ making the following diagram commute

\begin{displaymath}
\begin{array}{lllll}
\ex F_1 & \longleftarrow & \ex C_1 & \longrightarrow & \ex B\\ 
\downarrow f &  & \downarrow c & × & \downarrow\id \\ 
\ex F_2 & \longleftarrow & \ex C_2  & \longrightarrow & \ex B
\end{array}
\end{displaymath}
 
so that $c$ is a $j$ preserving isomorphism on fibers. 

The functor down to the category of $\C\infty1$ exploded manifolds is given by taking the base $\ex F$ of a family.
\end{defn}
Note that morphisms are not quite determined by the map  $f:\ex F_{1}\longrightarrow\ex  F_{2}$.  $\ex C_{1}$ is non-canonically isomorphic to the fiber product of $\ex C_{2}$ and $\ex F_{1}$ over $\ex F_{2}$.

\

This is a moduli stack in the sense that a  morphism $\St{\ex F}\longrightarrow \smod(\ex B)$ is equivalent to a  $\C\infty1$ family  curves  $\ex F\longleftarrow\ex C\longrightarrow\ex B$ (this is the family which is the image of the identity map $\ex F\longrightarrow \ex F$ considered as an object in $\St{\ex F}$).

Recall that
a holomorphic curve $\ex C\longrightarrow\ex B$ is stable if it has a finite number of automorphisms, and is not a nontrivial refinement of another holomorphic curve. (If $\ex B$ is basic, this is equivalent to all smooth components of $\ex C$ which are mapped to a point in $\totl{\ex B}$ being  stable as punctured Riemann surfaces.)

\begin{defn}[Moduli stack of stable holomorphic curves]
Given an almost complex structure $J$ on $\ex B$, a $\C\infty1$ family of stable holomorphic curves in $\ex B$ is a $\C\infty1$ family of  curves so that the map restricted to fibers is holomorphic and stable. The moduli stack of stable holomorphic curves in $\ex B$ is the substack  $\mathcal M(\ex B)\subset\smod(\ex B)$ with objects consisting of all families of stable holomorphic curves, and morphisms the same as in $\smod$.
\end{defn}

%

%

 
 It is useful to be able to make  statements about holomorphic curves in families, so we generalize the above definitions for a family $\hat{\ex B}\longrightarrow\ex G$ as follows: 
 
 \begin{defn}[Moduli stack of curves in a family]
 The moduli stack of $\C\infty1$ curves in a family $\hat{\ex B}\longrightarrow\ex G$, $\smod(\hat{\ex B}\rightarrow\ex G)$ is the substack of $\smod(\hat{\ex B})$ which is the full subcategory which has as objects families which admit commutative diagrams
 \[\begin{array}{cll}(\ex C,j)&\longrightarrow&\hat{\ex B}
 \\ \downarrow& &\downarrow
 \\ \ex F&\longrightarrow & \ex G
 \end{array}\] 
 \end{defn}

The moduli stack of stable holomorphic curves $\mathcal M(\hat{\ex B}\rightarrow\ex G)$ is then defined as the appropriate substack  of $\smod(\hat{\ex B}\rightarrow\ex G)$. Note that there is a morphism $\smod(\ex B\rightarrow\ex G)\longrightarrow\St{\ex G}$ which sends the object given by the diagram above to $\ex F\longrightarrow\ex G$. The appropriate compactness theorem for families states that if we restrict to the part of the moduli space with appropriately bounded combinatorial and topological data,  the map $\mathcal M(\hat {\ex B}\rightarrow\ex G)\longrightarrow\St{\ex G}$ is  proper.

\

 We can put a topology on $\St{\ex B}$  and  $\smod $  as follows: Consider the set of points in $\smod (\ex B)$, or (isomorphism classes of) maps from a point considered as a stack to $\smod (\ex B)$, 
 
 \[\St{p}\longrightarrow \smod (\ex B)\]
 A map  such as the one above is equivalent to a single $\C\infty1$ curve in $\ex B$ which is the image of $p$ considered as an object in $\St p$. Therefore the `set of points' in $\smod (\ex B)$ corresponds to the set of isomorphism classes of $\C\infty 1$ curves in $\ex B$. 
 
 The set of points in $\St {\ex B}$ is equal to the set of points in $\ex B$, so we can give the set of points in $\St{\ex B}$ the same topology as $\ex B$.  Below, we define a topology on this set of points in $\smod(\ex B)$ by defining convergence of a sequence of points. This topology will be non-Hausdorff in the same way as the topology on $\ex B$ is non-Hausdorff. 

\

There are two notions of convergence that a topology on $\smod(\ex B)$ should take account of:
\begin{enumerate}
\item $\C\infty1$ convergence of a sequence of of maps from a fixed domain, 

\item and convergence of fibers within a finite dimensional $\C\infty1$ family of maps.
\end{enumerate}

It is easy to construct an example of a sequence of curves converging in either one of the above senses, but not converging in the other sense. It is natural to do analysis with the first notion of convergence, and the second notion of convergence is a very natural notion of convergence within a moduli stack which records the structure of finite dimensional $\C\infty1$ families of maps. The following definition of convergence mixes the above two notions, (so it is weaker than both of them).

\begin{defn}[$\C\infty1$ convergence of a sequence of curves]
A sequence of $\C\infty1$ curves  \[f^{i}:\ex C^{i}\longrightarrow \ex B\] converges to a given curve  \[ f:\ex C\longrightarrow\ex B\]  if there exists a sequence of $\C\infty 1$ families 
  
  \[\ex F\longleftarrow (\ex {\hat C},j_i)\xrightarrow{\hat f^i}\ex B\]

and a sequence of points $p^{i}$ in $\ex F$ so that

\begin{enumerate}\item  
  this sequence of families converges in $\C\infty 1$ as a sequence of maps  to 
  
  \[\ex F\longleftarrow (\ex {\hat C},j)\xrightarrow{\hat f}\ex B\]
   
 \item $p^i\rightarrow p$ in $\ex F$
\item  
  $f^{i}$ is the map given by the restriction of $\hat f^{i}$ to the fiber over $p^{i}$, and $f$ is given by the restriction of $f$ to the fiber over $p$. 
  
  \end{enumerate}
  
  \
  
  A sequence of points converge in $\smod$ if the corresponding sequence of $\C\infty1$ curves converges.

\end{defn}

\

We say that a sequence of points in $\smod(\hat{\ex B}\rightarrow\ex G)$ converge  if they converge in $\smod(\hat{\ex B})$. Similarly, define the topology on the set of points in $\mathcal M\subset\smod$ to be the subspace topology, so a sequence of holomorphic curves converge  if they converge as $\C\infty1$ curves.  

\

For compactness results about $\mathcal M(\ex B)$, see \cite{cem}. For further discussion of the structure of $\mathcal M(\ex B)$ and $\mathcal M^{\infty,\underline 1}(\ex B)$, see \cite{egw}.

\bibliographystyle{plain}
\bibliography{ref}
\end{document}